\documentclass{amsart}
\pdfoutput=1

\usepackage{amssymb}
\usepackage{amsthm}
\usepackage{amsfonts}
\usepackage{amsmath}
\usepackage{pmboxdraw}
\usepackage{verbatim} 
\usepackage{graphicx}
\usepackage{color}
\usepackage[colorlinks=true, citecolor=blue, filecolor=black, linkcolor=black, urlcolor=black]{hyperref}
\usepackage{cite}
\usepackage[normalem]{ulem}
\usepackage{subcaption}
\usepackage{bbm}
\usepackage{bm}
\usepackage{mathtools}
\usepackage{todonotes}
\usepackage{kantlipsum}
\allowdisplaybreaks

\usepackage{tikz}
\usepackage{xcolor}
\usepackage{graphicx}
\usetikzlibrary{arrows.meta,calc,positioning,fit,backgrounds,shadows.blur}

\usepackage{colortbl}
\usetikzlibrary{decorations, arrows.meta, patterns}
\usetikzlibrary{decorations.pathreplacing}

\usepackage{tabularray}
\UseTblrLibrary{booktabs}

\newcommand{\RN}[1]{%
	\textup{\uppercase\expandafter{\romannumeral#1}}%
}

\def\R{\mathbb{R}}

\def\q{\boldsymbol{q}}

\newcommand{\Geom}{\operatorname{Geom}}
\newcommand{\Exp}{\operatorname{Exp}}

\newcommand{\Ai}{\operatorname{Ai}}

\newcommand{\re}{\operatorname{Re}}
\newcommand{\im}{\operatorname{Im}}

\def\Re{\operatorname{Re}}
\def\Im{\operatorname{Im}}

\newcommand{\supp}{\operatorname{supp}}

\newcommand{\qbinom}[2]{{\begin{bmatrix}#1\\#2\end{bmatrix}}_q}
\newcommand{\Li}{\operatorname{Li}_{2}}
\newcommand{\Id}{\operatorname{Id}}

\theoremstyle{plain}
\newtheorem{thm}{Theorem}[section]

\newtheorem{cor}[thm]{Corollary}
\newtheorem{lem}[thm]{Lemma}

\newtheorem{prop}[thm]{Proposition}

\theoremstyle{remark}

\newtheorem{rem}{Remark}

\numberwithin{equation}{section}

\begin{document}

\title[Inhomogeneous LPP and the $q$-JUE]{Inhomogeneous Geometric Last Passage Percolation \\ and the $q$-deformed Jacobi Unitary Ensemble}
\author{Sung-Soo Byun}
\address{Department of Mathematical Sciences and Research Institute of Mathematics, Seoul National University, Seoul 08826, Republic of Korea;
School of Mathematics, Korea Institute for Advanced Study, 85 Hoegi-ro, Dongdaemun-gu, Seoul 02455, Republic of Korea}
\email{sungsoobyun@snu.ac.kr}

\author{Yeong-Gwang Jung}
\address{Department of Mathematical Sciences, Seoul National University, Seoul 08826, Republic of Korea}
\email{wollow21@snu.ac.kr} 

\author{Nick Simm}
\address{Department of Mathematics, University of Sussex, Brighton, BN1 9RH, United Kingdom}
\email{n.j.simm@sussex.ac.uk}

\begin{abstract}
We study directed last passage percolation with geometric weights in a spatially inhomogeneous environment. We derive an explicit limit shape and establish fluctuation limits for the passage time, showing that, depending on the strength of the inhomogeneity, the limiting fluctuations are governed by Fredholm determinants associated with either the Airy or Bessel kernel.   Our approach is based on an exact relation between the passage time and the smallest particle of the $q$-deformed Jacobi unitary ensemble ($q$-JUE). As part of our analysis, we derive the macroscopic limiting density of the $q$-JUE, which provides a natural $q$-deformation of the classical Wachter distribution, as well as of the limiting distribution of the Meixner ensemble. We further establish Airy and Bessel scaling limits for the smallest particle, which are of independent interest. 
\end{abstract}

\maketitle

\section{Introduction} 

Random growth models and random matrices provide two of the most celebrated examples of universal phenomena in modern probability theory. Over the past three decades, a sequence of remarkable discoveries has uncovered intimate connections between these seemingly unrelated subjects; see e.g. \cite{BDJ99,Ok00,BOO00,Joh00,TW01,Joh01,PS02} for some early work. These developments have shaped many of the central ideas of integrable probability and have led to a detailed understanding of a broad class of exactly solvable stochastic systems.

A landmark result in this direction is due to Johansson \cite{Joh00}, who showed that the passage time in last passage percolation (LPP) with i.i.d. geometric weights has the same distribution as the largest particle in the Meixner orthogonal polynomial ensemble.
This correspondence transformed a problem in a stochastic growth model into one of random matrix theory, allowing the powerful asymptotic machinery of orthogonal polynomial ensembles to be applied to the analysis of LPP. As a consequence, it led to explicit limit shapes, Tracy--Widom fluctuations, and a detailed understanding of homogeneous LPP.
More broadly, it established random matrix theory as one of the central tools of integrable probability.   

Since then, many stochastic growth models, interacting particle systems and related integrable systems have been shown to possess analogous algebraic representations through symmetric functions and algebraic combinatorics, leading to determinantal or Pfaffian point processes amenable to asymptotic analysis; see e.g. \cite{OR03,BR05,ACQ11,OCo12,BC14,BCS14,BP18,MQR21,DOV22,DLM23,DLM24}, as well as the surveys \cite{Cor12,BG16,BP14,Rem22}. More recently, these ideas have culminated in universal scaling objects such as the KPZ fixed point and the directed landscape \cite{MQR21,DOV22,ACH24a,ACH24b,DV24}. A common feature underlying these developments is the existence of an exactly solvable algebraic structure that serves as a dual description of the original stochastic model.
 
Despite these remarkable advances, explicit dualities of this kind remain largely confined to homogeneous systems. In many applications, however, one naturally encounters growth in heterogeneous media, where the local growth mechanism varies across space. Such spatial inhomogeneity can drastically alter both the macroscopic geometry and the fluctuation behaviour of the model, while simultaneously obscuring the algebraic structures responsible for exact solvability. Consequently, even determining the limiting shape often becomes a challenging problem.

Motivated by these challenges, various inhomogeneous versions of directed LPP have been investigated. Integrable examples with parameter-dependent geometric weights were studied in \cite{Joh08,FR05,BBP05}, while more general inhomogeneous environments, including macroscopically varying and discontinuous media, have been studied through hydrodynamic and variational methods \cite{RT08,Em16,Cald15,CG21,BCD21,EJS21, JR22}. Related developments include \cite{BEMSS26,KPS19,BGS26}. Nevertheless, beyond a small number of special cases, explicit asymptotic descriptions of both the limit shape and the fluctuations remain remarkably scarce.

\medskip 

From the perspective of the general philosophy of integrable probability discussed above, this naturally leads to the following question. 
\begin{center}
\emph{
What is the appropriate dual model for an inhomogeneous LPP?
}
\end{center}
Answering this question requires considerably more than finding another random matrix analogue. The desired dual model should retain enough algebraic structure to admit an exact finite-dimensional correspondence while remaining sufficiently tractable to allow a complete asymptotic analysis.   

The principal result of the present paper is that this role is played by
the family of $q$-deformed Jacobi unitary ensembles ($q$-JUE), a model in the class of $q$-deformed random matrix models that have recently been actively studied; see e.g. \cite{BFO26,BJM26,BF25a,FLSY23,FL20,Fo21,Ol21,JS26} and references therein. Specifically, we prove
that the passage time in the inhomogeneous geometric LPP \cite{Joh08,FR05,BBP05} is related, through an explicit logarithmic transformation, to the \emph{smallest particle} of a suitably parameterised $q$-JUE; see Proposition~\ref{Prop_duality Gnm q-JUE}. 

This exact correspondence allows the asymptotic theories of the two models to be developed in parallel. On the LPP side, it yields explicit limit shapes together with a complete description of the passage-time fluctuations. On the random matrix side, it leads to a detailed asymptotic theory of the $q$-JUE, including both its macroscopic spectral density and its edge statistics. Moreover, the strength of the spatial inhomogeneity is reflected directly in the spectral behaviour of the ensemble: weak inhomogeneity corresponds to a soft-edge regime with Airy asymptotics, whereas sufficiently strong inhomogeneity gives rise to a hard-edge regime with Bessel asymptotics. 
The present correspondence naturally extends the classical homogeneous picture. As $q\to1$, corresponding to the disappearance of spatial inhomogeneity, the $q$-JUE degenerates, under an appropriate scaling, to the Meixner ensemble appearing in Johansson's correspondence for homogeneous geometric LPP \cite{Joh00}.

\smallskip 

The rest of this paper is organised as follows.
Section~\ref{Section_LPP} introduces the inhomogeneous geometric LPP model and presents our main results on its law of large numbers and fluctuation behaviour.
In Section~\ref{Section_qJUE}, we formulate the corresponding results for the $q$-JUE.
The exact duality between the two models is established in Section~\ref{Section_duality}, and the remaining sections are devoted to the proofs.

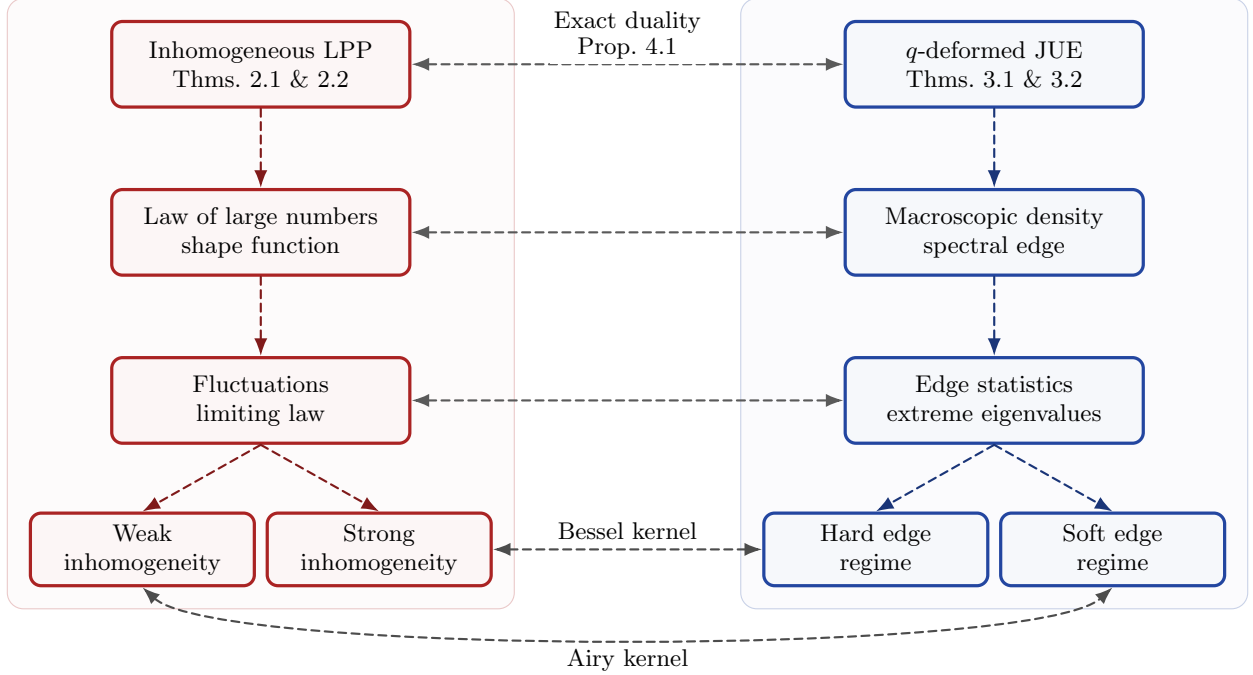
\begin{figure}[t]
\centering

\definecolor{lppcol}{RGB}{170,35,35}
\definecolor{juecol}{RGB}{35,70,160}
\definecolor{bridgecol}{RGB}{90,90,90}
\definecolor{airycol}{RGB}{90,90,90}
\definecolor{besselcol}{RGB}{90,90,90}

\resizebox{0.98\textwidth}{!}{
\begin{tikzpicture}[
    line cap=round,
    line join=round,
    >=Latex,
    font=\small,
    base/.style={
        rectangle,
        rounded corners=4pt,
        very thick,
        align=center,
        inner xsep=5pt,
        inner ysep=1pt,
        text width=3.55cm,
        minimum height=1.12cm
    },
    lppbox/.style={
        base,
        draw=lppcol,
        fill=lppcol!4
    },
    juebox/.style={
        base,
        draw=juecol,
        fill=juecol!4
    },
    smallbranchbox/.style={
        base,
        text width=2.65cm,
        minimum height=0.95cm,
        inner xsep=4pt,
        inner ysep=1pt
    },
    maparrow/.style={
        <->,
        thick,
        densely dashed,
        draw=bridgecol
    },
    imply/.style={
        ->,
        thick,
        draw=bridgecol
    },
    sidearrowL/.style={
        ->,
        thick,
        densely dashed,
        draw=lppcol!75!black
    },
    sidearrowR/.style={
        ->,
        thick,
        densely dashed,
        draw=juecol!75!black
    },
    airyarrow/.style={
        <->,
        thick,
        densely dashed,
        draw=airycol!85!black
    },
    besselarrow/.style={
        <->,
        thick,
        densely dashed,
        draw=besselcol!85!black
    }
]

\def\xL{-4.80}
\def\xR{4.80}

\def\ytop{0}
\def\ymidA{-2.20}
\def\ymidB{-4.40}
\def\yreg{-6.35}
\def\yairycorr{-7.75}

\def\xWeak{-6.35}
\def\xStrong{-3.25}
\def\xHard{3.25}
\def\xSoft{6.35}

\node[lppbox] (lpp) at (\xL,\ytop) {
    Inhomogeneous LPP\\
    Thms.~\ref{Thm_LPP geometric weak} \& \ref{Thm_LPP geometric strong}
};

\node[juebox] (jue) at (\xR,\ytop) {
    $q$-deformed JUE\\
    Thms.~\ref{thm:limiting density} \& \ref{Thm_qJUE scaling limits}
};

\draw[maparrow]
    (lpp.east) --
    node[
        above,
        font=\small,
        text=black,
        fill=white,
        inner sep=1.5pt,
        align=center
    ] {Exact duality\\[-1pt] Prop.~\ref{Prop_duality Gnm q-JUE}}
    (jue.west);

\node[lppbox] (lln) at (\xL,\ymidA) {
    Law of large numbers\\
    shape function
};

\node[juebox] (macro) at (\xR,\ymidA) {
    Macroscopic density\\
    spectral edge
};

\draw[sidearrowL] (lpp.south) -- (lln.north);
\draw[sidearrowR] (jue.south) -- (macro.north);

\draw[maparrow]
    (lln.east) --
    (macro.west);

\node[lppbox] (fluct) at (\xL,\ymidB) {
    Fluctuations\\
    limiting law
};

\node[juebox] (edge) at (\xR,\ymidB) {
    Edge statistics\\
    extreme eigenvalues
};

\draw[sidearrowL] (lln.south) -- (fluct.north);
\draw[sidearrowR] (macro.south) -- (edge.north);

\draw[maparrow]
    (fluct.east) --
    (edge.west);

\node[
    lppbox,
    smallbranchbox
] (weak) at (\xWeak,\yreg) {
    Weak\\
    inhomogeneity
};

\node[
    lppbox,
    smallbranchbox
] (strong) at (\xStrong,\yreg) {
    Strong\\
    inhomogeneity
};

\node[
    juebox,
    smallbranchbox
] (hard) at (\xHard,\yreg) {
    Hard edge \\ regime 
};

\node[
    juebox,
    smallbranchbox
] (soft) at (\xSoft,\yreg) {
    Soft edge \\ regime 
};

\draw[sidearrowL] (fluct.south) -- (weak.north);
\draw[sidearrowL] (fluct.south) -- (strong.north);

\draw[sidearrowR] (edge.south) -- (hard.north);
\draw[sidearrowR] (edge.south) -- (soft.north);

\draw[besselarrow]
    (strong.east) --
    node[
        above,
        yshift=2pt,
        font=\small,
        text=black,
        fill=white,
        inner sep=1.5pt
    ] {Bessel kernel}
    (hard.west);

\draw[airyarrow]
    (weak.south)
    .. controls (-5.35,\yairycorr) and (5.35,\yairycorr) ..
    node[
        pos=0.5,
        below,
        yshift=-2pt,
        font=\small,
        text=black,
        fill=white,
        inner sep=1.5pt
    ] {Airy kernel}
    (soft.south);

\begin{scope}[on background layer]
\node[
    draw=lppcol!25,
    fill=lppcol!2,
    rounded corners=6pt,
    inner sep=8pt,
    fit=(lpp)(lln)(fluct)(weak)(strong),
] {};

\node[
    draw=juecol!25,
    fill=juecol!2,
    rounded corners=6pt,
    inner sep=8pt,
    fit=(jue)(macro)(edge)(hard)(soft),
] {};
\end{scope}

\end{tikzpicture}
}

\caption{
Schematic overview of the duality framework and the corresponding asymptotic results for the inhomogeneous geometric LPP model and the $q$-JUE.
} 
\label{fig:duality-roadmap}
\end{figure}

Figure~\ref{fig:duality-roadmap} provides a schematic overview of the main results and the underlying duality framework of the present paper.

\section{Last passage percolation with inhomogeneous environment} \label{Section_LPP}

In this section, we introduce the spatially inhomogeneous LPP models \cite{Joh08,FR05,BBP05} and state our main results. We begin with
the geometric model, followed by its exponential analogue, since the
latter arises as an extremal limit of the former.

\subsection{Geometric LPP in the weak and strong inhomogeneity regimes} 

To incorporate spatial inhomogeneity, let $x_i,y_j\in[0,1)$ and assign to each lattice site $(i,j)$ an independent geometric random variable
\begin{equation}
\omega_{ij}\sim \Geom(x_i y_j),
\end{equation}
whose law is
\begin{equation}\label{def of omega geometric xi yj}
\mathbb{P}(\omega_{ij}=k)
=(1-x_i y_j)(x_i y_j)^k,
\qquad
k=0,1,2,\ldots.
\end{equation}
Throughout this paper, we consider the exponentially varying parameters
\begin{equation}\label{def of xi yj weights}
x_i=\alpha q^{\,i-1},
\qquad
y_j=\alpha q^{\,j-1},
\end{equation}
where $\alpha\in(0,1)$ and $q\in(0,1]$. 
When $q=1$, this reduces to the classical homogeneous LPP model \cite{Joh00,Sep98}. We note that the parameter $q$ will later be identified with the deformation parameter of the $q$-JUE.  

\begin{figure}[t]
	\begin{subfigure}{0.32\textwidth}
		\begin{center}	
			\includegraphics[width=0.9\textwidth]{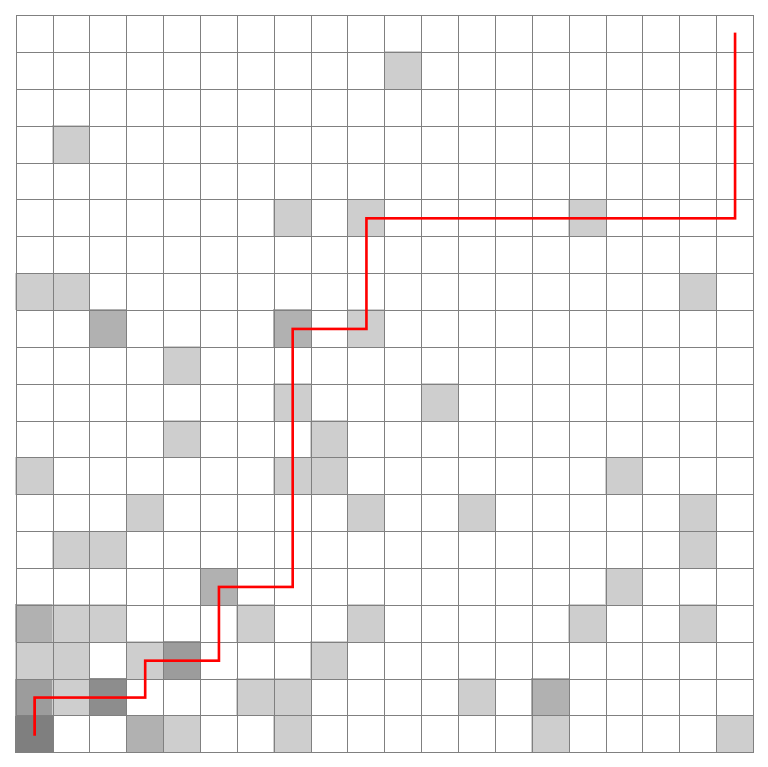}
		\end{center}
		\subcaption{$q=0.9$}
	\end{subfigure}	 
		\begin{subfigure}{0.32\textwidth}
		\begin{center}	
			\includegraphics[width=0.9\textwidth]{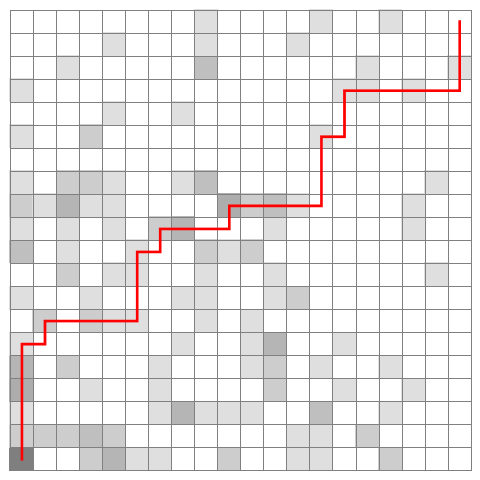}
		\end{center}
		\subcaption{$q=0.95$}
	\end{subfigure} 
    	\begin{subfigure}{0.32\textwidth}
		\begin{center}	
			\includegraphics[width=0.9\textwidth]{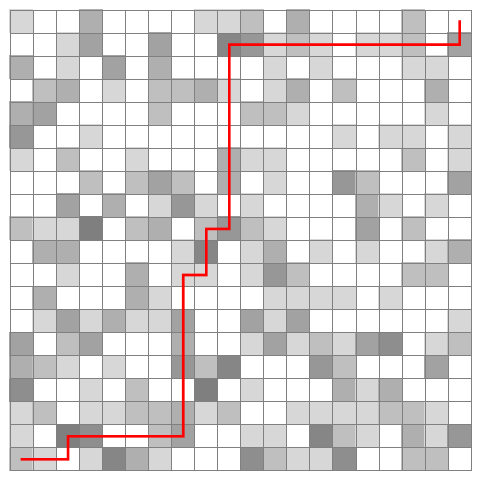}
		\end{center}
		\subcaption{$q=0.99$}
	\end{subfigure}  
	\caption{LPP realisation for $n=m=20$ and $\alpha=0.8$. The darker lattice sites indicate larger geometric weights.}   \label{Fig_LPP realisation}
\end{figure}

Consider the square lattice equipped with the random vertex weights
$\{\omega_{ij}\}$. Let $\Pi_{n,m}$ denote the set of all up-right lattice
paths from $(1,1)$ to $(n,m)$, where each step is either $(1,0)$ or
$(0,1)$. The corresponding last passage time is defined by
\begin{equation} \label{def of Gnm passage time}
G_{n,m}
=
\max_{\pi\in\Pi_{n,m}}
\sum_{(i,j)\in\pi}\omega_{ij}.
\end{equation}
Any path attaining the maximum is called a geodesic. The resulting random growth model is referred to as the \emph{spatially inhomogeneous geometric LPP} model; see Figure~\ref{Fig_LPP realisation}. 

As previously mentioned, such a spatially inhomogeneous directed LPP and related integrable growth models have been studied extensively. Integrable models with parameter-dependent weights were introduced in \cite{Joh08,FR05,BP08}, while variational and PDE approaches to limit shapes were developed in \cite{RT08,Em16,Cald15,CG21,EJS21}. Spatial inhomogeneity also gives rise to diverse fluctuation behaviours, including Tracy--Widom, Gaussian, and BBP-type limits in various models \cite{GTW02,Em16a,BP08,BC11,KPS19,BCD21}.   
Our model belongs to this class through a deterministic row--column dependent geometric parameter profile. 

In the asymptotic analysis of the inhomogeneous LPP, it is natural to
choose the parameter $q$ in \eqref{def of xi yj weights} according to
\begin{equation}\label{def of q scale}
q=e^{-\lambda/N},
\qquad
\lambda\in(0,\infty),
\end{equation}
where $N$ denotes the macroscopic scale of the lattice, so that
$n,m=O(N)$. 
Indeed, since $x_i$ and $y_j$ depend exponentially on $q$, keeping $q$ fixed would make them exponentially small for lattice sites $i,j=O(N)$. The scaling \eqref{def of q scale} precisely compensates for this effect, ensuring that both parameters remain of order $O(1)$ across the macroscopic lattice. Moreover, this double-scaling regime naturally appears in the asymptotic theory of $q$-deformed random matrix models; see \cite{BFO26,BJM26,BJO25,JS26}.

\begin{figure}[t]
\centering

\begin{tblr}{
    width = \textwidth,
    colspec = {Q[c,m,2cm] X[c,m] X[c,m]},
    rows = {rowsep = 6pt},
    row{1} = {bg = gray!12, font = \small},
    column{1} = {bg = gray!8},
    cell{1}{1} = {bg = gray!22},
    hlines = {gray!55, 0.5pt},
    vlines = {gray!55, 0.5pt},
    row{2} = {ht = 4.5cm},
    row{3,4} = {ht = 0.75cm},
}
&
\shortstack{
    $\alpha \in (0,1)$ fixed\\[2pt]
    \emph{Weak inhomogeneity} (Thm.~\ref{Thm_LPP geometric weak})
}
&
\shortstack{
    $\alpha = e^{-\eta/N}$ with $\eta \in (0,\infty)$\\[2pt]
    \emph{Strong inhomogeneity}  (Thm.~\ref{Thm_LPP geometric strong})
}
\\
\raisebox{1.8cm}{\shortstack{Weight\\on lattice}}
&
\includegraphics[width=.88\linewidth]{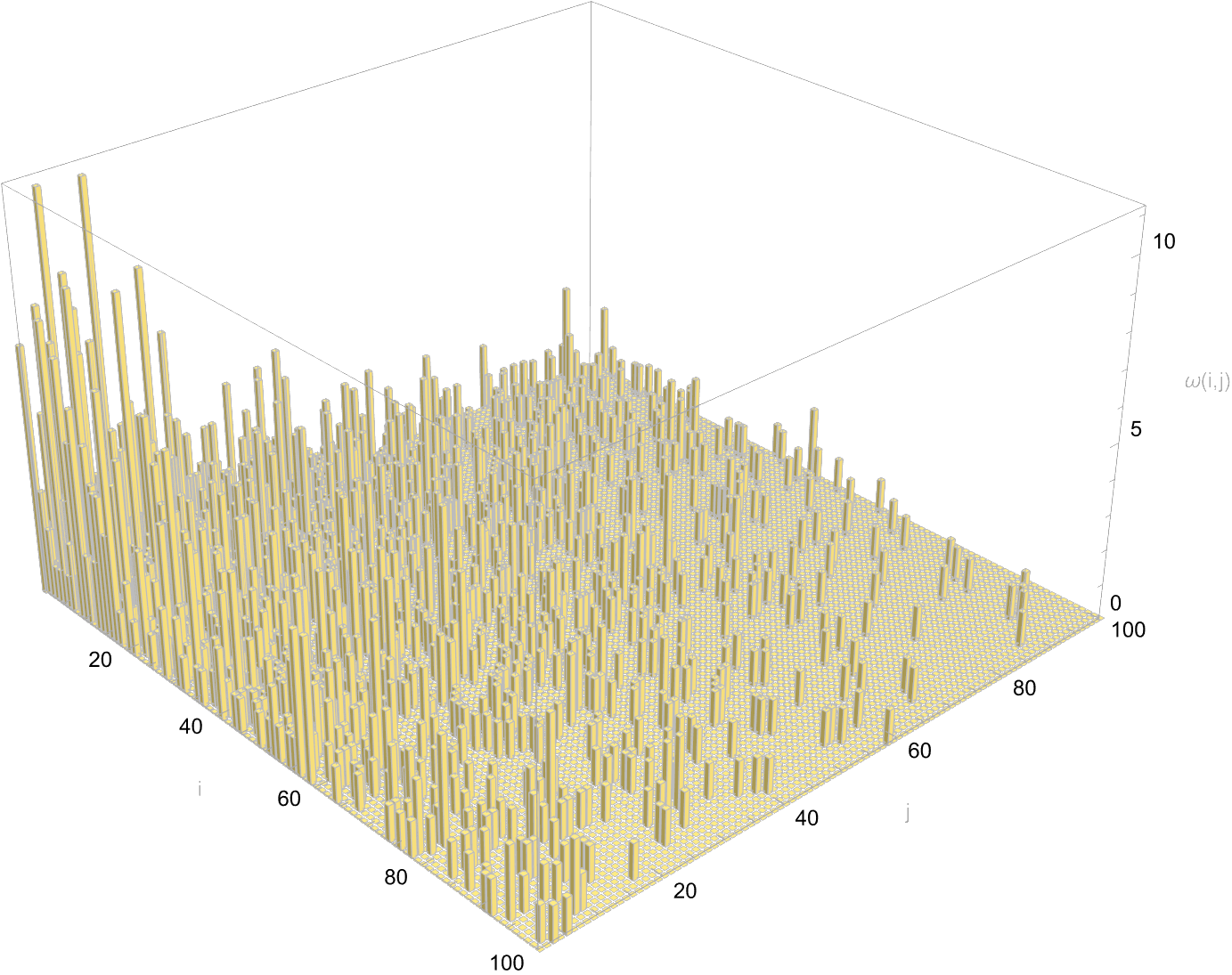}
&
\includegraphics[width=.88\linewidth]{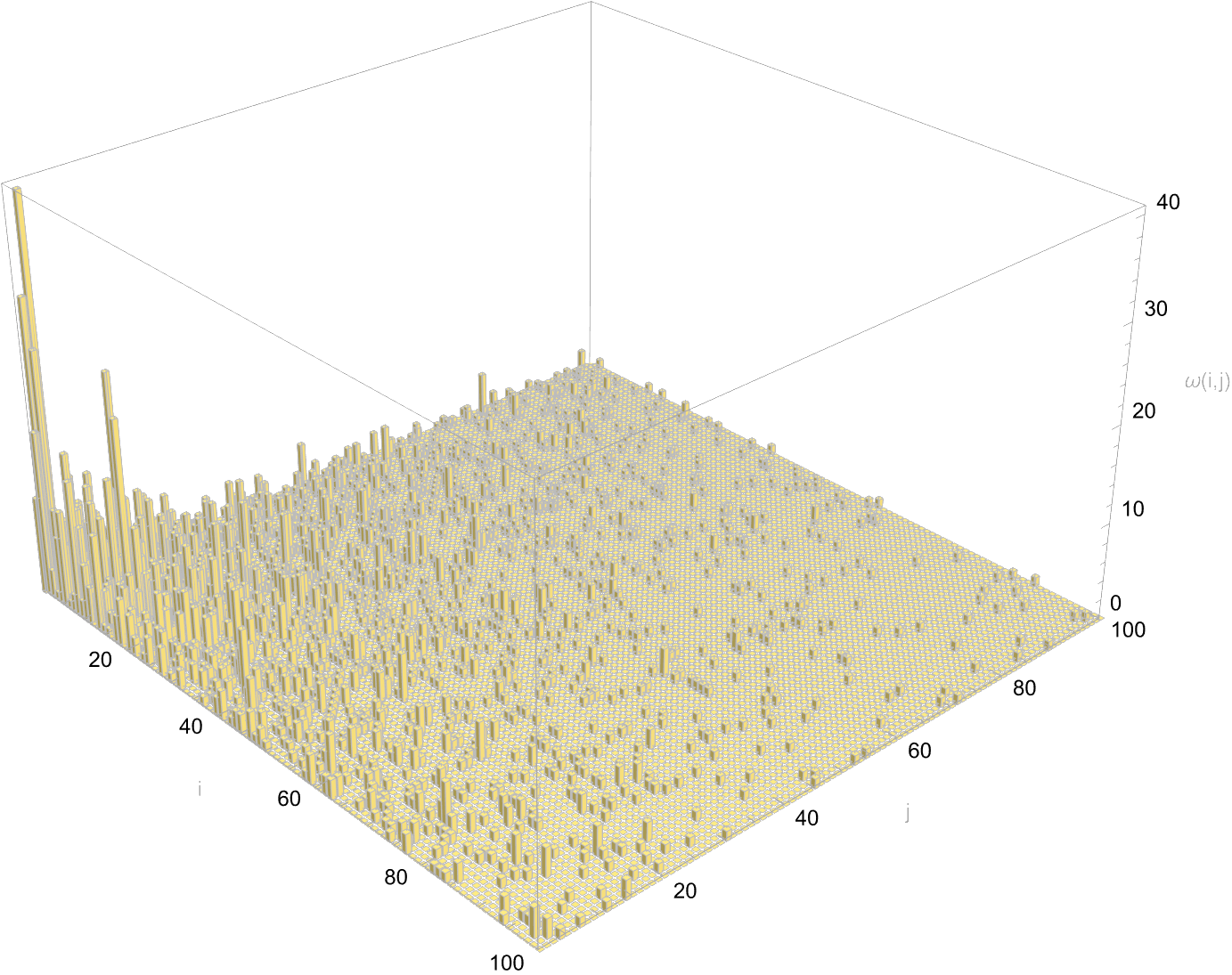}
\\

Time scale
&
$G_{\lfloor \gamma N \rfloor,N} = O(N)$
&
$G_{\lfloor \gamma N \rfloor,N} = O(N\log N)$
\\

Fluctuations
&
Airy Fredholm determinant
&
Bessel Fredholm determinant
\end{tblr} 
\caption{Weak and strong inhomogeneity regimes in spatially inhomogeneous LPP.
The first row shows sample realisations of the random environment
$\{\omega_{ij}\}$ on the lattice $\{1,\ldots,100\}^2$.}
\label{Fig_Strong/weak}
\end{figure}

In the scaling regime \eqref{def of q scale}, the parameter $q$
determines the macroscopic spatial variation of the environment.
The remaining parameter $\alpha$ controls the overall strength of the
inhomogeneity. As $\alpha$ increases, the weights become larger, while
their mass becomes increasingly concentrated near the origin; see
Figure~\ref{Fig_Strong/weak}. This naturally gives rise to two distinct
asymptotic regimes.
\begin{itemize}
    \item \emph{Weak inhomogeneity.} Here $\alpha\in(0,1)$ is fixed.
    The weights remain of order one throughout the macroscopic lattice,
    including both the origin and the endpoint $(n,m)$. Thus the
    inhomogeneity is present but remains moderate across the entire
    lattice.
    \smallskip
    \item \emph{Strong inhomogeneity.} Here
    $\alpha=e^{-\eta/N}$ for some fixed $\eta\in(0,\infty)$. In this
    regime, the weights become significantly larger and are strongly
    concentrated near the origin, leading to a markedly different
    asymptotic behaviour.
\end{itemize}
See Figure~\ref{Fig_Strong/weak} for an overview of these two regimes and the corresponding results. 

We now state our first result concerning the weak inhomogeneity regime.
For this, we recall that the GUE Tracy--Widom distribution
$F_{\rm TW}(s)$ is given by
\begin{equation} \label{def of TW}
F_{\rm TW}(s) = \det\bigl(I-K_{\Ai}\bigr)_{L^{2}(s,\infty)} = \exp\biggl( -\int_s^{\infty}(x-s)u(x)^2\,dx \biggr),
\end{equation}
where
\begin{equation}
K_{\Ai}(x,y) = \frac{\Ai(x)\Ai'(y)-\Ai'(x)\Ai(y)}{x-y}
\end{equation}
is the Airy kernel, and $u$ denotes the Hastings--McLeod solution of the
Painlev\'e II equation.  

We introduce the functions which appear in the limit theorem. 
Fix $\gamma\ge 1$, $\alpha\in(0,1)$, and $\lambda>0$. 
Define 
\begin{equation} \label{def of shape function weak}
\omega(\gamma,\alpha,\lambda)
:=
-\frac{2}{\lambda}
\log
\Bigg(
\frac{
\sqrt{(1-\alpha^2 e^{-\lambda})(1-\alpha^2 e^{-\gamma\lambda})}
-
\sqrt{\alpha^2(1-e^{-\lambda})(1-e^{-\gamma\lambda})}
}{
1-\alpha^2 e^{-(\gamma+1)\lambda}
}
\Bigg).
\end{equation}
We also write
\begin{align}
\label{def of a weak LPP}
   \mathfrak{a} \equiv   \mathfrak{a}_{\gamma,\alpha,\lambda}
    &:=
    e^{-\lambda}
    \Bigg(
    \frac{
    \sqrt{(1-\alpha^2 e^{-\lambda})(1-\alpha^2 e^{-\gamma\lambda})}
    -
    \sqrt{\alpha^2(1-e^{-\lambda})(1-e^{-\gamma\lambda})}
    }{
    1-\alpha^2 e^{-(\gamma+1)\lambda}
    }
    \Bigg)^2, 
    \\
   \mathfrak{b} \equiv   \mathfrak{b}_{\gamma,\alpha,\lambda} \label{def of b weak LPP}
    &:=
    e^{-\lambda}
    \Bigg(
    \frac{
    \sqrt{(1-\alpha^2 e^{-\lambda})(1-\alpha^2 e^{-\gamma\lambda})}
    +
    \sqrt{\alpha^2(1-e^{-\lambda})(1-e^{-\gamma\lambda})}
    }{
    1-\alpha^2 e^{-(\gamma+1)\lambda}
    }
    \Bigg)^2 
\end{align} 
and define 
\begin{equation} \label{def of theta weak LPP}
\theta_{\gamma,\alpha,\lambda}
:=
\begin{cases}
\displaystyle
\frac{
\sqrt{(1-\mathfrak{a})
(1-e^{-(\gamma-1)\lambda}\mathfrak{b})}
-
\sqrt{(1-\mathfrak{b})
(1-e^{-(\gamma-1)\lambda}\mathfrak{a})}
}{
\lambda\mathfrak{a}
\sqrt{
(1-\mathfrak{a})
(1-e^{-(\gamma-1)\lambda}\mathfrak{a})
(\mathfrak{b}-\mathfrak{a})
}
},
& \textup{if } \Xi_{\gamma,\alpha,\lambda}\ge 1,
\smallskip\\
\displaystyle
\frac{
\sqrt{(1-\mathfrak{a})
(1-e^{-(\gamma-1)\lambda}\mathfrak{b})}
+
\sqrt{(1-\mathfrak{b})
(1-e^{-(\gamma-1)\lambda}\mathfrak{a})}
}{
\lambda\mathfrak{a}
\sqrt{
(1-\mathfrak{a})
(1-e^{-(\gamma-1)\lambda}\mathfrak{a})
(\mathfrak{b}-\mathfrak{a})
}
},
& \textup{if } \Xi_{\gamma,\alpha,\lambda}< 1 ,
\end{cases}
\end{equation}
where 
\begin{equation} \label{def of critical quantity weak LPP}
    \Xi_{\gamma,\alpha,\lambda}
    :=
    e^{-\lambda}
    \bigl(1+\alpha^2-\alpha^2e^{-\gamma\lambda}\bigr) . 
\end{equation}
We then set
\begin{equation} \label{def of sigma weak LPP}
    \sigma(\gamma,\alpha,\lambda)
    :=
    \frac{1}{ \lambda\, \mathfrak{a} \, \theta_{\gamma,\alpha,\lambda}^{2/3}} .
\end{equation}
Let $\lfloor x\rfloor$ be the largest integer less than or equal to $x$. We are now ready to present our result.

\begin{thm}[\textbf{Geometric LPP at weak inhomogeneity}] \label{Thm_LPP geometric weak}
Let $\alpha \in (0,1)$ be fixed and $\gamma \ge 1$. Set $q = e^{-\lambda/N}$ for some $\lambda \in (0,\infty)$. Let $G_{n,m}$ denote the last passage time with weights given by~\eqref{def of xi yj weights}. Then as $N \to \infty$, 
\begin{equation} \label{LLN for inhomo LPP}
  \frac{ G_{\lfloor\gamma N \rfloor, N}  }{N}       \to \omega(\gamma,\alpha,\lambda) 
\end{equation}
in probability, where the shape function $\omega(\gamma,\alpha,\lambda)$ is given by \eqref{def of shape function weak}. Furthermore, for $s \in \R$, 
\begin{equation} \label{Fluctuation for inhomo LPP weak}
\lim_{ N \to \infty } \mathbb{P} \bigg[  \frac{ G_{\lfloor\gamma N \rfloor, N} - N \,\omega(\gamma,\alpha,\lambda) }{ \sigma(\gamma,\alpha,\lambda) N^{1/3}    }  \le s \bigg] = F_{ \rm TW }(s),  
\end{equation} 
where $\sigma( \gamma,\alpha,\lambda )$ is given by \eqref{def of sigma weak LPP}. 
\end{thm}

\begin{figure}[t]  
	\begin{subfigure}{0.45\textwidth}
		\begin{center}	
			\includegraphics[width=\textwidth]{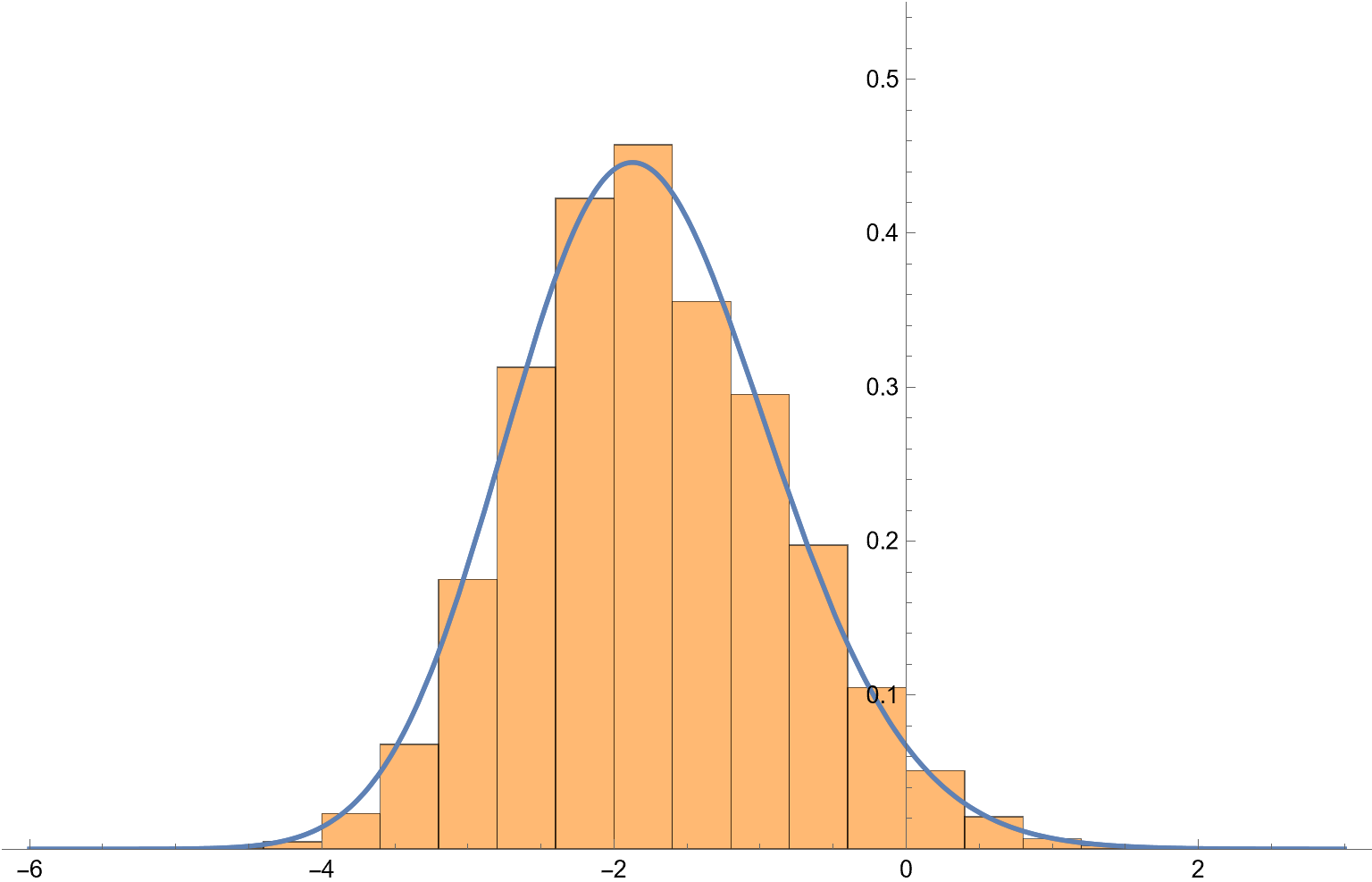}
		\end{center}
		\subcaption{$\lambda=1, \alpha=0.8, \gamma=1$}
	\end{subfigure}	\qquad 
		\begin{subfigure}{0.45\textwidth}
		\begin{center}	
			\includegraphics[width=\textwidth]{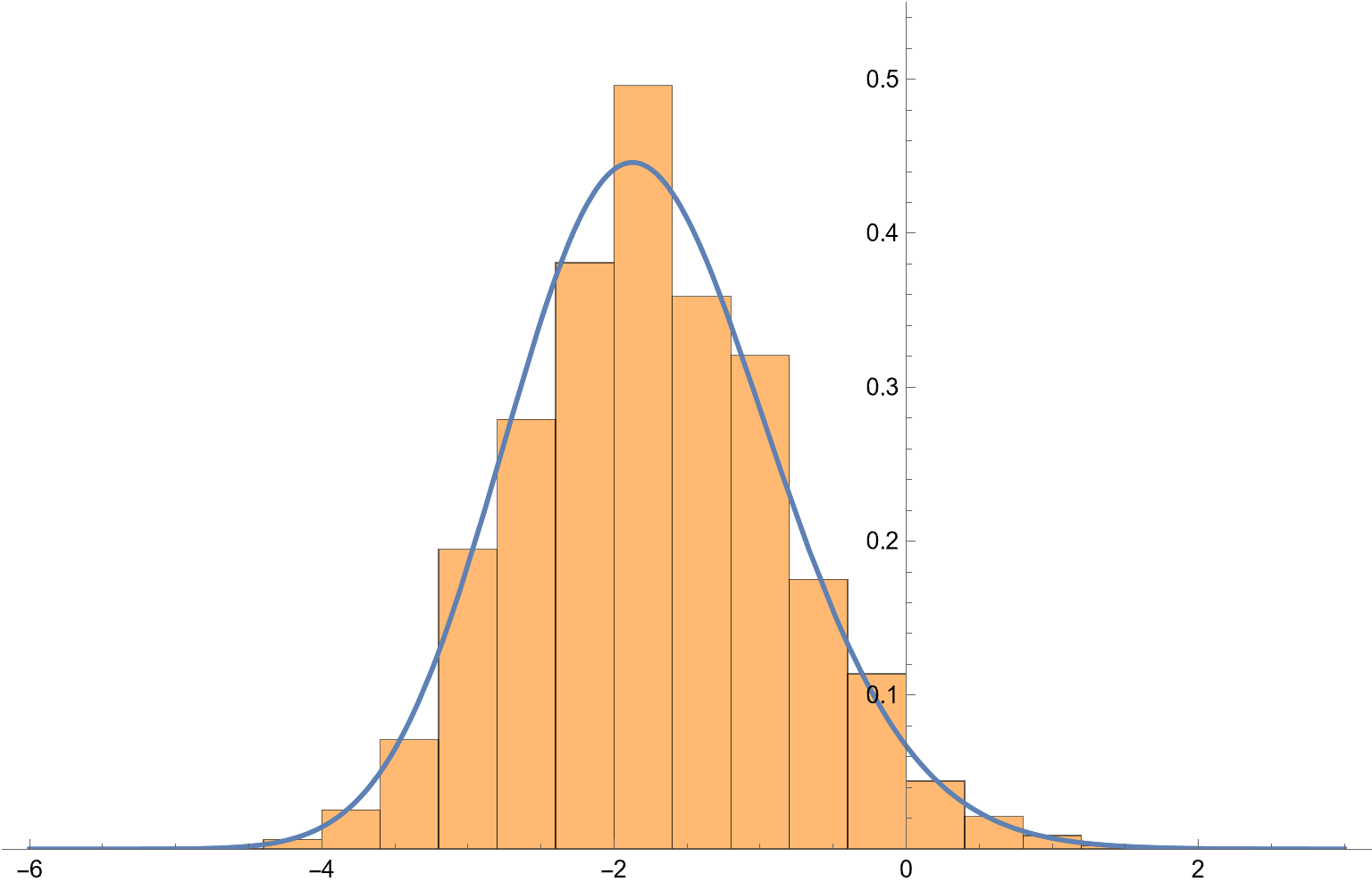}
		\end{center}
		\subcaption{$\lambda=0.5, \alpha=0.6, \gamma=2$}
	\end{subfigure} 
	\caption{ The histogram was made by computing 20000 samples of $\frac{G_{ \lfloor \gamma N \rfloor, N} - N \omega(\gamma, \alpha,\lambda)}{\sigma(\gamma, \alpha, \lambda) N^{1/3}}$ with $N=400$ for different choice of parameters. The blue curve is the density $F_{\mathrm{TW}}'(x)$ of the Tracy-Widom distribution.  }  \label{Fig_LPT weak}    
\end{figure}

See Figure~\ref{Fig_LPT weak} for a numerical illustration of
Theorem~\ref{Thm_LPP geometric weak}.

We emphasise that \eqref{LLN for inhomo LPP} provides one of the very few
known explicit shape functions for LPP, especially in a spatially
inhomogeneous environment. Furthermore, the fluctuation formula \eqref{Fluctuation for inhomo LPP weak} establishes the Tracy--Widom fluctuations predicted by KPZ universality; cf. \cite[Conjecture 1.9]{Joh00}.

\begin{rem}[Monotonicity of the shape function]
It is straightforward to verify that the shape function
$\omega(\gamma,\alpha,\lambda)$ is decreasing in $\lambda$ and increasing
in $\alpha$. This is consistent with the microscopic structure of the
model. Indeed, the geometric weights satisfy
\begin{equation}
\mathbb{E}[\omega_{ij}]
=
\frac{x_i y_j}{1-x_i y_j}
=
\frac{1}{1-\alpha^2q^{\,i+j-2}}-1.
\end{equation}
Since $q=e^{-\lambda/N}$ decreases as $\lambda$ increases, so does
$x_i y_j=\alpha^2q^{\,i+j-2}$. Consequently, the mean weight at each
vertex decreases with $\lambda$, and this monotonicity is inherited by
the macroscopic shape function. Likewise,
$\omega(\gamma,\alpha,\lambda)$ is increasing in~$\alpha$.
\end{rem}

\begin{rem}[Homogeneous limit]
In the homogeneous limit $\lambda\downarrow0$, equivalently $q\uparrow1$, the
inhomogeneous parameters in \eqref{def of xi yj weights} become
$x_i=y_j=\alpha$. Thus the model reduces to homogeneous geometric LPP. 
In this limit, the shape function satisfies
\begin{equation} \label{homogeneous limit of shape weak LPP}
    \lim_{\lambda\downarrow0}
    \omega(\gamma,\alpha,\lambda)
    =
    \frac{(1+\alpha\sqrt{\gamma})^2}{1-\alpha^2}-1,
\end{equation}
which agrees with \cite[Theorem~1.1]{Joh00}. 
Moreover, a direct expansion of
\eqref{def of sigma weak LPP} gives
\begin{equation} \label{homogeneous limit of sigma weak LPP}
    \lim_{\lambda\downarrow0}
    \sigma(\gamma,\alpha,\lambda)
    =
    \frac{
    \alpha^{1/3}
    \bigl((1+\alpha\sqrt{\gamma})(\sqrt{\gamma}+\alpha)\bigr)^{2/3}
    }{
    \gamma^{1/6}(1-\alpha^2)
    }.
\end{equation}
This agrees with the fluctuation scale for homogeneous geometric LPP in \cite[Theorem~1.2]{Joh00}.
\end{rem} 

We now turn to the strong inhomogeneity regime. For this, recall that the Bessel kernel is given by
\begin{equation}\label{def of K Bessel standard}
K^{\rm Bessel}_\nu(x,y)
:=
\frac{
J_\nu(\sqrt{x})\,\sqrt{y}\,J_\nu'(\sqrt{y})
-
\sqrt{x}\,J_\nu'(\sqrt{x})\,J_\nu(\sqrt{y})
}{
2(x-y)
}.
\end{equation} 
For the strong inhomogeneity regime, we also define
\begin{equation} \label{def of kappa strong LPP}
    \kappa_{\gamma,\lambda}
    := \frac{2}{\lambda}
    \sqrt{
    e^{\lambda}(1-e^{-\lambda})(1-e^{-\gamma\lambda})
    } .
\end{equation}
Then we have the following result. 

\begin{thm}[\textbf{Geometric LPP at strong inhomogeneity}] \label{Thm_LPP geometric strong} 
Let $\alpha =e^{-\eta/N}$ for some $\eta \in (0,\infty)$ and $\gamma \ge 1$. Set $q = e^{-\lambda/N}$ for some $\lambda \in (0,\infty)$. Let $G_{n,m}$ denote the last passage time with weights given by~\eqref{def of xi yj weights}. Then as $N \to \infty$,  
\begin{equation}  \label{LLN for inhomo LPP alpha close to 1}
  \frac{  G_{\lfloor\gamma N \rfloor, N}  }{N \log N}   \to   \frac{2}{\lambda}   
\end{equation}
in probability. 
Furthermore, for $s \in \R$, 
\begin{equation} \label{Fluctuation for inhomo LPP strong}
    \lim_{N\to\infty}
    \mathbb{P}
    \bigg[
    \frac{
    G_{\lfloor\gamma N\rfloor,N}
    -
    \frac{2}{\lambda}N\log N
    }{
    N/\lambda
    }
    \le s
    \bigg]
    =
    \det\bigl(I-K_{\mathsf{a}}^{\rm Bessel}\bigr)_{
    L^2(0,\kappa_{\gamma,\lambda}^{2}e^{-\lambda}e^{-s})}, 
\end{equation}
where $\kappa_{\gamma,\lambda}$ is given by \eqref{def of kappa strong LPP} and 
\begin{align}
\mathsf{a}
= \frac{2\eta}{\lambda} -1.  
\end{align}
\end{thm}

\begin{figure}[t]  
	\begin{subfigure}{0.45\textwidth}
		\begin{center}	
			\includegraphics[width=\textwidth]{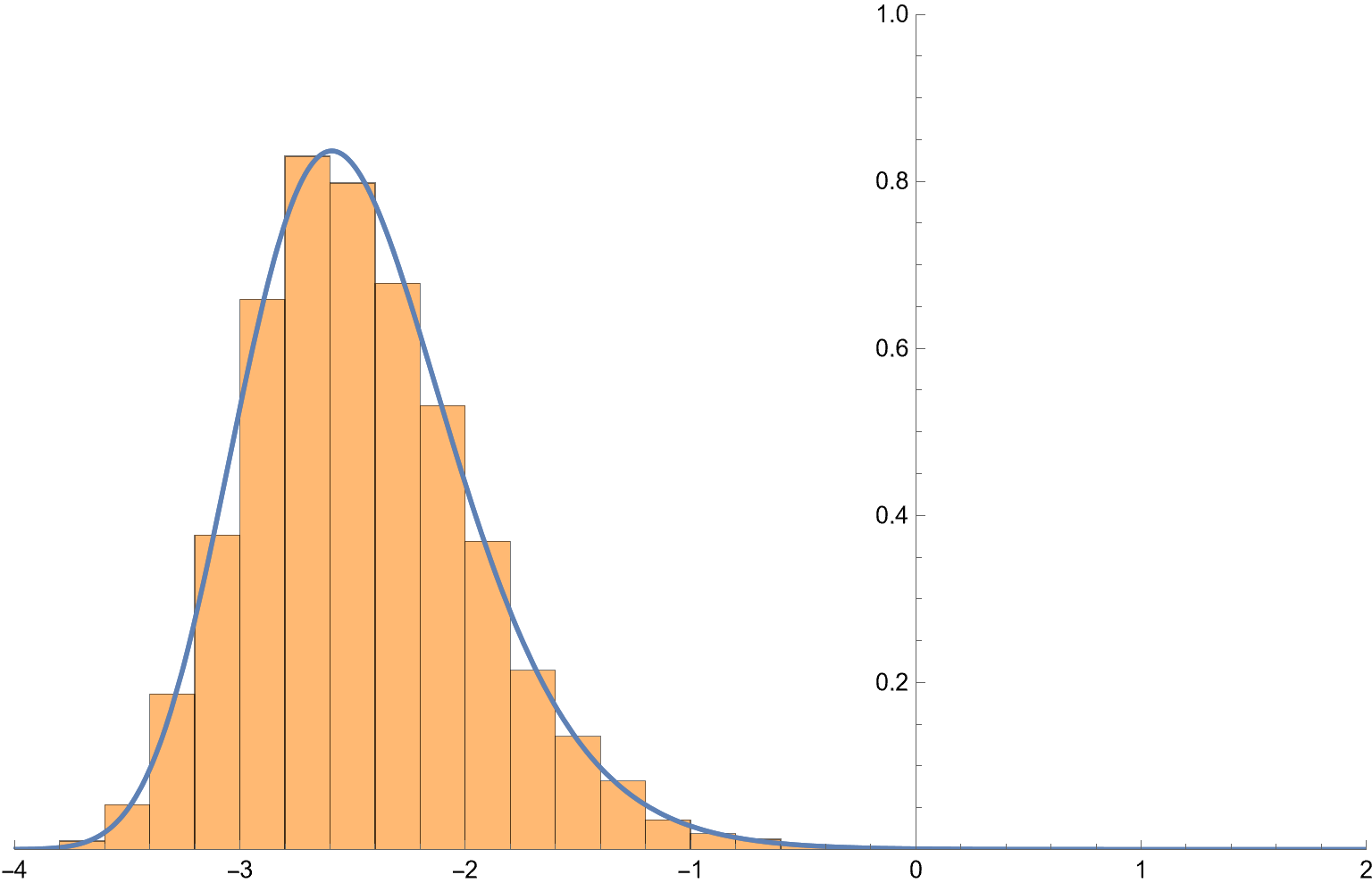}
		\end{center}
		\subcaption{$\lambda=0.5, \eta=1, \gamma=1$}
	\end{subfigure}	\qquad 
		\begin{subfigure}{0.45\textwidth}
		\begin{center}	
			\includegraphics[width=\textwidth]{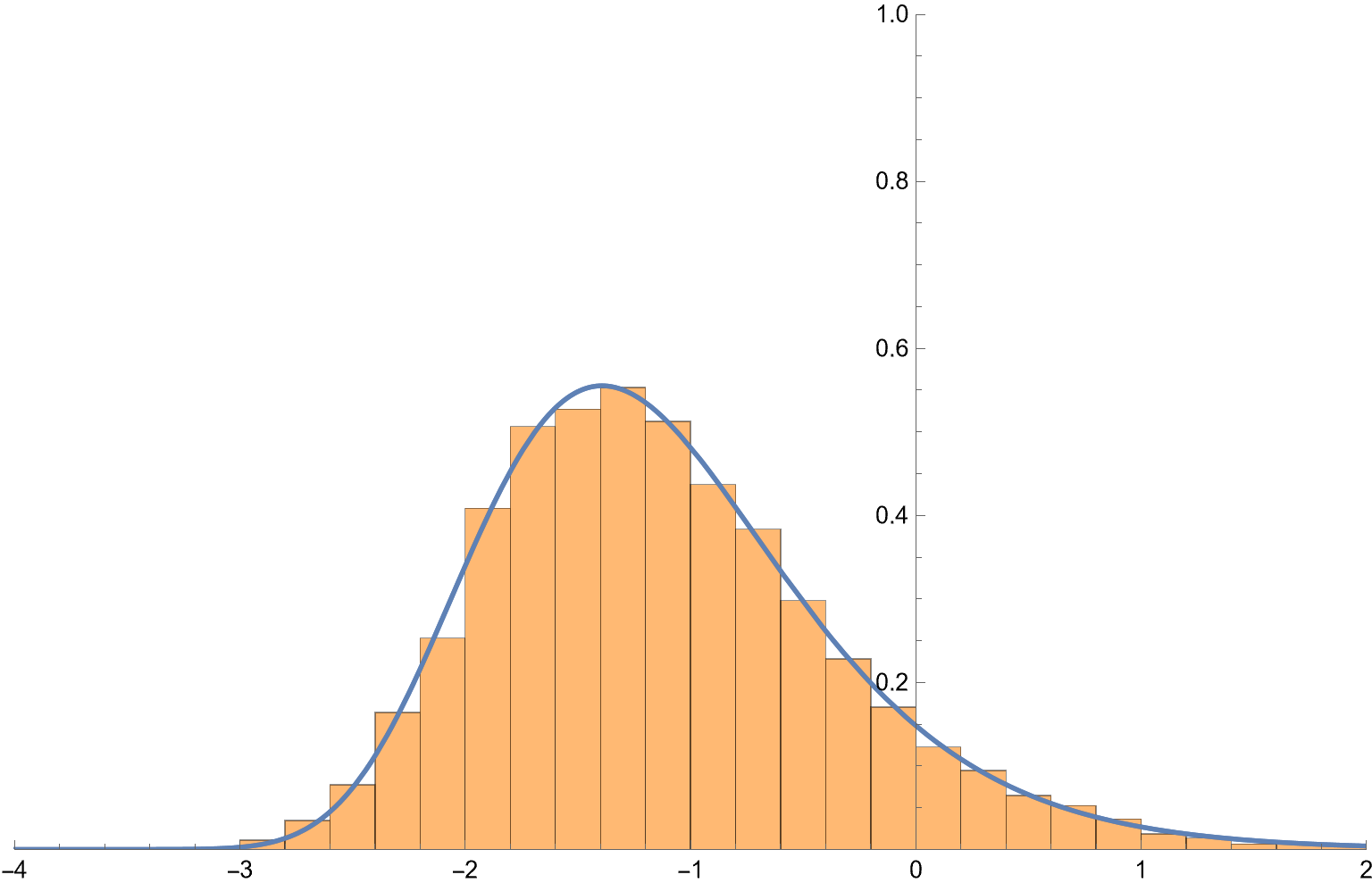}
		\end{center}
		\subcaption{$\lambda=0.8, \eta=0.8, \gamma=2$}
	\end{subfigure} 
	\caption{ The histogram was made by computing 10000 samples of $\frac{G_{ \lfloor \gamma N \rfloor, N} - \frac{2}{\lambda} N\log N }{ N/\lambda }$ with $N=400$ for different choice of parameters. The blue curve is the density of the Bessel Fredholm determinant, the right-hand side of \eqref{Fluctuation for inhomo LPP strong}.  }  \label{Fig_LPT strong}    
\end{figure}

See Figure~\ref{Fig_LPT strong} for a numerical illustration of
Theorem~\ref{Thm_LPP geometric strong}.

A particularly notable feature of Theorem~\ref{Thm_LPP geometric strong}
is the emergence of hard-edge Bessel fluctuations. While Tracy--Widom
fluctuations are ubiquitous in integrable probability, limiting laws
governed by the hard-edge Bessel kernel are comparatively rare.
Related Bessel fluctuations on an infinite quadrant were obtained in \cite[Theorem~2.21]{BO24}.

\begin{rem}[Non-standard time scale and locality of the shape function]
The $O(N\log N)$ scaling is at least formally consistent with the weak inhomogeneity result in Theorem~\ref{Thm_LPP geometric weak}. To see this, note that the shape function in \eqref{def of shape function weak} satisfies
\begin{equation}\label{asymp of omega as alpha to 1}
\omega(\gamma,\alpha,\lambda)
=
-\frac{2}{\lambda}\log(1-\alpha)
+\frac{1}{\lambda}
\log\Bigl((1-e^{-\lambda})(1-e^{-\gamma\lambda})\Bigr)
+O(1-\alpha),
\qquad
\alpha\to1.
\end{equation}
Hence, under the scaling $\alpha=e^{-\eta/N}$, the leading-order term of the passage time is naturally of order $N\log N$, in contrast to the $O(N)$ behaviour in the weak inhomogeneity regime. We remark that such a non-standard logarithmic scaling also appears in other growth models, for example in long-range first passage percolation~\cite{CD16}. 

Another notable feature of the strong inhomogeneity regime is that
the shape function \eqref{LLN for inhomo LPP alpha close to 1} is
independent of the aspect ratio $\gamma$ of the rectangular lattice.
This reflects the fact that, as illustrated in
Figure~\ref{Fig_Strong/weak}, the weights are overwhelmingly
concentrated near the starting point $(1,1)$. Consequently, the shape
function exhibits a \textit{locality} property: it is asymptotically determined by the environment near the starting point, rendering it insensitive to the global geometry of the lattice.
\end{rem}

\subsection{Exponential inhomogeneous LPP} 

We now consider the spatially inhomogeneous LPP with
exponential weights. As in the homogeneous setting \cite{Joh00}, the exponential model arises as a scaling limit of the geometric model, which usually leads to a considerably simpler theory.  

We now introduce the model and state the corresponding limit theorems.  
Let $\Exp(t)$ denote the exponential distribution with rate parameter
$t>0$, i.e. 
\begin{equation}
\mathbb{P}(X\in dx)
=
t e^{-tx}\,dx,
\qquad x\ge0.
\end{equation}
Assign to each lattice site $(i,j)$ an independent random variable
\begin{equation}
\widetilde{\omega}_{ij}
\sim
\Exp(t_i+t_j).
\end{equation}
For given $\lambda,\eta>0$, we consider the linearly varying parameters
\begin{equation}\label{def of tk exponential}
t_k=\lambda(k+\beta), \qquad \beta:= \frac{\eta}{\lambda}-1.
\end{equation} 
The associated last passage time is defined by
\begin{equation}\label{def of Hnm exponential LPP time}
H_{n,m}
=
\max_{\pi\in\Pi_{n,m}}
\sum_{(i,j)\in\pi}
\widetilde{\omega}_{ij}.
\end{equation} 
We refer to this as the \emph{spatially inhomogeneous exponential LPP} model. The special case $\gamma=\lambda=1$ was investigated in
\cite{Joh08}; see also \cite{BP08}. 
For this model, we have the following. 

\begin{thm}[\textbf{LPP with inhomogeneous exponential weight}]\label{Thm_exponential}
Let $\eta , \lambda >0$ be fixed and $\gamma \ge 1$.  
Let $H_{n,m}$ denote the last passage time given in \eqref{def of Hnm exponential LPP time}.  Then as $N \to \infty,$ 
\begin{equation} \label{LLN for inhomo LPP exp}
 \frac{  H_{\lfloor\gamma N \rfloor, N} }{\log N} \to \frac{2}{\lambda}    
\end{equation}
in probability.  
Furthermore, for $\xi \in \R$, 
\begin{equation} \label{Fluctuation limit for exponential LPP}
\lim_{ N \to \infty } \mathbb{P}\Big[ H_{\lfloor\gamma N \rfloor, N} -\frac{2}{\lambda} \log N \le \frac{\xi+\log \gamma }{ \lambda} \Big]  =  \det (I-K_\beta)_{ L^2(\xi,\infty) },  
\end{equation}
where
\begin{equation}\label{def of K Bessel exponential}
K_\beta (x,y) = 4 e^{ -\frac{x+y}{2} } K_{ 2\beta+1 }^{ \rm Bessel }( 4 e^{-x},4e^{-y} ).
\end{equation}
Here,  $K_\nu^{ \rm Bessel }$ is the standard Bessel kernel \eqref{def of K Bessel standard}. 
\end{thm}

\begin{figure}[t]  
\begin{subfigure}{0.45\textwidth}
		\begin{center}	
			\includegraphics[width=\textwidth]{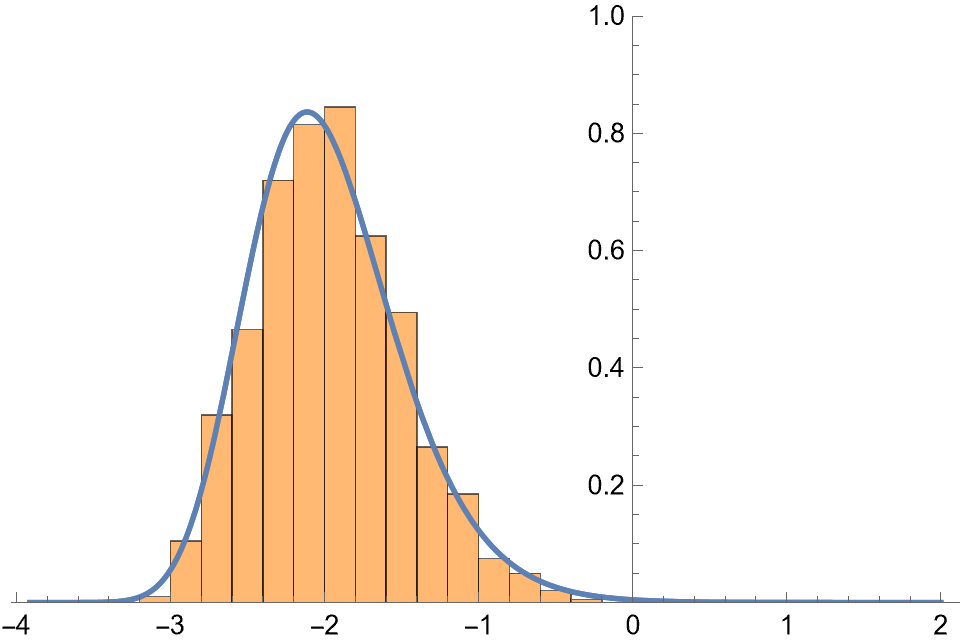}
		\end{center}
		\subcaption{$\eta=2, \beta=1$}
	\end{subfigure} \qquad 
	\begin{subfigure}{0.45\textwidth}
		\begin{center}	
			\includegraphics[width=\textwidth]{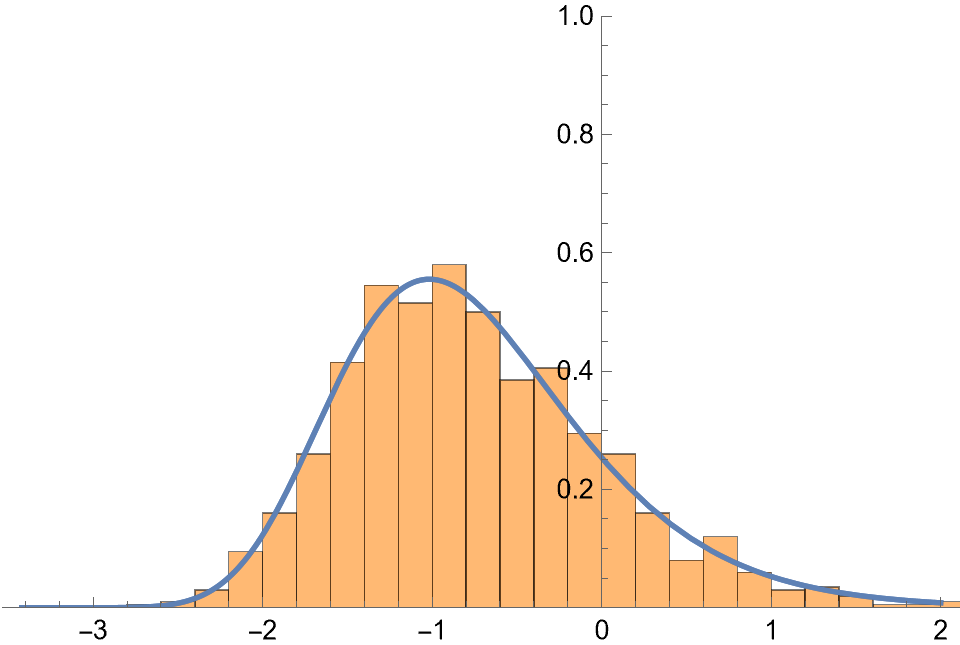}
		\end{center}
		\subcaption{$\eta=1, \beta=0$}
	\end{subfigure}	
	\caption{ The histogram was made by computing 1000 samples of $\lambda(H_{\lfloor\gamma N \rfloor, N} -\frac{2}{\lambda} \log N) -\log \gamma $ with $N=400$. Here, $\gamma=\lambda=1$. The blue curve is the density of $\det (I-K_\beta)_{ L^2(\xi,\infty) }.$  }  \label{Fig_exponential fluctuation simulation}
\end{figure}

See Figure~\ref{Fig_exponential fluctuation simulation} for a numerical illustration of Theorem~\ref{Thm_exponential}.

As in Theorem~\ref{Thm_LPP geometric strong}, the independence of the
shape function from the aspect ratio $\gamma$ may be interpreted as a
manifestation of the locality phenomenon. We also note that the right-hand side of \eqref{Fluctuation limit for exponential LPP} interpolates between the Gumbel and Tracy--Widom distributions under the appropriate scaling limits; see \cite{BF03,Joh08}. 

For the special case $\gamma=\lambda=1$, the formula in
Theorem~\ref{Thm_exponential} already appears in
\cite[Eq.~(1.6)]{Joh08}, where Johansson attributes it to Forrester as a private communication; see also \cite[Section~2.5]{FR05} for related discussions.
For the square case and $\beta\ge -1/2$, a proof also follows from \cite[Theorem 2.25]{BO24} after a logarithmic change of variables. Here we give a direct derivation for all $\gamma\ge1$ and $\beta>-1$ through the duality with the classical JUE.

We note that, unlike the geometric case considered in this paper, Theorem~\ref{Thm_exponential} can also be deduced from earlier results, notably \cite{BP08} (see also \cite{BS25}), via a suitable duality for the inhomogeneous model; see Remark~\ref{Rem_BP extended Wishart}.
Nevertheless, as a by-product of our analysis of the geometric model, we obtain a direct and self-contained proof of Theorem~\ref{Thm_exponential}, which is of independent interest.

\section{\texorpdfstring{$q$}{q}-deformed Jacobi unitary ensembles} \label{Section_qJUE}
 
This section is devoted to the results for the $q$-JUE.   

In general, the concept of $q$-deformation plays a prominent role in several areas of mathematics.
In random matrix theory, $q$-deformed ensembles arise as natural discrete counterparts of the classical random matrix ensembles. Their study is motivated from several complementary perspectives, including
$q$-Selberg type integrals \cite{As80,FLSY23}, Macdonald
polynomial theory \cite{BF00,Cu18,Kon96,Mat07,Ol21},
$q$-Virasoro constraints
\cite{CLPZ19,CZ22,KK22,LPSZ20,NZ17,SKAO96,MZ26}, moment and combinatorial
formulas \cite{BFO26,Co21,FLSY23,MPS20,Ven14,Wim14}, stochastic
processes and martingale structures \cite{BW05,BW10}, tiling models and
Gaussian free field fluctuations \cite{DK19,FR05,GH22}, and
connections with supersymmetric gauge theories
\cite{AOSV05,CLZ20,Sha15,ST13,LWYZ24}. We refer the reader to \cite{BFO26,BF25a,BJM26} and the references therein for further background. 
\medskip 

The $q$-JUE admits two natural degenerations as $q\to1$. Before
introducing the $q$-JUE itself, we first recall these two limiting
ensembles. 

The first is the classical JUE \cite[Chapter~3]{Fo10}. It is a system
$\boldsymbol x=\{x_j\}_{j=1}^N$ with $0<x_j<1$, whose joint probability
density is given by
\begin{equation} \label{def of classical JUE}
 \mathbb P_N^{\rm J}(d\boldsymbol x)
    =
    \frac{1}{Z_N^{\rm J}} \prod_{1\le j<k\le N}(x_j-x_k)^2
\prod_{j=1}^{N} x_j^{ \mathsf{a} } (1-x_j)^{ \mathsf{b} }\,dx_j,  
\end{equation}
where $\mathsf a,\mathsf b>-1$ are real parameters and
$Z_N^{\rm J}$ is the normalisation constant, which can be evaluated
explicitly by means of the Selberg integral \cite{FW08}.

The second example is the Meixner orthogonal polynomial ensemble \cite{Joh00}. For $\alpha\in(0,1)$ and $\beta>0$, the Meixner ensemble $\boldsymbol h=\{h_j\}_{j=1}^N$ is a discrete orthogonal polynomial ensemble on $\mathbb Z_{\geq 0}$, with probability mass function, relative to the counting measure on $\mathbb Z_{\geq 0}^N$, given by 
\begin{equation}
\label{eqn:def of Meixner ensemble general}
    \mathbb P_N^{\rm M}(\boldsymbol h)
    =
    \frac{1}{Z_N^{\rm M}}
    \prod_{1\leq j<k\leq N}(h_j-h_k)^2
    \prod_{j=1}^N
    \frac{(\beta)_{h_j}}{h_j!}\alpha^{2h_j}, 
\end{equation}
where $(\beta)_h$ denotes the Pochhammer symbol and $Z_N^{\rm M}$ is the corresponding normalisation constant. 

We now introduce the $q$-JUE. Here and throughout, we use the $q$-Pochhammer symbols
\begin{equation}
    (z;q)_k := \prod_{j=0}^{k-1}(1-zq^j),
    \qquad
    (z;q)_\infty := \prod_{j=0}^{\infty}(1-zq^j),
    \label{q-poch}
\end{equation}
with the convention $(z;q)_0=1$. For parameters
$\mathsf{a},\mathsf{b}>-1$, let
\begin{equation}\label{eqn:def of qjacobi weight}
w^{(\mathsf{a},\mathsf{b})}(x;q)
=
x^{\mathsf{a}}
\frac{(qx;q)_\infty}{(q^{\mathsf{b}+1}x;q)_\infty},
\end{equation}
which is the little $q$-Jacobi weight; see e.g. \cite[Section~4.1]{FLSY23}. This weight belongs to the
Askey--Wilson scheme of orthogonal polynomials and is associated with
the little $q$-Jacobi polynomials.  
We consider a system of particles
$\boldsymbol{x}=\{x_j\}_{j=1}^{N}$ supported on the one-sided exponential
$q$-lattice $\{1,q,q^2,\ldots\}$. Equivalently, one may write
$x_j=q^{h_j}$ with $h_j\in\mathbb Z_{\ge0}$. Let
\begin{equation} \label{def of Jackson measure}
    d_qx=(1-q)\sum_{h=0}^{\infty}q^h\delta_{q^h}(dx)
\end{equation}
denote the Jackson measure on $[0,1]$. The joint distribution of
$\boldsymbol{x}$ is
\begin{equation}\label{eqn:jpdf of qJUE}
\mathbb{P}_{N}(d\boldsymbol{x})
=
\frac{1}{Z_{N}}
\prod_{1\le j<k\le N}(x_j-x_k)^2
\prod_{j=1}^{N}
w^{(\mathsf{a},\mathsf{b})}(x_j;q)\,d_qx_j, 
\end{equation}
where $Z_N$ denotes the normalisation constant.  
Equivalently, under the parametrisation $x_j=q^{h_j}$, the induced
distribution of $\boldsymbol h=\{h_j\}_{j=1}^N$ is a discrete probability
measure on $\mathbb Z_{\ge0}^N$. More precisely, with respect to the
counting measure on $\mathbb Z_{\ge0}^N$, its probability mass function
is
\begin{equation}\label{eqn:jpdf of qJUE in h}
\mathbb P_N(\boldsymbol h)
=
\frac{1}{\widetilde Z_N}
\prod_{1\le j<k\le N}(q^{h_j}-q^{h_k})^2
\prod_{j=1}^N
\frac{(q^{\mathsf b+1};q)_{h_j}}{(q;q)_{h_j}}
q^{(\mathsf a+1)h_j}, 
\end{equation}
where $\widetilde Z_N$ is the corresponding normalisation constant.

\begin{figure}[t]
    \centering
    \resizebox{0.76\textwidth}{!}{
    \begin{tikzpicture}[
        every node/.style={align=center},
        mainbox/.style={
            draw=black,
            line width=0.75pt,
            rounded corners=2pt,
            fill=black!2,
             minimum width=7.4cm,
            minimum height=0.75cm,
            inner sep=5pt,
        },
        subbox/.style={
            draw=black,
            line width=0.75pt,
            rounded corners=2pt,
            fill=black!2,
            minimum width=3.7cm,
            minimum height=0.75cm,
            inner sep=5pt,
        },
        branchlabel/.style={
            font=\small,
            align=center,
            inner sep=0pt
        },
        brancharrow/.style={
            -{Latex[length=2.2mm,width=1.5mm]},
            line width=0.8pt,
            draw=black
        }
    ]

    \node[mainbox] (qjue) at (0,0) {
        \begin{tabular}{c}
            \(q\)-JUE $\{x_j=q^{h_j}\}_{j=1}^N$ ; \eqref{eqn:jpdf of qJUE}, 
            \eqref{eqn:jpdf of qJUE in h} 
        \end{tabular}
    };

    \node[subbox] (jue) at (-4.5,-3.2) {
        Classical JUE; \eqref{def of classical JUE}
    };

    \node[subbox] (meixner) at (4.5,-3.2) {
        Meixner ensemble; \eqref{eqn:def of Meixner ensemble general}
    };

    \coordinate (branch) at ($(qjue.south)+(0,-7mm)$);

    \draw[line width=0.8pt]
        (qjue.south) -- (branch);

    \draw[brancharrow]
        (branch) -- (jue.north);

    \draw[brancharrow]
        (branch) -- (meixner.north);

    \node[branchlabel] at (-3.85,-1.55) {
        $x_j$ fixed
        \\[0.15em]
        $w^{(\mathsf a,\mathsf b)}(x;q)
        \to
        x^{\mathsf a}(1-x)^{\mathsf b}$
    };

    \node[branchlabel] at (3.85,-1.55) {
        $h_j$ fixed
        \\[0.15em]
        $\mathsf b\to\beta-1,
        \quad
        q^{\mathsf a+1}\to\alpha^2$
    };

    \end{tikzpicture}
    }
    \caption{Two degenerations of the $q$-JUE as $q\to1$.}
    \label{fig:qJUE_degenerations}
\end{figure}
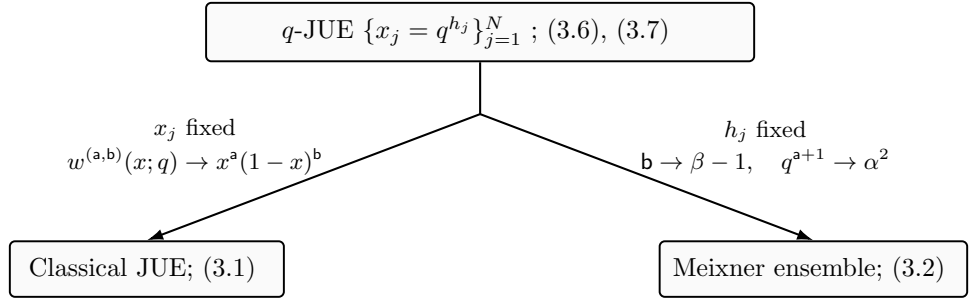

Since
\begin{equation}
    \lim_{q\to1}
    w^{(\mathsf a,\mathsf b)}(x;q)
    =
    x^{\mathsf a}(1-x)^{\mathsf b},
\end{equation}
the classical JUE is recovered from \eqref{eqn:jpdf of qJUE} by
letting $q\to1$ with $\mathsf a$ and $\mathsf b$ fixed. This limit is
taken in the continuous variables $x_j\in(0,1)$. 
The Meixner degeneration is obtained from the representation
\eqref{eqn:jpdf of qJUE in h} by instead keeping the lattice variables
$h_j$ fixed and scaling the parameters so that
\begin{equation}
\label{eqn:qJUE to Meixner finite scaling}
    \mathsf b\to\beta-1,
    \qquad
    q^{\mathsf a+1}\to\alpha^2,
    \qquad
    q\to1.
\end{equation}
For instance, one may take
\begin{equation}
    \mathsf a
    =
    \frac{\log\alpha^2}{\log q}-1,
    \qquad
    \mathsf b=\beta-1.
\end{equation}
Under this scaling, as $q \to 1$, we have 
\begin{equation*}
    \frac{(q^{\mathsf b+1};q)_h}{(q;q)_h}
    \to
    \frac{(\beta)_h}{h!},
    \qquad
    q^{(\mathsf a+1)h}
    \to 
    \alpha^{2h}, \qquad  \frac{q^{h_j}-q^{h_k}}{1-q}
    \to 
    -(h_j-h_k).
\end{equation*} 
The probability mass function \eqref{eqn:jpdf of qJUE in h} therefore converges to the Meixner ensemble \eqref{eqn:def of Meixner ensemble general}. 
Figure~\ref{fig:qJUE_degenerations} summarises these two degenerations of the $q$-JUE.

\medskip 

Unlike the classical JUE, the asymptotic theory of the
$q$-JUE has remained largely unexplored. In particular, its limiting spectral distribution, which may be viewed as a $q$-deformation of the Wachter distribution \cite{Wac80}, does not appear to have been identified in the existing literature. This should be compared with the recent identification of the $q$-deformed semicircle distribution \cite{BFO26} and the $q$-deformed Marchenko--Pastur law \cite{BJM26}. 

In order to describe the limiting macroscopic behaviour, we first define the empirical spectral distribution 
\begin{equation}
  \mu_N(dx) := \frac{1}{N} \sum_{j=1}^N \delta_{x_j} .
 \end{equation}
As in the classical JUE case, when describing the large-$N$ behaviour, one needs to consider the proper scaling of the Jacobi weight parameters $\mathsf{a}$ and $\mathsf{b}$. With sufficient generality we consider the scaling 
\begin{equation}\label{eqn:scale of alphabet}
    \mathsf{a}=c_{1}N+d_{1},\qquad \mathsf{b}=c_{2}N+d_{2}. 
\end{equation}
Here, $c_1,c_2\ge 0$, and we require $d_i>-1$ whenever $c_i=0$, for $i=1,2$.
Note that the parameters $c_1,c_2$ determine the macroscopic picture while $d_1,d_2$ determine the local spectral statistics. 

While the derivation of limiting spectral distributions for classical random matrix ensembles is by now well understood, the corresponding theory for $q$-deformed ensembles has emerged only recently. For the
$q$-deformed Gaussian unitary ensemble, the moment method, combined with
the Flajolet--Viennot combinatorial theory of orthogonal polynomials,
was employed in \cite{BFO26} to derive the $q$-deformed semicircle law.
This approach was subsequently extended in \cite{BJO25} to the
Al-Salam--Carlitz ensemble, where the limiting spectral measure was
shown to coincide with the limiting zero distribution of the associated
orthogonal polynomials. More recently, \cite{BJM26} combined these
algebraic-combinatorial techniques with potential theory \cite{HMO25} to establish the global asymptotics of the $q$-deformed Laguerre unitary ensemble, leading to the $q$-deformation of the Marchenko--Pastur law.

Building on these developments, we derive the limiting spectral
distribution of the $q$-JUE.  
For fixed $c_1,c_2 \ge 0$, we define 
\begin{equation}\label{eqn:def of endpoints ab}
    \begin{split}
        a&=\mathsf{s}\Bigg(\frac{\sqrt{(1-\mathsf{s}^{c_{1}+1})(1-\mathsf{s}^{c_{1}+c_{2}+1})}- \sqrt{\mathsf{s}^{c_{1}}(1-\mathsf{s})(1-\mathsf{s}^{c_{2}+1})}}{1-\mathsf{s}^{c_{1}+c_{2}+2}}\Bigg)^{2},
        \\
        b&=\mathsf{s}\Bigg(\frac{\sqrt{(1-\mathsf{s}^{c_{1}+1})(1-\mathsf{s}^{c_{1}+c_{2}+1})}+ \sqrt{\mathsf{s}^{c_{1}}(1-\mathsf{s})(1-\mathsf{s}^{c_{2}+1})}}{1-\mathsf{s}^{c_{1}+c_{2}+2}}\Bigg)^{2},
    \end{split}
    \end{equation}
where 
\begin{equation} \label{def of mathsf s}
\mathsf{s}:= e^{-\lambda}, \qquad \lambda \in (0,\infty). 
\end{equation}
Then we have the following. 

\begin{thm}[\textbf{$q$-deformed Wachter distribution}] \label{thm:limiting density} 
    Let $c_1,c_2\geq0$ be fixed and let $q=e^{-\lambda/N}$ for a fixed  $\lambda \in (0,\infty)$. Let the parameters $\mathsf{a}$ and $\mathsf{b}$ be scaled as in \eqref{eqn:scale of alphabet}. For $c_2>0$, let $\lambda_{c}$ denote the unique positive solution of
    \begin{equation}\label{eqn:lambdac condition}
        \mathsf{s}(1+\mathsf{s}^{c_1}-\mathsf{s}^{c_1+c_2+1})=1,
    \end{equation}
    where $\mathsf{s}$ is given by \eqref{def of mathsf s}. For $c_2=0$, set $\lambda_c=0$.  
    Then as $N \to \infty$, we have 
    \begin{equation} 
    \mu_N(dx)   \to \rho^{(c_1,c_2)}(x)\,dx \qquad \text{almost surely},
    \end{equation} 
    where 
        \begin{equation}
    \begin{split} \label{def of limiting law}
   \rho^{(c_{1},c_{2})}(x)
        :=\frac{1}{\pi\lambda x}
    \arccos \bigg(  \frac{1+\mathsf{s}^{c_1}-(\mathsf{s}^{-1}+\mathsf{s}^{c_1+c_2+1})x}
    {2\sqrt{\mathsf{s}^{c_1}(1-x)(1-\mathsf{s}^{c_2}x)}} \bigg) \mathbbm{1}_{[a,b]}(x)  + \begin{cases}
            0, & \textup{if }\lambda\leq\lambda_{c},\smallskip\\
            \displaystyle\frac{1}{\lambda x}\mathbbm{1}_{[b,1]}(x), & \textup{if }\lambda>\lambda_{c}.
        \end{cases}
    \end{split}
    \end{equation}
    Here, $a$ and $b$ are given by \eqref{eqn:def of endpoints ab}.   
\end{thm}

\begin{figure}[t]   \centering
    \begin{tabular}{ | c | c | c  | } 
        \hline
        \cellcolor{gray!25}  & \cellcolor{gray!10} \parbox[c]{5.0cm}{\centering \vspace{.5\baselineskip} $ \lambda < \lambda_c $ \\ \emph{band region}   \vspace{.5\baselineskip} }  & \cellcolor{gray!10} \parbox[c]{5.0cm}{\centering \vspace{.5\baselineskip} $ \lambda > \lambda_c $ \\ \emph{band \& saturated region} \vspace{.5\baselineskip}}     
        \\
        \hline
        \cellcolor{gray!10} \parbox[c]{1.5cm}{\centering  $c_1=0$ \\ \emph{hard edge} \vspace{.5\baselineskip} }  
        &   \rule{0pt}{3cm} \raisebox{-0.4\height}{\includegraphics[width=0.4\linewidth]{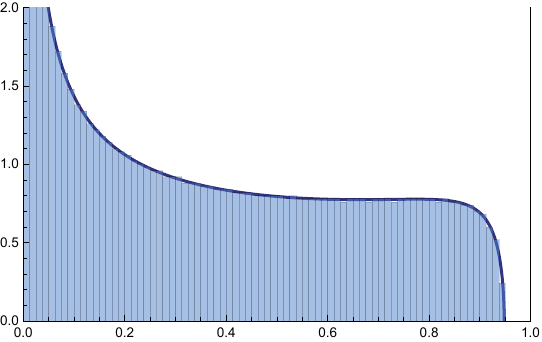}  }   
        &     \rule{0pt}{3cm} \raisebox{-0.4\height}{\includegraphics[width=0.4\linewidth]{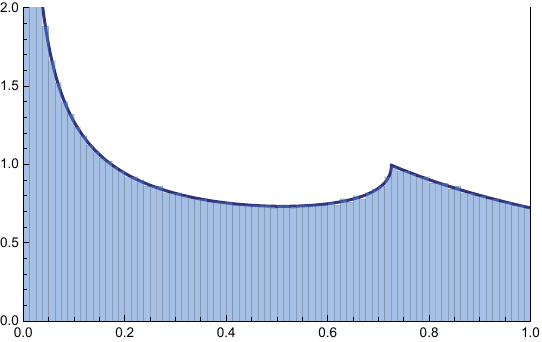} }  
         \\
        \hline
        \cellcolor{gray!10}  \parbox[c]{1.5cm}{\centering \vspace{-3em} $c_1>0$ \\ \emph{soft edge}  }    &    \rule{0pt}{3cm} \raisebox{-0.4\height}{\includegraphics[width=0.4\linewidth]{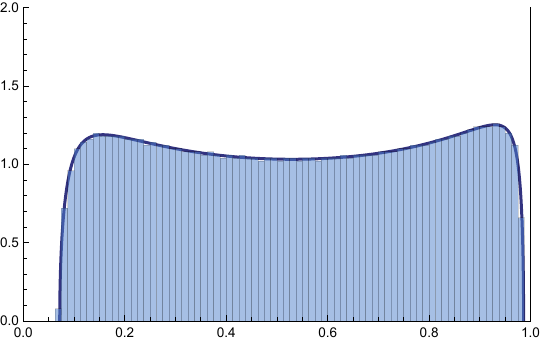}   }     &    \rule{0pt}{3cm} \raisebox{-0.4\height}{\includegraphics[width=0.4\linewidth]{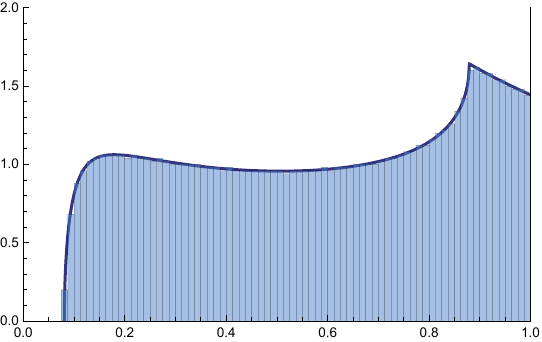} }    
        \\
        \hline
    \end{tabular}
    \caption{Histograms of the zero sets of the corresponding little \(q\)-Jacobi orthogonal polynomials, together with the limiting density. The left column shows the regime \(\lambda<\lambda_c\), where only the band region is present, while the right column shows the regime \(\lambda>\lambda_c\), where a saturated region appears. The top and bottom rows correspond respectively to the hard-edge case \(c_1=0\) and the soft-edge case \(c_1>0\).} 
    \label{Fig_limiting density}
\end{figure}

See Figure~\ref{Fig_limiting density} for a numerical illustration of
Theorem~\ref{thm:limiting density}. Rather than sampling the \(q\)-JUE,
we plot the empirical distribution of the zeros of the corresponding little
\(q\)-Jacobi orthogonal polynomials. This deterministic zero distribution may
be regarded as a low-temperature limit of the ensemble which preserves the
macroscopic density; see \cite{BJM26,BJO25} for related discussions.

\begin{rem}[Band-to-saturated transition] \label{Rem_band to saturated}
A notable feature of the limiting density is the presence of a nontrivial
phase transition, which we refer to as the \emph{band-to-saturated transition}.
Here the saturated region is the part of the support on which the density in \eqref{def of limiting law} attains the upper constraint \(1/(\lambda x)\). In the present setting, this region appears only in the regime \(\lambda>\lambda_c\). Such saturated regions are a familiar feature of discrete orthogonal polynomial ensembles; see \cite{BKMM07,DD22,HMO25}. We also refer to \cite{BJM26} for further discussion and intuition explaining why analogous phase transitions arise in \(q\)-deformed random matrix models.

An important aspect of our result, in line with the series of works
\cite{BFO26,BJO25,BJM26}, is that both the limiting density and the exact
critical transition point can be identified explicitly. Another notable feature is that the transition is absent when \(c_2=0\). 
Indeed, in this case
\eqref{eqn:lambdac condition} reduces to
\[
    \mathsf{s}\bigl(1+\mathsf{s}^{c_1}-\mathsf{s}^{c_1+1}\bigr)
    =
    \mathsf{s}+(1-\mathsf{s})\mathsf{s}^{c_1+1}
    =
    1-(1-\mathsf{s})(1-\mathsf{s}^{c_1+1}).
\]
Under the standing assumptions on the parameters, this quantity is strictly
less than \(1\) for \(0<\mathsf{s}<1\). Hence the only solution in
\(0<\mathsf{s}\leq 1\) is \(\mathsf{s}=1\), equivalently \(\lambda=0\).
Thus there is no positive critical point when \(c_2=0\), and a nontrivial
band-to-saturated transition can occur only for \(c_2>0\).  
\end{rem}

\begin{rem}[Degenerations to the limiting laws of the JUE and the Meixner ensemble]
The limiting law \eqref{def of limiting law} recovers the macroscopic
limits of both the classical JUE and the Meixner ensemble under two
different limiting procedures discussed above.

We first recall that the limiting spectral distribution of the
classical JUE is the Wachter distribution \cite{Wac80}
  \begin{equation}\label{def of Wachter distribution}
 \rho_{\rm J}(x)  =  (c_1 + c_2+2) \frac{ \sqrt{(x - a_{ \rm J } ) ( b_{\rm J} - x)}}{ 2 \pi x (1 - x)}
  \end{equation}   
  supported on $(a_{ \rm J },b_{ \rm J })$ with these endpoints specified by
   \begin{equation} \label{def of Wachter endpts}
 \frac{1}{(c_1 + c_2 + 2)^2} \Big ( \sqrt{(c_1 + 1) (c_1 + c_2 + 1)} \pm \sqrt{c_2 + 1} \Big )^2.
  \end{equation}  
This distribution is recovered from
\eqref{def of limiting law} by keeping $c_1$ and $c_2$ fixed and
letting $\lambda\to0$. Indeed, the endpoints $a$ and $b$ defined in
\eqref{eqn:def of endpoints ab} satisfy
\begin{equation}
    \lim_{\lambda\to0}a=a_{\rm J},
    \qquad
    \lim_{\lambda\to0}b=b_{\rm J},
\end{equation}
and the density converges pointwise on the interior of the limiting
support:
\begin{equation}
    \lim_{\lambda\to0}
    \rho^{(c_1,c_2)}(x)
    =
    \rho_{\rm J}(x),
    \qquad
    x\in(a_{\rm J},b_{\rm J}).
\end{equation}

The limiting distribution of the Meixner ensemble is obtained through
a different scaling. Let $\alpha\in(0,1)$ and $Q\geq0$, and define
\begin{equation}\label{eqn:Meixner limiting endpoints}
    a_{\rm M}
    =
    \frac{
    \bigl(1-\alpha\sqrt{Q+1}\bigr)^2
    }{
    1-\alpha^2
    },
    \qquad
    b_{\rm M}
    =
    \frac{
    \bigl(1+\alpha\sqrt{Q+1}\bigr)^2
    }{
    1-\alpha^2
    }.
\end{equation}
For $a_{\rm M}<u<b_{\rm M}$, set
\begin{equation}
    \begin{split}
        \vartheta_0(u)
   =
        \frac{1}{\pi}
        \arctan
        \sqrt{
        \frac{
        a_{\rm M}(b_{\rm M}-u)
        }{
        b_{\rm M}(u-a_{\rm M})
        }
        },
        \qquad 
        \vartheta_Q(u)
        =
        \frac{1}{\pi}
        \arctan
        \sqrt{
        \frac{
        (a_{\rm M}+Q)(b_{\rm M}-u)
        }{
        (b_{\rm M}+Q)(u-a_{\rm M})
        }
        }.
    \end{split}
\end{equation}
The associated limiting density is then
\begin{equation}\label{eqn:Meixner limiting density}
    \rho_{\rm M}(u)
    =
    \begin{cases}
        \displaystyle
        \mathbbm 1_{[0,a_{\rm M}]}(u)
        +
        \bigl(
        \vartheta_0(u)+\vartheta_Q(u)
        \bigr)
        \mathbbm 1_{(a_{\rm M},b_{\rm M})}(u),
        &\textup{if }
        \alpha^2(Q+1)<1,
        \smallskip
        \\
        \displaystyle
        \bigl(
        -\vartheta_0(u)+\vartheta_Q(u)
        \bigr)
        \mathbbm 1_{(a_{\rm M},b_{\rm M})}(u),
        &\textup{if }
        \alpha^2(Q+1)> 1.
    \end{cases}
\end{equation}
The critical case $\alpha^2(Q+1)=1$, for which $a_{\rm M}=0$, is obtained by continuity. 
This limiting density was obtained by Johansson \cite[Eqs.~(2.19) and (2.20), Section~6]{Joh00}; see also \cite[Propositions~2.1 and 2.2]{CP15}.

To recover \eqref{eqn:Meixner limiting density} from the $q$-Wachter
distribution, we pass from the exponential variable $x$ to the
macroscopic lattice variable of the Meixner ensemble. More precisely,
set
\begin{equation}
    c_1=-\frac{\log\alpha^2}{\lambda},
    \qquad
    c_2=Q,
\end{equation}
and introduce the change of variables
\begin{equation}
    u=\frac{1}{\lambda}\log \frac{1}{x},
    \qquad
    x=e^{-\lambda u}.
\end{equation}
Let $a_\lambda$ and $b_\lambda$ denote the endpoints in
\eqref{eqn:def of endpoints ab} under this scaling. Then
\begin{equation}
   \lim_{ \lambda \to 0 } \frac{1}{\lambda}\log \frac{1}{b_\lambda}
    = a_{\rm M},
    \qquad
    \lim_{ \lambda \to 0 } \frac{1}{\lambda}\log \frac{1}{a_\lambda}
    = b_{\rm M}.
\end{equation}
The reversal of the endpoints is due to the fact that the map
$u\mapsto e^{-\lambda u}$ is decreasing. 
At the level of densities, we define
\begin{equation}
    \rho_{{\rm M},\lambda}(u)
    :=
    \lambda e^{-\lambda u}
    \rho^{(-\log\alpha^2/\lambda,Q)}(e^{-\lambda u}).
\end{equation}
Then it follows from straightforward computations that 
\begin{equation}
    \lim_{\lambda\to0}
    \rho_{{\rm M},\lambda}(u)
    =
    \rho_{\rm M}(u),
    \qquad
    u\in(a_{\rm M},b_{\rm M}).
\end{equation} 
\end{rem}

\begin{rem}[Degenerations to the $q$-deformed Marchenko--Pastur and semicircle laws]
The $q$-Wachter distribution also contains the previously known
$q$-deformed Marchenko--Pastur \cite{BJM26} and semicircle laws \cite{BFO26} as natural degenerations. Fix $\lambda>0$ and let
$\rho_{\rm qMP}^{(c)}$ denote the $q$-deformed Marchenko--Pastur
density obtained in \cite[Eq.~(1.5) and Theorem~1.2]{BJM26}. If $c_1=c$ and $c_2 \to \infty$, then $\mathsf s^{c_2}\to0$, and the endpoints and the density in \eqref{eqn:def of endpoints ab} and \eqref{def of limiting law} reduce to their $q$-deformed Marchenko--Pastur counterparts. In
particular,
\begin{equation}
    \lim_{c_2\to\infty}
    \rho^{(c,c_2)}(x)
    =
    \rho_{\rm qMP}^{(c)}(x)
\end{equation}
at every continuity point of the limiting density. The corresponding
critical equation reduces to $ \mathsf s\bigl(1+\mathsf s^c\bigr)=1,$ which is precisely the phase transition condition for the
$q$-deformed Marchenko--Pastur law. 
The $q$-deformed semicircle law in \cite[Theorem~2.6]{BFO26} is recovered from the special case $c=0$ via the quadratic transform; see \cite[Remark~3]{BJM26}. 
\end{rem}

Beyond the global limiting distribution, we study the local scaling
limits of the $q$-JUE, with particular emphasis on the fluctuation of its
smallest particle. This follows the same philosophy as Johansson's
analysis in \cite{Joh00}, where the fluctuation of the passage time in
homogeneous LPP was obtained through the largest particle distribution
of the Meixner ensemble.

The Riemann--Hilbert theory of discrete orthogonal polynomials
\cite{BKMM07} treats discrete orthogonal polynomial ensembles on finite
sets of nodes whose empirical distributions converge to a positive
analytic density on a compact interval, together with analytic varying
weights. In this setting, one obtains universal limits such as the
discrete sine kernel in bands and the Airy kernel at generic band edges.
Nevertheless, the hypotheses of this theory do not apply directly to the
present $q$-JUE, which is supported on a one-sided geometric $q$-lattice;
consequently, a separate asymptotic analysis is required. 

Motivated by this general picture, we derive the relevant scaling limits
for the $q$-JUE. As in classical random matrix theory, the local
behaviour near the smallest particle depends on the nature of the edge: one must distinguish between the soft edge and hard edge regimes, according to whether $c_1>0$ or $c_1=0$.  

To identify the correct edge scaling, we first need the asymptotic
behaviour of the limiting density near the relevant edge.
A hard edge occurs at $x=0$ when $c_1=0$, that is, when
$\mathsf{a}=d_1$ remains of order one. In this case, the
$q$-deformed Wachter density has the inverse square-root behaviour
\begin{equation}\label{def of kappa}
\rho^{(0,c_2)}(x)
=
\frac{\kappa}{2\pi\sqrt{x}}\bigl(1+o(1)\bigr),
\qquad \kappa
:= \frac{2}{\lambda}
\sqrt{
\frac{(1-\mathsf{s})(1-\mathsf{s}^{c_2+1})}{\mathsf{s}}
} 
\end{equation}
as $x \downarrow 0$.  
When $c_1>0$, the left endpoint $a$ defined in
\eqref{eqn:def of endpoints ab} is a soft edge. In this case, the
$q$-deformed Wachter density exhibits the square-root behaviour
\begin{equation} 
\rho^{(c_1,c_2)}(x)
=
\frac{\theta}{\pi}(x-a)^{1/2}\bigl(1+o(1)\bigr),
\qquad x\downarrow a,
\end{equation}
where
\begin{equation}   \label{def of theta}
\theta
:= \begin{cases}
\displaystyle
\frac{
\sqrt{(1-a)(1-\mathsf{s}^{c_2}b)}
-
\sqrt{(1-b)(1-\mathsf{s}^{c_2}a)}
}{
\lambda a
\sqrt{(1-a)(1-\mathsf{s}^{c_2}a)(b-a)}
},
& \textup{if } \lambda\le \lambda_c,
\smallskip \\
\displaystyle
\frac{
\sqrt{(1-a)(1-\mathsf{s}^{c_2}b)}
+
\sqrt{(1-b)(1-\mathsf{s}^{c_2}a)}
}{
\lambda a
\sqrt{(1-a)(1-\mathsf{s}^{c_2}a)(b-a)}
},
& \textup{if } \lambda>\lambda_c.
\end{cases}
\end{equation} 
Let $x_{\min}$ be the smallest particle of the $q$-JUE. 

\begin{thm}[\textbf{Smallest particle distribution of the $q$-JUE}]
\label{Thm_qJUE scaling limits}
Let the parameters $\mathsf{a}$ and $\mathsf{b}$ be scaled as in
\eqref{eqn:scale of alphabet}. Fix $c_2\ge0$, and let
$q=e^{-\lambda/N}$ for some fixed $\lambda\in(0,\infty)$. Then the
following edge limits hold.

\begin{enumerate}
    \item[\textup{(i)}] \textbf{\textup{(Hard edge)}}
    Suppose that $c_1=0$. Then, for $s>0$,
    \begin{equation} \label{hard edge scaling limit main}
    \lim_{N\to\infty}
    \mathbb{P}\Big(
    x_{\min}>
    \frac{s}{\kappa^2N^2}
    \Big)
    =
    \det\bigl(I-K_{d_1}^{\rm Bessel}\bigr)_{L^2(0,s)},
    \end{equation}
    where $K_{d_1}^{\rm Bessel}$ is the standard Bessel kernel
    \eqref{def of K Bessel standard}, and $\kappa$ is the hard edge
    constant defined in \eqref{def of kappa}.

    \smallskip

    \item[\textup{(ii)}] \textbf{\textup{(Soft edge)}}
    Suppose that $c_1>0$. Then, for $s\in\mathbb{R}$,
    \begin{equation}  \label{soft edge scaling limit main}
    \lim_{N\to\infty}
    \mathbb{P}\Big(
    x_{\min}\ge
    a-\frac{s}{(\theta N)^{2/3}}
    \Big) = F_{\rm TW}(s),
    \end{equation}
    where $F_{\rm TW}$ is the GUE Tracy--Widom distribution
    \eqref{def of TW}, $a$ is the left endpoint defined in
    \eqref{eqn:def of endpoints ab}, and $\theta$ is the soft edge
    constant given in \eqref{def of theta}.  
\end{enumerate}
For $c_2>0$, both limits hold uniformly for $d_2$ in compact subsets of $\mathbb{R}$.
\end{thm}

Theorem~\ref{Thm_qJUE scaling limits} shows that, although the global limiting distribution of the $q$-JUE is a genuinely $q$-deformed object, its edge fluctuations are governed by the same universal kernels as in classical random matrix theory. In particular, the nature of the left edge is determined by the parameter $c_1$: when $c_1=0$, the limiting support reaches the hard wall at the origin, giving rise to the Bessel kernel; when $c_1>0$, the support detaches from the origin, and the Airy kernel, equivalently the Tracy--Widom distribution, emerges.

In fact, edge scaling limits for the $q$-Gaussian and $q$-Laguerre unitary ensembles have not previously appeared in the literature. In this respect,
Theorem~\ref{Thm_qJUE scaling limits} goes beyond the earlier works
\cite{BJO25,BFO26,BJM26}, which were concerned primarily with global
limiting distributions, by also establishing the local edge statistics.

Through the duality with inhomogeneous LPP established in the following section, the two edge regimes in Theorem~\ref{Thm_qJUE scaling limits} correspond to the strong and weak inhomogeneity regimes discussed in the previous section. Thus, the $q$-JUE provides a unified random matrix framework in which the two fluctuation behaviours of inhomogeneous LPP are realised as the hard and soft edge limits of the smallest particle.

\section{Duality and its consequences} \label{Section_duality}

This section is devoted to the duality relations. 
In Proposition~\ref{Prop_duality Gnm q-JUE}, we establish an exact duality
between inhomogeneous geometric LPP and the $q$-deformed JUE. 
In Corollary~\ref{Cor_duality Hnm JUE}, we state the analogous duality
between inhomogeneous exponential LPP and the classical JUE.  
These dualities provide the bridge that allows us to transfer asymptotic
results for the corresponding ensembles to the LPP setting. 
In particular, assuming the $q$-JUE asymptotics in
Theorems~\ref{thm:limiting density} and~\ref{Thm_qJUE scaling limits}, we
derive the LPP results stated in
Theorems~\ref{Thm_LPP geometric weak}, \ref{Thm_LPP geometric strong},
and~\ref{Thm_exponential}.

\subsection{Duality}

The duality relation, which is the main structural ingredient of the present
paper, is the following. We use the convention $\mathbb{N}=\mathbb{Z}_{\geq 0}$. 

\begin{prop}[\textbf{Duality between inhomogeneous geometric LPP and the $q$-JUE}]\label{Prop_duality Gnm q-JUE}
    Assume that $0<q<1$ and $n \geq m$, and let $G_{n,m}$ denote the last passage time with weights given by~\eqref{def of xi yj weights}. Then, for each $\ell \geq 0$, we have
    \begin{equation} \label{eqn for duality v0}
        \mathbb{P}(G_{n,m} \leq \ell)
        =\frac{1}{\mathcal{C}_{q}(n,m;\alpha)}
        \sum_{\substack{h \in \mathbb{N}^{m}\\\max h_i \leq m+\ell-1}}\prod_{1 \leq i < j \leq m}(q^{h_j}-q^{h_i})^{2}\prod_{j=1}^{m}\frac{(q;q)_{h_j+n-m}}{(q;q)_{h_j}}\alpha^{2h_j},
    \end{equation}
    where $\mathcal{C}_{q}(n,m;\alpha)$ is the normalising constant.
    
    Equivalently, let $\{x_j\}_{j=1}^{m}$ be distributed according to the $q$-JUE in~\eqref{eqn:jpdf of qJUE}, with parameters
    \begin{equation} \label{matching of a b qJUE parameter}
        \mathsf{a}= \frac{2\log \alpha}{\log q}-1,
\qquad
\mathsf{b}= n-m .
\end{equation}
Then we have 
\begin{equation} \label{eqn for duality v1}
\mathbb{P}(G_{n,m} \leq \ell)
=
\mathbb{P}\Big( \min_{1 \leq j \leq m} x_j
\geq q^{m+\ell-1} \Big).
\end{equation}
\end{prop}

\begin{rem}[Homogeneous limit and the Meixner ensemble]
The duality in Proposition~\ref{Prop_duality Gnm q-JUE} can be viewed as a $q$-deformation of Johansson's
duality in the homogeneous case. More precisely, in \cite[Proposition~1.3]{Joh00}, it was shown that, when $q=1$,
\begin{equation} \label{Johansson Meixner duality}
\mathbb{P}(G_{n,m} \leq \ell)\Big|_{q=1}
\propto
\sum_{\substack{h \in \mathbb{N}^{m} \\
\max h_i \leq m+\ell-1}}
\prod_{1 \leq i < j \leq m} (h_j-h_i)^2
\prod_{j=1}^{m}
\binom{h_j+n-m}{h_j}\alpha^{2h_j}.
\end{equation}
The right-hand side is the distribution function of the largest particle in
the Meixner ensemble.
As explained in Section~\ref{Section_qJUE}, the Meixner ensemble arises as
a degeneration of the $q$-JUE, and accordingly
\eqref{Johansson Meixner duality} is recovered as a limiting case of
\eqref{eqn for duality v0}. 
\end{rem} 

Let us mention a related appearance of the $q$-Jacobi structure in \cite{BO24}. For the probability measure on plane partitions considered in \cite[Eq.~(2.1)]{BO24}, the associated point process is identified in \cite[Theorem~2.1]{BO24} with a discrete Muttalib--Borodin ensemble.
The specialisation in \cite[Remark~2.3]{BO24} yields a little $q$-Jacobi ensemble, while their sampling construction connects the underlying partition to last-passage percolation; see \cite[Theorem~5.4]{BO24}.

The correspondence in Proposition~\ref{Prop_duality Gnm q-JUE} can also be obtained from \cite[Proposition~2]{FR05} by sending the parameter of the additional column to zero. We give a direct, self-contained derivation below, including the explicit identification of the $q$-JUE parameters.

\begin{proof}[Proof of Proposition~\ref{Prop_duality Gnm q-JUE}] 
We first fix the notation for Schur polynomials with a prescribed number of
variables. 
A partition is a weakly decreasing sequence
$\lambda=(\lambda_1,\lambda_2,\dots)$ of non-negative integers with only
finitely many non-zero parts. Its length is denoted by $l(\lambda):=\#\{j:\lambda_j>0\}.$ 
For an $r$-tuple $\bm z=(z_1,\dots,z_r)$, we write
\[
    s_\lambda^{(r)}(\bm z)
    =
    s_\lambda(z_1,\dots,z_r)
\]
for the Schur polynomial in $r$ variables. 
If $l(\lambda)\le r$, it is given by
\begin{equation} \label{def of Schur polynomial}
    s_\lambda^{(r)}(\bm z)
    =
    \frac{1}{\Delta(\bm z)}
    \det\Big\{z_i^{\lambda_j+r-j}\Big\}_{i,j=1}^{r},
    \qquad
    \Delta(\bm z)
    =
    \prod_{1\le i<j\le r}(z_i-z_j).
\end{equation}
We use the convention that
\begin{equation}
 s_\lambda^{(r)}(\bm z)=0
    \qquad \text{if } l(\lambda)>r .
\end{equation} 
When the number of variables is clear from the context, we suppress the
superscript and simply write $s_\lambda(\bm z)$.

Let $  \bm x^{(n)}=(x_1,\dots,x_n) $ and $ \bm y^{(m)}=(y_1,\dots,y_m)$ be two finite alphabets satisfying $0\le x_i y_j<1$. 
By the Robinson--Schensted--Knuth correspondence, the inhomogeneous
geometric LPP satisfies
\begin{equation} \label{schur-meas}
    \mathbb{P}(G_{n,m}\le \ell)
    =
    \prod_{i=1}^{n}\prod_{j=1}^{m}(1-x_i y_j)
    \sum_{\lambda:\lambda_1\le \ell}
    s_\lambda^{(n)}(\bm x^{(n)})
    s_\lambda^{(m)}(\bm y^{(m)}),
\end{equation}
where the sum is over all partitions $\lambda$ with $\lambda_1\le \ell$;
see e.g. \cite[Section 6]{BBP05} and \cite[Section 2]{FR05}. 
By the above convention, the summand vanishes unless $    l(\lambda)\le \min\{n,m\}.$  
In particular, if $n\ge m$, then only partitions with $l(\lambda)\le m$
contribute. 
We also emphasise that the parameters $x_i$ and $y_j$ in
\eqref{schur-meas} are general admissible parameters; they need not yet be
specialised as in \eqref{def of xi yj weights}.

From now on, assume that $n\ge m$. 
We first consider the principal specialisation
\[
    \bm x^{(n)}=(1,q,\dots,q^{n-1}).
\]
For any partition $\lambda$ contributing to \eqref{schur-meas}, we have
$l(\lambda)\le m\le n$, and hence \eqref{def of Schur polynomial} gives
\begin{align} \label{Schur q exponential}
\begin{split}
    s_{\lambda}^{(n)}(1,q,\dots, q^{n-1})& = \frac{1}{\Delta(1,q,\dots, q^{n-1})} \det\Big\{(q^{\lambda_{j}+n-j})^{i-1}\Big\}_{i,j=1}^{n}
    \\
    &=  \frac{1}{\Delta(1,q,\dots, q^{n-1})}  \prod_{1 \leq i < j \leq n}(q^{\lambda_{j}+n-j}-q^{\lambda_{i}+n-i}). 
\end{split}
\end{align} 
Here the last equality follows from the Vandermonde determinant identity.

Since only partitions with $l(\lambda)\le m$ contribute, we have $ \lambda_j=0$ for all $j>m$.  
Therefore the product in \eqref{Schur q exponential} can be separated into
the part involving the first $m$ indices and the remaining factors: 
\begin{equation} \label{product decomposition q Schur}
    \prod_{1 \leq i < j \leq n}(q^{\lambda_{j}+n-j}-q^{\lambda_{i}+n-i}) =  \prod_{1 \leq i < j \leq m}(q^{\lambda_{j}+n-j}-q^{\lambda_{i}+n-i})\prod_{j=m+1}^{n}\prod_{i=1}^{j-1}(q^{ n-j}-q^{\lambda_{i}+n-i}),
\end{equation}
where we used $\lambda_j=0$ for all $j>m$.

We next rewrite the second product in terms of the $q$-Pochhammer symbols defined in \eqref{q-poch}. In the remaining steps of the proof, we shall write $C$ for a positive constant which is
independent of $\lambda$ and $\ell$. 
It may depend on the fixed parameters $q,n,m$, and, after the
specialisation below, also on $\alpha$. 
By a slight abuse of notation, the value of $C$ may change from one
occurrence to the next.
For each fixed $i=1,\dots,m$, we have
\begin{equation}
\begin{split}
    \prod_{j=m+1}^{n}
    \bigl(q^{n-j}-q^{\lambda_i+n-i}\bigr)
    &=
    q^{\sum_{j=m+1}^{n}(n-j)}
    \prod_{j=m+1}^{n}
    \bigl(1-q^{\lambda_i+j-i}\bigr)   =
    q^{\binom{n-m}{2}}
    \frac{(q;q)_{\lambda_i+n-i}}
         {(q;q)_{\lambda_i+m-i}} .
\end{split}
\end{equation}
The remaining factors in the second product of
\eqref{product decomposition q Schur}, namely those with $i>m$, are
independent of $\lambda$. Hence
\begin{equation} \label{tail product up to constant}
    \prod_{j=m+1}^{n}\prod_{i=1}^{j-1}
    \bigl(q^{n-j}-q^{\lambda_i+n-i}\bigr)
    =
    C
    \prod_{i=1}^{m}
    \frac{(q;q)_{\lambda_i+n-i}}
         {(q;q)_{\lambda_i+m-i}} .
\end{equation}
Moreover,
\begin{equation} \label{first product up to constant}
    \prod_{1\le i<j\le m}
    \bigl(q^{\lambda_j+n-j}-q^{\lambda_i+n-i}\bigr)
    =
    C
    \prod_{1\le i<j\le m}
    \bigl(q^{\lambda_j+m-j}-q^{\lambda_i+m-i}\bigr).
\end{equation}
Combining \eqref{Schur q exponential},
\eqref{product decomposition q Schur},
\eqref{tail product up to constant}, and
\eqref{first product up to constant} with the analogous formula for
$s_\lambda^{(m)}(1,q,\dots,q^{m-1})$, we obtain
\begin{equation} \label{Schur product q-JUE form}
\begin{split}
    &s_{\lambda}^{(n)}(1,q,\dots,q^{n-1})
    s_{\lambda}^{(m)}(1,q,\dots,q^{m-1})   =
    C
    \prod_{1\le i<j\le m}
    \bigl(q^{\lambda_j+m-j}-q^{\lambda_i+m-i}\bigr)^2
    \prod_{j=1}^{m}
    \frac{(q;q)_{\lambda_j+n-j}}
         {(q;q)_{\lambda_j+m-j}} .
\end{split}
\end{equation}

Now we put the parameters $x_i=\alpha q^{i-1}, y_j=\alpha q^{j-1}$ as in \eqref{def of xi yj weights}. 
Using the homogeneity relation 
\begin{equation}
s_{\lambda}(\alpha\bm x) = \alpha^{ |\lambda| } s_{\lambda}(\bm x)
\end{equation}
and absorbing the prefactor $\prod_{i=1}^{n}\prod_{j=1}^{m}
    (1-\alpha^2q^{i+j-2})$  in \eqref{schur-meas} into the constant $C$, we obtain  
\begin{equation} \label{LPP Schur sum before h}
\begin{split}
    \mathbb{P}(G_{n,m}\le \ell)
    &=
    C
    \sum_{\substack{\lambda:\lambda_1\le \ell\\ l(\lambda)\le m}}
    \alpha^{2|\lambda|}
    \prod_{1\le i<j\le m}
    \bigl(q^{\lambda_j+m-j}-q^{\lambda_i+m-i}\bigr)^2 
    \prod_{j=1}^{m}
    \frac{(q;q)_{\lambda_j+n-j}}
         {(q;q)_{\lambda_j+m-j}} .
\end{split}
\end{equation}

We now make the change of variables
\[
    h_j=\lambda_j+m-j,
    \qquad j=1,\dots,m.
\]
Then \eqref{LPP Schur sum before h} becomes 
\begin{equation} \label{unordered h sum}
\begin{split}
    \mathbb{P}(G_{n,m}\le \ell)
    &=
    C
    \sum_{\substack{h \in \mathbb{N}^{m} \\
    \max h_i \leq m+\ell-1}}
    \prod_{1 \leq i < j \leq m}
    (q^{h_j}-q^{h_i})^{2}
    \prod_{j=1}^{m}
    \frac{(q;q)_{h_j+n-m}}{(q;q)_{h_j}}
    \alpha^{2h_j}.
\end{split}
\end{equation}
Letting $\ell\to\infty$ and using $\mathbb P(G_{n,m}<\infty)=1$ fixes the normalisation constant and completes the proof of \eqref{eqn for duality v0}.

\medskip

It remains to identify the right-hand side of
\eqref{eqn for duality v0} with the $q$-JUE distribution.
Recall that the $q$-JUE is defined as a density with respect to the
Jackson measure \eqref{def of Jackson measure} on the $q$-lattice.
Thus, when the particles are written as \(x_j=q^{h_j}\), the
corresponding probability mass function with respect to the counting
measure on \(\mathbb Z_{\geq 0}^m\) is proportional to
\[
    \prod_{1\le i<j\le m}(q^{h_i}-q^{h_j})^2
    \prod_{j=1}^m q^{h_j}
    w^{(\mathsf a,\mathsf b)}(q^{h_j};q).
\]

For $x=q^h$, $h\in\mathbb N$, the $q$-Jacobi weight in
\eqref{eqn:def of qjacobi weight} can be written in finite
$q$-Pochhammer form. Indeed, using
\[
    \frac{(a;q)_\infty}{(aq^r;q)_\infty}
    =
    (a;q)_r,
    \qquad r\in \mathbb Z_{\geq 0},
\]
the corresponding one-particle mass at \(x=q^h\) is
\begin{equation} \label{def of qJUE weight v2}
   q^h w^{(\mathsf{a},\mathsf{b})}(q^h;q)
    =
    q^{(\mathsf{a}+1)h}(q^{h+1};q)_{\mathsf{b}} .
\end{equation}  
With the parameter matching  \eqref{matching of a b qJUE parameter}, we have
\[
    q^{(\mathsf{a}+1)h}=\alpha^{2h},
    \qquad
    (q^{h+1};q)_{n-m}
    =
    \frac{(q;q)_{h+n-m}}{(q;q)_h}.
\]
Thus the one-particle weight in \eqref{eqn for duality v0} is precisely the
$q$-JUE weight with parameters \eqref{matching of a b qJUE parameter}. 
Since $0<q<1$, writing $x_j=q^{h_j}$, the identity \eqref{eqn for duality v0} is therefore equivalent to \eqref{eqn for duality v1}. 
This completes the proof. 
\end{proof}

We next record the continuum counterpart of Proposition~\ref{Prop_duality Gnm q-JUE}. 
In the homogeneous case, Johansson identified exponential LPP with the Wishart (Laguerre) ensemble whose asymptotic analysis is substantially simpler than that of the discrete Meixner ensemble associated with geometric LPP. 
A similar phenomenon occurs in the present inhomogeneous setting. While the geometric model is naturally related to the $q$-deformed JUE, the exponential model is instead dual to the classical JUE. As a consequence, many of the asymptotic results can be deduced directly from the well-developed theory of classical random matrix ensembles.
 
We denote the smallest eigenvalue of the JUE defined in
\eqref{def of classical JUE} by $  y_{\min}^{(N)}(\mathsf{a},\mathsf{b}) :=
    \min_{1\le j\le N} y_j$.  
The following is the exponential analogue of the duality in
Proposition~\ref{Prop_duality Gnm q-JUE}.

\begin{cor}[\textbf{Duality between inhomogeneous exponential LPP and the JUE}]  \label{Cor_duality Hnm JUE} Assume that $n\ge m$. 
Let $\eta,\lambda>0$ be fixed, and let $H_{n,m}$ denote the last passage
time defined in \eqref{def of Hnm exponential LPP time}, where the
parameters $t_k$ are given by \eqref{def of tk exponential}. 
Let $\beta$ be the parameter determined by this choice of $t_k$. 
Then, for each $\ell\ge 0$, we have 
\begin{equation} \label{duality Hnm and JUE smallest eigenvalue}
 \mathbb{P}( H_{n,m} \leq \ell) =  \mathbb{P}(y^{(m)}_{\mathrm{min}}(2\beta+1,n-m) \geq e^{-\lambda \ell}).   
\end{equation} 
\end{cor}

\begin{proof}
Let $q=e^{ -\epsilon \lambda }$ and $\alpha = e^{ -\epsilon \eta }$ for some small $\epsilon>0$. Recall that $x_i$ and $y_j$ are given by \eqref{def of xi yj weights}. Then as $\epsilon \to 0$, we have 
\[
x_i y_j = \alpha^2 q^{i+j-2}
= \exp\!\Big(- \epsilon (2\eta + \lambda(i+j-2) )\Big) = 1- \epsilon (2\eta + \lambda(i+j-2) ) +O( \epsilon^2 ) .
\]
Recalling that $\omega_{ij}$ is geometric with parameter $x_i y_j$, it follows that, as $\epsilon \to 0$, the rescaled random variable $\epsilon \omega_{ij}$ converges in distribution to an exponential random variable with rate $t_i + t_j$. 

Note that 
\begin{align*}
\frac{(q;q)_{h_j+n-m}}{(q;q)_{ h_j }} &= \prod_{k=1}^{n-m} (1 - q^{h_j+k}).
\end{align*}
Since $G_{n,m}$ is a continuous function of the weight $w_{i,j}$, by Proposition~\ref{Prop_duality Gnm q-JUE}, we have 
\begin{align*}
&\quad \mathbb{P}( H_{n,m} \leq \ell)  = \lim_{ \epsilon \to 0 } \mathbb{P}( G_{n,m} \le \ell/ \epsilon ) \bigg|_{ q=e^{ -\epsilon \lambda },  \alpha = e^{ -\epsilon \eta } }
\\
& = \lim_{ \epsilon \to 0 } \bigg[ \frac{1}{\mathcal{C}_{q}(n,m;\alpha)}\sum_{\substack{h \in \mathbb{N}^{m}\\\mathrm{max}\{h_{i}\}\leq m+\ell/\epsilon-1}}\prod_{1 \leq i < j \leq m}(q^{h_{j}}-q^{h_{i}})^{2}\prod_{j=1}^{m}\frac{(q;q)_{h_{j}+n-m}}{(q;q)_{h_{j}}}\alpha^{2h_j} \bigg] \bigg|_{ q=e^{ -\epsilon \lambda },  \alpha = e^{ -\epsilon \eta } }
\\
& = \lim_{ \epsilon \to 0 } \bigg[ \frac{1}{\mathcal{C}_{q}(n,m;\alpha)}\sum_{\substack{h \in \mathbb{N}^{m}\\\mathrm{max}\{h_{i}\}\leq m+\ell/\epsilon-1}}\prod_{1 \leq i < j \leq m}(e^{ -\epsilon \lambda h_{j}}-e^{-\epsilon \lambda h_{i}})^{2}\prod_{j=1}^{m} \bigg(   \prod_{k=1}^{n-m} (1 - e^{-\epsilon \lambda(h_j+k)} ) \bigg)  e^{-2\epsilon \eta h_j} \bigg]  
\\
&= \lim_{ \epsilon \to 0 } \bigg[ \frac{1}{\mathcal{C}_{q}(n,m;\alpha)}\sum_{\substack{ x \in  \epsilon \mathbb{N}^{m}\\\mathrm{max}\{x_{i}\}\leq \ell+ \epsilon (m-1)} }\prod_{1 \leq i < j \leq m}(e^{ - \lambda x_{j}}-e^{- \lambda x_{i}})^{2}\prod_{j=1}^{m} \bigg(   \prod_{k=1}^{n-m} (1 - e^{- \lambda x_j-\epsilon \lambda k } ) \bigg)  e^{-2 \eta x_j} \bigg]   ,
\end{align*} 
where $x_j=\epsilon h_j$. 
Therefore we obtain   
\begin{equation}
\mathbb{P}( H_{n,m} \leq \ell) \propto \int_{ [0,\ell]^m }\prod_{1 \leq i < j \leq m}  (e^{ -\lambda x_j }- e^{ -\lambda x_i })^2 \prod_{j=1}^m (1-e^{ -\lambda x_j })^{n-m} e^{ -2\eta x_j }\,dx_j, 
\end{equation} 
where the multiplicative normalisation constant is independent of $\ell$. 
Furthermore, by making the change of variables $y_j=e^{-\lambda x_j}$, since  
\begin{equation} \label{change of var x y exp}
x_j=-\frac{1}{\lambda}\log y_j,
\qquad
dx_j=-\frac{1}{\lambda y_j}\,dy_j, 
\end{equation} 
we have 
\begin{equation}
\mathbb{P}( H_{n,m} \leq \ell) \propto \int_{ [e^{-\lambda \ell},1]^m }
\prod_{1 \le i<j\le m}(y_j-y_i)^2
\prod_{j=1}^m y_j^{2\beta+1}(1-y_j)^{n-m}  \,dy_j.   
\end{equation} 
Thus we conclude \eqref{duality Hnm and JUE smallest eigenvalue}, where the normalisation constant is fixed by letting $\ell \to \infty$. 
\end{proof}

\begin{rem}[Extended Wishart correspondence] \label{Rem_BP extended Wishart}
A related matrix model appears in \cite{BP08}, where inhomogeneous
exponential LPP is connected to the so-called extended Wishart ensemble.
Let $X=(X_{ij})$ be an $N\times p$ random matrix with independent complex
Gaussian entries satisfying
\begin{equation} \label{def of Gaussian weight exponential inhomo}
    X_{ij}\sim \textup{N}_{\mathbb C}
    \Bigl(0,\frac{1}{t_i+t_j}\Bigr).
\end{equation}
Thus $X$ may be viewed as an inhomogeneous Ginibre-type matrix: the entries
remain independent, but their variances depend on the indices. 
In the homogeneous case, when $t_i+t_j$ is constant, this reduces, after a
deterministic rescaling, to the rectangular complex Ginibre ensemble \cite{BF25}.
The matrix $XX^*$ is then called the extended Wishart matrix.

It was shown in \cite[Theorem~1]{BP08} that the last passage time is related
to the largest eigenvalue of this matrix. More precisely, if
$\widetilde y_{\max}$ denotes the largest eigenvalue of $XX^*$, then
\begin{equation}
    \mathbb{P}(H_{p,N}\le x)
    =
    \mathbb{P}(\widetilde y_{\max}\le x).
\end{equation}
The change of variables used in \eqref{change of var x y exp} transforms
this extended Wishart description into the Jacobi form appearing in
Corollary~\ref{Cor_duality Hnm JUE}. 
Thus Corollary~\ref{Cor_duality Hnm JUE} can also be recovered from \cite{BP08}. 
\end{rem}

\subsection{Proof of Theorems~\ref{Thm_LPP geometric weak}, ~\ref{Thm_LPP geometric strong} and ~\ref{Thm_exponential}}

In this subsection, we derive the LPP results from the duality relations established in the previous subsection, namely Proposition~\ref{Prop_duality Gnm q-JUE} and Corollary~\ref{Cor_duality Hnm JUE}.
For this purpose, we take the $q$-JUE asymptotic results (Theorems~\ref{thm:limiting density} and \ref{Thm_qJUE scaling limits}) as given and use the duality relations to deduce the corresponding LPP results (Theorems~\ref{Thm_LPP geometric weak}, \ref{Thm_LPP geometric strong}, and \ref{Thm_exponential}).
Thus, the proofs in this subsection are essentially consequences of the duality, while the underlying asymptotic analysis of the $q$-JUE will be established in the subsequent sections. 
 
Before turning to the proofs, we note two minor technicalities.
Since the geometric LPP random variables are integer-valued, the exact threshold arising from the duality is $q^{N+\ell-1}$ rather than the canonical edge scaling in Theorem~\ref{Thm_qJUE scaling limits}.
This discrepancy is handled by a squeezing argument using monotonicity and continuity of the limiting Fredholm determinants.
Moreover, for $\gamma>1$, the right Jacobi parameter satisfies
\[
    \mathsf{b}_N
    =(\gamma-1)N+\big(\lfloor\gamma N\rfloor-\gamma N\big),
    \qquad
    \lfloor\gamma N\rfloor-\gamma N\in[-1,0].
\]
The uniformity in $d_2$ in
Theorem~\ref{Thm_qJUE scaling limits} therefore allows us to
apply that theorem with $c_2=\gamma-1$.
 
We first apply Proposition~\ref{Prop_duality Gnm q-JUE} to prove Theorems~\ref{Thm_LPP geometric weak} and \ref{Thm_LPP geometric strong}, assuming the $q$-JUE asymptotic results established in Theorems~\ref{thm:limiting density} and \ref{Thm_qJUE scaling limits}.

\begin{proof}[Proof of Theorems~\ref{Thm_LPP geometric weak} and ~\ref{Thm_LPP geometric strong}] 

By Proposition~\ref{Prop_duality Gnm q-JUE}, for each integer $\ell\ge0$,
\begin{equation} \label{duality weak proof}
    \mathbb{P}(G_{\lfloor \gamma N\rfloor,N}\le \ell)
    =
    \mathbb{P}\Big(
    x_{\min}^{(N)}
    \ge q^{N+\ell-1}
    \Big),
\end{equation}
where $x_{\min}^{(N)}$ is the smallest particle of the $q$-JUE of size $N$ with parameters
\[
    \mathsf a_N  = \frac{2\log\alpha}{\log q}-1 =-\frac{2\log\alpha}{\lambda}N-1,
    \qquad
    \mathsf b_N
    =
    \lfloor \gamma N\rfloor-N.
\]

We first prove Theorem~\ref{Thm_LPP geometric weak}. 
For notational simplicity, throughout the proof we write 
\[
    \omega=\omega(\gamma,\alpha,\lambda),
    \qquad
    \sigma=\sigma(\gamma,\alpha,\lambda),
    \qquad
    \theta=\theta_{\gamma,\alpha,\lambda}.
\]
Under the assumptions of Theorem~\ref{Thm_LPP geometric weak}, the
corresponding $q$-JUE lies in the soft-edge regime of Theorem~\ref{Thm_qJUE scaling limits}, with macroscopic parameters 
\begin{equation} \label{weak LPP parameter matching}
    c_1=-\frac{2\log\alpha}{\lambda}>0,
    \qquad 
    c_2=\gamma-1 .
\end{equation}
Then the left endpoint in \eqref{eqn:def of endpoints ab}
is precisely $\mathfrak a$ in \eqref{def of a weak LPP}, the right endpoint is $\mathfrak b$ in \eqref{def of b weak LPP}, and the
soft edge constant in \eqref{def of theta} is $\theta$ in \eqref{def of theta weak LPP}. Moreover,
$\mathsf{s}(1+\mathsf{s}^{c_1}-\mathsf{s}^{c_1+c_2+1})
=\Xi_{\gamma,\alpha,\lambda}$, so the two case distinctions in the
definition of the soft-edge constant agree. Hence
Theorem~\ref{Thm_qJUE scaling limits} \textup{(ii)} gives, for each $r\in\mathbb R$,
\begin{equation} \label{qJUE soft edge weak proof}
    \lim_{N\to\infty}
    \mathbb{P}\Big(
    x_{\min}^{(N)}
    \ge
    \mathfrak a-\frac{r}{(\theta N)^{2/3}}
    \Big)
    =
    F_{\rm TW}(r).
\end{equation}
Moreover, by the definition of $\mathfrak a$ and $\omega$, we have
\begin{equation} \label{a omega weak proof}
    \mathfrak a
    =
    e^{-\lambda(1+\omega)} .
\end{equation}

We first prove the law of large numbers. Fix $\varepsilon>0$ and set $\ell_N^+
    :=\lfloor N(\omega+\varepsilon)\rfloor .$ 
Then, using \eqref{a omega weak proof},
\[
    q^{N+\ell_N^+-1}
    =
    \exp\Big(
    -\frac{\lambda}{N}(N+\ell_N^+-1)
    \Big)
    =
    \mathfrak a e^{-\lambda\varepsilon+o(1)} .
\]
Since $e^{-\lambda\varepsilon}<1$, this threshold lies a fixed distance to
the left of $\mathfrak a$ for all sufficiently large $N$. Therefore, for
each fixed $R>0$ and all sufficiently large $N$,
\[
    q^{N+\ell_N^+-1}
    \le
    \mathfrak a-\frac{R}{(\theta N)^{2/3}} .
\]
Then by \eqref{duality weak proof} and \eqref{qJUE soft edge weak proof}, we obtain
\[
\begin{split}
    \liminf_{N\to\infty}
    \mathbb{P}
    \Big(
    G_{\lfloor \gamma N\rfloor,N}
    \le
    \ell_N^+
    \Big)
    &\ge
    \lim_{N\to\infty}
    \mathbb{P}\Big(
    x_{\min}^{(N)}
    \ge
    \mathfrak a-\frac{R}{(\theta N)^{2/3}}
    \Big)  = F_{\rm TW}(R).
\end{split}
\]
Letting $R\to\infty$ gives
\[
    \lim_{N\to\infty}
    \mathbb{P}
    \Big(
    G_{\lfloor \gamma N\rfloor,N}
    \le
    N(\omega+\varepsilon)
    \Big)
    =
    1 .
\] 
Similarly, it follows that 
\[
    \lim_{N\to\infty}
    \mathbb{P}
    \Big(
    G_{\lfloor \gamma N\rfloor,N}
    \le
    N(\omega-\varepsilon)
    \Big)
    =
    0 .
\]
Combining the two estimates gives \eqref{LLN for inhomo LPP}.

It remains to prove the fluctuation limit. Fix $s\in\mathbb R$ and define
\[
    \ell_N(s)
    :=
    \lfloor
    N\omega+\sigma N^{1/3}s
    \rfloor  =
    N\omega+\sigma N^{1/3}s+O(1).
\] 
It follows from \eqref{a omega weak proof} and the definition \eqref{def of sigma weak LPP} that
\begin{equation} \label{threshold matching weak proof}
    q^{N+\ell_N(s)-1}
    =
    \mathfrak a
    -
    \lambda\mathfrak a\sigma sN^{-2/3}
    +
    o(N^{-2/3})  =
    \mathfrak a
    -
    \frac{s}{(\theta N)^{2/3}}
    +
    o(N^{-2/3}).
\end{equation}
For any $\delta>0$, \eqref{threshold matching weak proof} implies that, for all
sufficiently large $N$,
\[
    \mathfrak a-\frac{s+\delta}{(\theta N)^{2/3}}
    \le
    q^{N+\ell_N(s)-1}
    \le
    \mathfrak a-\frac{s-\delta}{(\theta N)^{2/3}} .
\]
Using monotonicity and \eqref{qJUE soft edge weak proof}, we obtain 
\[
    F_{\rm TW}(s-\delta)
    \le
    \liminf_{N\to\infty}
    \mathbb{P}\Big(
    x_{\min}^{(N)}
    \ge
    q^{N+\ell_N(s)-1}
    \Big), \qquad  \limsup_{N\to\infty}
    \mathbb{P}\Big(
    x_{\min}^{(N)}
    \ge
    q^{N+\ell_N(s)-1}
    \Big)
    \le
    F_{\rm TW}(s+\delta).
\] 
Since $F_{\rm TW}$ is continuous, letting $\delta\downarrow0$ yields
\[
    \lim_{N\to\infty}
    \mathbb{P}\Big(
    x_{\min}^{(N)}
    \ge
    q^{N+\ell_N(s)-1}
    \Big)
    =
    F_{\rm TW}(s).
\]
By the duality \eqref{duality weak proof}, we obtain \eqref{Fluctuation for inhomo LPP weak}, which completes the proof of Theorem~\ref{Thm_LPP geometric weak}. 

Next, we prove Theorem~\ref{Thm_LPP geometric strong}.
Since $\log q=-\lambda/N$ and $\log\alpha_N=-\eta/N$, the corresponding
$q$-JUE lies in the hard-edge regime of
Theorem~\ref{Thm_qJUE scaling limits}, with parameters
\[
c_1=0,
\qquad
d_1=\mathsf{a}=\frac{2\eta}{\lambda}-1,
\qquad
c_2=\gamma-1.
\]
Under this matching, the hard-edge constant \eqref{def of kappa} becomes
precisely $\kappa_{\gamma,\lambda}$ in
\eqref{def of kappa strong LPP}. Therefore, for every $u>0$,
\begin{equation} \label{qJUE hard edge strong proof}
    \lim_{N\to\infty}
    \mathbb{P}\Big(
    x_{\min}^{(N)}
    >
    \frac{u}{\kappa_{\gamma,\lambda}^{2}N^{2}}
    \Big)
    =
    \det\bigl(I-K_{\mathsf{a}}^{\rm Bessel}\bigr)_{L^2(0,u)} .
\end{equation}  
The remainder of the argument is analogous to the soft-edge case above, with \eqref{qJUE hard edge strong proof} replacing \eqref{qJUE soft edge weak proof}, which completes the proof. 
\end{proof}

Finally, we prove Theorem~\ref{Thm_exponential}.
In contrast to Theorems~\ref{Thm_LPP geometric weak} and
\ref{Thm_LPP geometric strong}, no $q$-deformed asymptotic analysis is
needed here.
Indeed, by Corollary~\ref{Cor_duality Hnm JUE}, the dual model is the
classical JUE, whose hard-edge scaling limit is standard.
We shall not repeat the elementary monotonicity argument used in the
geometric cases.

\begin{proof}[Proof of Theorem~\ref{Thm_exponential}]
Write $\nu:=2\beta+1=2\eta/\lambda-1>-1$.
By Corollary~\ref{Cor_duality Hnm JUE}, for each $\ell\ge0$ we have
\begin{equation} \label{duality exponential proof}
    \mathbb{P}(H_{\lfloor \gamma N\rfloor,N}\le \ell)
    =
    \mathbb{P}\Big(
    y_{\min}^{(N)}(\nu,\lfloor \gamma N\rfloor-N)
    \ge e^{-\lambda \ell}
    \Big),
\end{equation}
where $y_{\min}^{(N)}(\nu,\lfloor \gamma N\rfloor-N)$ denotes the smallest eigenvalue of the JUE
of size $N$ with parameters $\nu$ and $\lfloor \gamma N\rfloor-N$.

We recall the hard-edge scaling for this JUE.
The limiting Wachter density with left parameter of order one and right
parameter $(\gamma-1)N+O(1)$ has the hard-edge behaviour
\[
    \rho_{\rm J}(x)
    =
    \frac{\sqrt{\gamma}}{\pi\sqrt{x}}\bigl(1+o(1)\bigr),
    \qquad x\downarrow0.
\]
Then the standard hard-edge scaling theorem for the JUE (see e.g. \cite[Theorem 9 and Section 4.2]{HM12})
gives, for every $u>0$,
\begin{equation} \label{JUE hard edge exponential proof}
    \lim_{N\to\infty}
    \mathbb{P}\Big(
    4\gamma N^2
    y_{\min}^{(N)}(\nu,\lfloor \gamma N\rfloor-N)
    \ge u
    \Big)
    =
    \det\bigl(I-K_{\nu}^{\rm Bessel}\bigr)_{L^2(0,u)} .
\end{equation}  

By \eqref{duality exponential proof}, we have
\begin{align}
\begin{split}  \label{exponential threshold proof}
\mathbb{P}\Big[
    H_{\lfloor \gamma N\rfloor,N}
    -
    \frac{2}{\lambda}\log N
    \le
    \frac{\xi+\log\gamma}{\lambda}
    \Big]
    & =
    \mathbb{P}\Big(
    \gamma N^2
    y_{\min}^{(N)}(\nu,\lfloor \gamma N\rfloor-N)
    \ge e^{-\xi}
    \Big)  \\
    & =
    \mathbb{P}\Big(
    4\gamma N^2
    y_{\min}^{(N)}(\nu,\lfloor \gamma N\rfloor-N)
    \ge 4e^{-\xi}
    \Big).
\end{split}
\end{align}
Applying \eqref{JUE hard edge exponential proof} with
$u=4e^{-\xi}$ yields
\begin{equation} \label{Bessel determinant before change variable}
    \lim_{N\to\infty}
    \mathbb{P}\Big[
    H_{\lfloor \gamma N\rfloor,N}
    -
    \frac{2}{\lambda}\log N
    \le
    \frac{\xi+\log\gamma}{\lambda}
    \Big]
    =
    \det\bigl(I-K_{\nu}^{\rm Bessel}\bigr)_{L^2(0,4e^{-\xi})}.
\end{equation}

It remains only to rewrite the Fredholm determinant in the form stated in
the theorem.
By the Fredholm expansion,
\[
\det\bigl(I-K_{\nu}^{\rm Bessel}\bigr)_{L^2(0,4e^{-\xi})}
=
1+
\sum_{k=1}^{\infty}
\frac{(-1)^k}{k!}
\int_{(0,4e^{-\xi})^k}
\det\bigl(K_{\nu}^{\rm Bessel}(x_i,x_j)\bigr)_{i,j=1}^k
\prod_{j=1}^{k}dx_j .
\]
In each integral, make the change of variables $ x_j=4e^{-u_j}$ $(j=1,\dots,k)$.  
Then we have 
\begin{align*}
    \int_{(0,4e^{-\xi})^k}
    \det\bigl(K_{\nu}^{\rm Bessel}(x_i,x_j)\bigr)_{i,j=1}^k
    \prod_{j=1}^{k}dx_j  & =
    \int_{(\xi,\infty)^k}
    \det\bigl(K_{\nu}^{\rm Bessel}(4e^{-u_i},4e^{-u_j})\bigr)_{i,j=1}^k
    \prod_{j=1}^{k}4e^{-u_j}\,du_j  \\
    & =
    \int_{(\xi,\infty)^k}
    \det\Bigl(
    4e^{-(u_i+u_j)/2}
    K_{\nu}^{\rm Bessel}(4e^{-u_i},4e^{-u_j})
    \Bigr)_{i,j=1}^k
    \prod_{j=1}^{k}du_j .
\end{align*}
Since $\nu=2\beta+1$, the desired fluctuation limit \eqref{Fluctuation limit for exponential LPP} follows from \eqref{Bessel determinant before change variable}.
Finally, as in the geometric case above, the law of large numbers
\eqref{LLN for inhomo LPP exp} follows immediately from the fluctuation
limit, since the centred fluctuations are of order $O(1)$, whereas the
leading-order term is $(2/\lambda)\log N$. This completes the proof. 
\end{proof}

\section{Equilibrium measure and proof of Theorem~\ref{thm:limiting density}}   \label{Sec_eq msr}

This section is devoted to the proof of Theorem~\ref{thm:limiting density}. Our approach follows the strategy developed in \cite{BJM26} for the $q$-deformed Laguerre unitary ensemble and is based on a constrained logarithmic equilibrium problem.
This is also in the spirit of Johansson's derivation \cite{Joh00} of the limiting density \eqref{eqn:Meixner limiting density} for the Meixner ensemble, although the form of the constraint is different due to the underlying lattice structure.

More precisely, the large-$N$ behaviour of the $q$-JUE is governed by the unique minimiser of an associated energy functional subject to an upper constraint arising from the underlying $q$-lattice. Once the candidate equilibrium measure is given explicitly, it therefore suffices to verify the corresponding Euler--Lagrange variational conditions. This is the approach we take in the present section: we show directly that the measure stated in Theorem~\ref{thm:limiting density} satisfies these conditions, which, by uniqueness of the equilibrium measure, completes the proof.

While this verification provides a complete proof of the theorem, it does not explain how the explicit formula for the limiting density is found in the first place. For completeness, and to make the origin of this formula transparent, we present its derivation in Appendix~\ref{Appendix_spectral moment} based on the moment method.

\medskip 
  
We now introduce the equilibrium measure problem. Recall that the weight function $w^{(\mathsf{a},\mathsf{b})}$ is defined in \eqref{eqn:def of qjacobi weight}. Under the scaling \eqref{eqn:scale of alphabet}, the asymptotic formula for the $q$-Pochhammer symbol \cite[Theorem~4]{Mc99} yields
\begin{equation}\label{eqn:log weight asymptotic}
    \log w^{(\mathsf{a},\mathsf{b})}(x;q)
    =-NV(x)+O(1),
    \qquad N\to\infty,
\end{equation}
locally uniformly for $x\in(0,1)$, where 
\begin{equation}\label{eqn:def of potential}
    V(x)
  :  =\frac{1}{\lambda}
    \big(\Li(x)-\Li(\mathsf{s}^{c_{2}}x)\big)
    -c_{1}\log x, \qquad x \in (0,1). 
\end{equation}
Here, $\Li$ denotes the dilogarithm function defined by
\begin{equation} \label{def of Li}
    \Li(x):=-\int_0^x\frac{\log(1-t)}{t}\,dt. 
\end{equation}
The associated energy functional is 
\begin{equation}\label{eqn:def of I}
    I[\mu]
    =-\int_{0}^{1}\int_{0}^{1}\log|x-y|\,\mu(dx)\,\mu(dy)
    +\int_{0}^{1}V(x)\,\mu(dx) 
\end{equation}
over the class of probability measures satisfying the upper constraint
\begin{equation}
    \mathcal{P}_{\lambda}(0,1)
    :=\Big\{\mu(dx)\Big\vert\supp \mu\subset[0,1],~\frac{d\mu}{dx}\leq\frac{1}{\lambda x}\Big\}.
\end{equation} 
Here, the upper constraint reflects the spacing of the underlying $q$-lattice \cite[Remark 7]{BJM26}; see also \cite{VV21} for its appearance in the study of limiting zero distributions of orthogonal polynomials.

Following the large-deviation argument of \cite[Theorem~1.2]{HMO25} with $\eta=\theta=1$, as in the proof of \cite[Theorem~1.2(B)]{BJM26}, the empirical measures concentrate exponentially at speed $N^2$ around the unique minimiser of $I$ over $\mathcal{P}_{\lambda}(0,1)$. 
The same argument accommodates the fixed shifts $d_1,d_2$, which do not affect the limiting rate function. The Borel--Cantelli lemma then yields almost-sure convergence to this minimiser. Thus, it remains to identify the minimiser with $\rho^{(c_1,c_2)}(x)\,dx$.

A standard way to identify the equilibrium measure is to verify the
associated Euler--Lagrange variational conditions (see e.g. \cite[Section~2.1.4]{BKMM07}). Accordingly, the proof of Theorem~\ref{thm:limiting density} reduces to the following proposition. 

\begin{prop} \label{Prop_EL conditions}
    Let $\mu(dx)=\rho^{(c_{1},c_{2})}(x)\,dx$, where $\rho^{(c_1,c_2)}$ is given in \eqref{def of limiting law}, and write
    \begin{equation}\label{eqn:def of d}
        d:=
        \begin{cases}
            b&\textup{if }\lambda \leq\lambda_{c}, \smallskip \\
            1&\textup{if }\lambda >\lambda_{c}.
        \end{cases}
    \end{equation}
  Then there exists a constant $\Omega\in\mathbb{R}$ such that 
    \begin{equation}\label{eqn:Euler-Lagrange}
        -2\int_{a}^{d}\log|x-y|\,\mu(dy)+V(x)
        \begin{cases}
            =-\Omega & \textup{if } x\in(a,b),\smallskip\\
            \geq-\Omega & \textup{if } x\in(0,a)\cup(d,1),\smallskip\\
            \leq-\Omega & \textup{if } x\in(b,d).
        \end{cases}
    \end{equation}
\end{prop}

Before the proof of Proposition~\ref{Prop_EL conditions}, we first evaluate the Cauchy transform 
\begin{equation}
  g(z) :=\int_a^d\frac{\mu(dx)}{z-x}.
\end{equation} 
For this, recalling that $a$ and $b$ are given by \eqref{eqn:def of endpoints ab}, we introduce 
\begin{equation}\label{eqn:EL discriminant}
    R(z)
    :=\sqrt{(Bz-1-\mathsf{s}^{c_{1}})^{2}-4\mathsf{s}^{c_{1}}(1-z)(1-\mathsf{s}^{c_{2}}z)}
    =(\mathsf{s}^{-1}-\mathsf{s}^{c_1+c_2+1})\sqrt{(z-a)(z-b)},
\end{equation}
where 
\begin{equation} \label{def of B in terms of s}
B:=\mathsf{s}^{-1}+\mathsf{s}^{c_1+c_2+1}. 
\end{equation}
Here, the branch is analytic in $\mathbb{C}\setminus[a,b]$ and
normalised by $R(z) \sim (\mathsf{s}^{-1}-\mathsf{s}^{c_1+c_2+1}) z$ as $z \to \infty.$ 

\begin{lem} \label{lem:Cauchy transform}
We have 
\begin{equation}\label{eqn:EL Cauchy transform}
\begin{split}
g(z) =\frac{1}{\lambda z} \log\frac{Bz-1-\mathsf{s}^{c_1}+R(z)}{2(z-1)}, \qquad z\in\mathbb{C}\setminus[a,d].
\end{split}
\end{equation}
The logarithm is chosen by analytic continuation from $z>1$, where its argument is positive.
Moreover, 
\begin{equation}\label{eqn:EL g band values}
    g_\pm(x)
    =
    \frac{1}{2\lambda x}
    \log\frac{\mathsf{s}^{c_1}(1-\mathsf{s}^{c_2}x)}{1-x}
    \mp\pi i\,\rho^{(c_1,c_2)}(x),
    \qquad x\in(a,b),
\end{equation}
where $g_\pm$ denote boundary values from the upper and lower half-planes, respectively.
\end{lem}

\begin{proof}
By \eqref{eqn:EL discriminant}, our choice of branch implies that $R(x)<0$ for $x<a$, $R(x)>0$ for $x>b$, and
\begin{equation}\label{eqn:EL R boundary}
    R_+(x)=-R_-(x)
    =
    i\sqrt{
        4\mathsf{s}^{c_1}(1-x)(1-\mathsf{s}^{c_2}x)
        -(Bx-1-\mathsf{s}^{c_1})^2
    },
    \qquad x\in(a,b).
\end{equation}

To establish \eqref{eqn:EL Cauchy transform}, it suffices to show that
\begin{equation} \label{eqn:EL H definition}
\int_a^d\frac{\mu(dx)}{z-x} = \frac{\log H(z)}{ \lambda z }, \qquad H(z):=\frac{Bz-1-\mathsf{s}^{c_1}+R(z)}{2(z-1)}.
\end{equation} 
We verify this identity by computing the boundary jump of $\log H(z)/(\lambda z)$ across $[a,d]$.

We first examine the analyticity of $\log H(z)/(\lambda z)$.
Note that by \eqref{eqn:EL discriminant},
\begin{equation}\label{eqn:EL product identity}
\begin{split}
\big(Bz-1-\mathsf{s}^{c_1}+R(z)\big) \big(Bz-1-\mathsf{s}^{c_1}-R(z)\big) =4\mathsf{s}^{c_1}(1-z)(1-\mathsf{s}^{c_2}z).
\end{split}
\end{equation}
Consequently, a zero of the numerator of $H$ can occur only at
$z=1$ or $z=\mathsf{s}^{-c_2}$.
When $c_2>0$, we have
\begin{equation*}
\begin{split}
    B\mathsf{s}^{-c_2}-1-\mathsf{s}^{c_1}
    &=
    (1-\mathsf{s})(\mathsf{s}^{-1}-\mathsf{s}^{c_1})
    +\mathsf{s}^{-1}(\mathsf{s}^{-c_2}-1)>0.
\end{split}
\end{equation*}
Since $c_2>0$, we have $\mathsf{s}^{-c_2}>1\geq b$, which implies that $R(\mathsf{s}^{-c_2})>0$. Thus the numerator does not vanish there. 
To examine $z=1$, notice first that by \eqref{eqn:lambdac condition}, we have 
\begin{equation} \label{sign of B}
    B
    \begin{cases}
        <1+\mathsf{s}^{c_1}, &\textup{if } \lambda<\lambda_c,
        \smallskip 
        \\
        =1+\mathsf{s}^{c_1},  &\textup{if } \lambda=\lambda_c,
        \smallskip 
        \\
        >1+\mathsf{s}^{c_1},  &\textup{if } \lambda>\lambda_c.
    \end{cases}
\end{equation}
If $\lambda<\lambda_c$, then
$R(1)=1+\mathsf{s}^{c_1}-B$.
By \eqref{eqn:EL H definition} and \eqref{eqn:EL product identity}, we have 
\begin{equation*}
    H(z)
    =
    \frac{2\mathsf{s}^{c_1}(1-\mathsf{s}^{c_2}z)}
    {1+\mathsf{s}^{c_1}-Bz+R(z)}.
\end{equation*}
Therefore, the singularity at $z=1$ is removable, with $H(1)>0.$ 
In the remaining cases, $1\in[a,d]$.
It follows that $H$ is analytic and non-vanishing in
$\mathbb{C}\setminus[a,d]$.

On the other hand, the normalisation $R(z)\sim(\mathsf{s}^{-1}-\mathsf{s}^{c_1+c_2+1})z$ as $z\to\infty$ gives
\begin{equation}\label{eqn:EL auxiliary infinity}
   H(z)=\mathsf{s}^{-1}+O(z^{-1}), \qquad  \log H(z) 
    =\lambda +O(z^{-1}),
    \qquad z\to\infty.
\end{equation}
If $a>0$, then $c_1>0$ and
\begin{equation*}
    R(0)=-(1-\mathsf{s}^{c_1}),
    \qquad H(0)=1.
\end{equation*}
Moreover, $H(x)>0$ for $x<0$, so the chosen logarithm satisfies
$\log H(0)=0$.
Thus, $\log H(z)/z$ has a removable singularity at $z=0$.

We next determine the boundary values of $\log H(z)$. 
When $\lambda<\lambda_c$, the function $H$ is positive on $(b,\infty)$, since it is real and non-vanishing there and tends to $\mathsf{s}^{-1}>0$ at infinity.
When $\lambda\geq\lambda_c$, as $z \to 1$, we have
\begin{equation}\label{eqn:EL H at one}
    H(z)
    =
    \begin{cases}
        \displaystyle
        \frac{\sqrt{\mathsf{s}^{c_1}(1-\mathsf{s}^{c_2})}}
        {\sqrt{z-1}}+O(1),
        & \textup{if } \lambda=\lambda_c,
        \smallskip
        \\
        \displaystyle
        \frac{B-1-\mathsf{s}^{c_1}}{z-1}+O(1),
        & \textup{if } \lambda>\lambda_c.
    \end{cases}
\end{equation}
When $\lambda=\lambda_c$, we have $b=1$, and continuation through the
upper half-plane gives $\arg H_+(x)\to-\pi/2$ as $x\uparrow1$.
On the other hand, when $\lambda>\lambda_c$, the pole in \eqref{eqn:EL H at one} has positive residue. Since $H$ is real and non-vanishing on $(b,1)$, it is negative throughout this interval, and continuation from $z>1$ gives 
\begin{equation}\label{eqn:EL saturated logarithm}
    \log H_\pm(x)=\log|H(x)|\mp\pi i,
    \qquad x\in(b,d).
\end{equation}

For $x\in(a,b)$, by using \eqref{eqn:EL R boundary}, we have 
\begin{equation*}
    H_+(x)
    =
    \frac{1+\mathsf{s}^{c_1}-Bx-R_+(x)}{2(1-x)},
    \qquad
    H_-(x)=\overline{H_+(x)}.
\end{equation*}
In particular, $\im H_+(x)<0$, and
\begin{equation*}
\begin{split}
    |H_+(x)|^2
    &=
    \frac{(Bx-1-\mathsf{s}^{c_1})^2+|R_+(x)|^2}
    {4(1-x)^2} =
    \frac{\mathsf{s}^{c_1}(1-\mathsf{s}^{c_2}x)}{1-x}.
\end{split}
\end{equation*}
The branch determination above shows that the argument of $H_+(x)$
lies in $(-\pi,0)$ throughout $(a,b)$.
Therefore,
\begin{equation*}
\begin{split}
    \arg H_+(x)
    &=
    -\arccos
    \frac{1+\mathsf{s}^{c_1}-Bx}
    {2\sqrt{\mathsf{s}^{c_1}(1-x)(1-\mathsf{s}^{c_2}x)}} =-\pi\lambda x\,\rho^{(c_1,c_2)}(x),
\end{split}
\end{equation*}
where the last equality follows from
\eqref{def of limiting law}.
Combining the above, we obtain
\begin{equation}\label{eqn:EL band logarithm}
    \log H_\pm(x)
    =
    \frac{1}{2}
    \log\frac{\mathsf{s}^{c_1}(1-\mathsf{s}^{c_2}x)}{1-x}
    \mp\pi i\lambda x\,\rho^{(c_1,c_2)}(x),
    \qquad x\in(a,b).
\end{equation}
Together with \eqref{eqn:EL saturated logarithm} and the identity
$\rho^{(c_1,c_2)}(x)=1/(\lambda x)$ on $(b,d)$, this gives
\begin{equation}\label{eqn:EL logarithm jump}
    \frac{\log H_+(x)-\log H_-(x)}{\lambda x}
    =-2\pi i\,\rho^{(c_1,c_2)}(x),
    \qquad x\in(a,d)\setminus\{b\}.
\end{equation}
Thus, applying the Sokhotski--Plemelj formula together with \eqref{eqn:EL auxiliary infinity} and \eqref{eqn:EL logarithm jump}, we obtain 
\begin{equation*}
\begin{split}
    \frac{\log H(z)}{\lambda z}
    &=
    \frac{1}{2\pi i}
    \int_a^d
    \frac{\log H_+(x)-\log H_-(x)}
    {\lambda x(x-z)}\,dx =
    \int_a^d\frac{\rho^{(c_1,c_2)}(x)}{z-x}\,dx.
\end{split}
\end{equation*}
This proves \eqref{eqn:EL Cauchy transform}, while
\eqref{eqn:EL g band values} follows from
\eqref{eqn:EL band logarithm}. 
\end{proof}

Using Lemma~\ref{lem:Cauchy transform}, we now prove Proposition~\ref{Prop_EL conditions}. 

\begin{proof}[Proof of Proposition~\ref{Prop_EL conditions}]
We use the notation of Lemma~\ref{lem:Cauchy transform} and set
\begin{equation}
    \mathcal{E}(x)
    :=-2\int_a^d\log|x-y|\,\mu(dy)+V(x).
\end{equation}
By using $\log\mathsf{s}=-\lambda$ and \eqref{def of Li}, we have 
\begin{equation}\label{eqn:EL V derivative}
\begin{split}
    V'(x) =
    \frac{1}{\lambda x}
    \log\frac{\mathsf{s}^{c_1}(1-\mathsf{s}^{c_2}x)}{1-x}.
\end{split}
\end{equation}
Thus, the boundary values \eqref{eqn:EL g band values} satisfy
\begin{equation*}
    g_+(x)+g_-(x)=V'(x),
    \qquad x\in(a,b).
\end{equation*}
By the Sokhotski--Plemelj formulas, it follows that for $x \in (a,b)$, 
\begin{equation*}
\begin{split}
    \mathcal{E}'(x)
    =
    V'(x)-2\,\operatorname{p.v.}
    \int_a^d\frac{\mu(dy)}{x-y} =V'(x)-g_+(x)-g_-(x)=0. 
\end{split}
\end{equation*}
Consequently, there exists a constant $\Omega\in\mathbb{R}$ such that $ \mathcal{E}(x)=-\Omega$ for $x \in (a,b)$. This gives the equality part of the Euler--Lagrange conditions \eqref{eqn:Euler-Lagrange}. 

We next verify the inequalities in \eqref{eqn:Euler-Lagrange}.
For $x\in(0,a)\cup(b,1)$, the function $R(x)$ is real.
Combining \eqref{eqn:EL discriminant}, \eqref{eqn:EL Cauchy transform} and \eqref{eqn:EL V derivative}, we obtain 
\begin{equation}\label{eqn:EL potential derivative}
\begin{split}
    \mathcal{E}'(x)
    &=V'(x)-2\re g_+(x) =
    \frac{1}{\lambda x}
    \log
    \frac{4\mathsf{s}^{c_1}(1-x)(1-\mathsf{s}^{c_2}x)}
    {(Bx-1-\mathsf{s}^{c_1}+R(x))^2}=
    \frac{1}{\lambda x}
    \log
    \frac{Bx-1-\mathsf{s}^{c_1}-R(x)}
    {Bx-1-\mathsf{s}^{c_1}+R(x)}.
\end{split}
\end{equation} 
Moreover, by using $\mathsf{s}^{c_1}(1-x)(1-\mathsf{s}^{c_2}x)>0$ for $x\in(0,1)$ and \eqref{eqn:EL discriminant}, we have 
\begin{equation}\label{eqn:EL R comparison}
    0<|R(x)|<|Bx-1-\mathsf{s}^{c_1}|,
    \qquad x\in(0,a)\cup(b,1).
\end{equation}
In particular, the argument of the logarithm in
\eqref{eqn:EL potential derivative} is positive.

Now it remains to determine the sign of $\mathcal{E}'$.
When $\lambda\leq\lambda_c$, by \eqref{sign of B}, we have $Bx-1-\mathsf{s}^{c_1}<0$ for every $x\in(0,1)$. On the other hand, when $\lambda>\lambda_c$, the point $(1+\mathsf{s}^{c_1})/B$ belongs to $(0,1)$, and the polynomial
under the square root in \eqref{eqn:EL discriminant} is negative
at this point.
Then by \eqref{eqn:EL discriminant}, we have 
\begin{equation*}
    a<\frac{1+\mathsf{s}^{c_1}}{B}<b.
\end{equation*}
Together with the definition \eqref{eqn:def of d} of $d$, these observations give
\begin{equation*}
    Bx-1-\mathsf{s}^{c_1}
    \begin{cases}
        <0, &\textup{if } x\in(0,a)\cup(d,1),
        \smallskip 
        \\
        >0, &\textup{if } x\in(b,d).
    \end{cases}
\end{equation*} 
On $(0,a)$, both $Bx-1-\mathsf{s}^{c_1}$ and $R(x)$ are negative,
so \eqref{eqn:EL R comparison} shows that the ratio in
\eqref{eqn:EL potential derivative} lies in $(0,1)$.
On $(d,1)$, the first quantity is negative and $R(x)>0$,
so the ratio is greater than $1$.
Finally, on $(b,d)$, both quantities are positive,
and the ratio again lies in $(0,1)$.
Therefore, we obtain 
\begin{equation*}
    \mathcal{E}'(x)
    \begin{cases}
        <0, &\textup{if } x\in(0,a)\cup(b,d),
        \smallskip 
        \\
        >0, &\textup{if } x\in(d,1).
    \end{cases}
\end{equation*} 
By continuity, $\mathcal{E}(a)=\mathcal{E}(b)=-\Omega$, and the desired
inequalities in \eqref{eqn:Euler-Lagrange} now follow. 
\end{proof}

\section{Fluctuations and proof of Theorem~\ref{Thm_qJUE scaling limits}}

In this section, we prove Theorem \ref{Thm_qJUE scaling limits}. 
When $c_2>0$, all estimates below, including those in Appendix~\ref{app:soft edge diagonal}, are uniform for $d_2$ in a fixed compact subset of $\mathbb{R}$. Indeed, the $q$-Pochhammer asymptotics are uniform for bounded shifts, while the limiting phase and the contours used in the soft-edge analysis are independent of $d_2$. The constants in the bounds and error estimates can therefore be chosen uniformly on such compact sets.

We first recall the determinantal structure of the $q$-JUE. The little $q$-Jacobi polynomials form a family of $q$-orthogonal polynomials $\{p^{(\mathsf{a},\mathsf{b})}_{n}\}_{n\geq0}$ with respect to the little $q$-Jacobi weight \eqref{eqn:def of qjacobi weight} with orthogonality
\begin{equation}
    \int_{0}^{1}p_{m}^{(\mathsf{a},\mathsf{b})}(x;q)p_{n}^{(\mathsf{a},\mathsf{b})}(x;q)w^{(\mathsf{a},\mathsf{b})}(x;q)\,d_{q}x
    =h_{n}^{(\mathsf{a},\mathsf{b})}\delta_{mn}.
\end{equation}
Here, $d_{q}x$ indicates the Jackson $q$-measure \eqref{def of Jackson measure}, and
\begin{equation}\label{eqn:def of hn}
    h_{n}^{(\mathsf{a},\mathsf{b})}
    =\frac{q^{n(\mathsf{a}+1)}(1-q)}{1-q^{2n+\mathsf{a}+\mathsf{b}+1}}\frac{(q,q^{\mathsf{a}+\mathsf{b}+1};q)_{\infty}}{(q^{\mathsf{a}+1},q^{\mathsf{b}+1};q)_{\infty}}\frac{(q,q^{\mathsf{b}+1};q)_{n}}{(q^{\mathsf{a}+1},q^{\mathsf{a}+\mathsf{b}+1};q)_{n}}
\end{equation}
is the squared orthogonal norm. (Expressions with removable singularities in the parameters are understood by continuity.)
The little $q$-Jacobi polynomial admits the basic hypergeometric series representation
\begin{equation}\label{eqn:littleq jacobi hypergeometric}
    p_{n}^{(\mathsf{a},\mathsf{b})}(x;q)
    ={}_{2}\phi_{1}\Big(\begin{matrix}q^{-n},q^{n+\mathsf{a}+\mathsf{b}+1}\\q^{\mathsf{a}+1}\end{matrix};q,qx\Big)
    =\sum_{k=0}^{n}\frac{(q^{-n},q^{n+\mathsf{a}+\mathsf{b}+1};q)_{k}}{(q,q^{\mathsf{a}+1};q)_{k}}(qx)^{k};
\end{equation}
see \cite[Eq. (14.12.1)]{KLS10}.

The $q$-JUE forms a determinantal point process with respect to the Jackson $q$-measure whose correlation kernel is
\begin{equation} \label{def of kernel OP sum}
    K_{N}^{(\mathsf{a},\mathsf{b})}(x,y;q)
    :=\sqrt{w^{(\mathsf{a},\mathsf{b})}(x;q)w^{(\mathsf{a},\mathsf{b})}(y;q)}\sum_{j=0}^{N-1}\frac{1}{h_{j}^{(\mathsf{a},\mathsf{b})}}p_{j}^{(\mathsf{a},\mathsf{b})}(x;q)p_{j}^{(\mathsf{a},\mathsf{b})}(y;q).
\end{equation}
The Christoffel-Darboux formula yields the representation of the correlation kernel  
\begin{equation}\label{eqn:CD kernel}
    K_{N}^{(\mathsf{a},\mathsf{b})}(x,y;q)
    =\frac{\sqrt{w^{(\mathsf{a},\mathsf{b})}(x;q)w^{(\mathsf{a},\mathsf{b})}(y;q)}}{h_{N-1}^{(\mathsf{a},\mathsf{b})}}\frac{\gamma_{N-1}^{(\mathsf{a},\mathsf{b})}}{\gamma_{N}^{(\mathsf{a},\mathsf{b})}}\frac{p_{N}^{(\mathsf{a},\mathsf{b})}(x;q)p_{N-1}^{(\mathsf{a},\mathsf{b})}(y;q)-p_{N-1}^{(\mathsf{a},\mathsf{b})}(x;q)p_{N}^{(\mathsf{a},\mathsf{b})}(y;q)}{x-y}.
\end{equation}
Here,
\begin{equation}\label{eqn:def of gammaN}
    \gamma_{N}^{(\mathsf{a},\mathsf{b})}
    =\frac{(-1)^{N}q^{-\frac{N(N-1)}{2}}(q^{N+\mathsf{a}+\mathsf{b}+1};q)_{N}}{(q^{\mathsf{a}+1};q)_{N}}
\end{equation}
is the leading coefficient of the little $q$-Jacobi polynomial $p_{N}^{(\mathsf{a},\mathsf{b})}(x)$.

\subsection{Hard edge scaling limit}

Throughout this subsection, we assume $c_{1}=0$ and prove Theorem \ref{Thm_qJUE scaling limits} (i). 

We begin with the following lemma.

\begin{lem}
    Let $q=e^{-\lambda/N}$ for a fixed $\lambda>0$, and assume
    $\mathsf{a}=d_1>-1$. For each fixed $j\in\mathbb{Z}_{\geq0}$,
    as $N\to\infty$,
    \begin{equation}\label{eqn:hard edge polynomial}
        p_{N-j}^{(\mathsf{a},\mathsf{b})}
        \Big(\frac{u}{(\kappa N)^2};q\Big)
        =
        \Gamma(d_1+1)2^{d_1}u^{-d_1/2}J_{d_1}(\sqrt{u})+o(1),
    \end{equation}
    locally uniformly for $u\in\mathbb{C}$, together with derivatives
    of any fixed order with respect to $u$.  
\end{lem}

\begin{proof}
    The basic hypergeometric representation
    \eqref{eqn:littleq jacobi hypergeometric} gives
    \begin{equation*}
        p_{N-j}^{(\mathsf{a},\mathsf{b})}
        \Big(\frac{u}{(\kappa N)^2};q\Big)
        =
        \sum_{k=0}^{N-j}
        \frac{(q^{-N+j},q^{N+\mathsf{a}+\mathsf{b}-j+1};q)_k}
             {(q,q^{\mathsf{a}+1};q)_k}
        \frac{q^k u^k}{(\kappa N)^{2k}}.
    \end{equation*}
    For each fixed $k$, a direct calculation using the definition
    \eqref{def of kappa} of $\kappa$ yields
    \begin{equation*}
        \frac{(q^{-N+j},q^{N+\mathsf{a}+\mathsf{b}-j+1};q)_k}
             {(q,q^{\mathsf{a}+1};q)_k}
        \frac{q^k}{(\kappa N)^{2k}}
        =
        \Big(-\frac14\Big)^k
        \frac{\Gamma(d_1+1)}{k!\Gamma(d_1+k+1)}
        +O(N^{-1}).
    \end{equation*}

    To justify passage to the limit in the sum, we derive a summable
    bound independent of $N$. Set $c=\min\{1,d_1+1\}>0$.
    Since $\mathsf{a}+\ell\ge c\ell$ for every $\ell\ge1$,
    monotonicity and concavity give
    \begin{equation*}
        1-q^{\mathsf{a}+\ell}
        \ge 1-q^{c\ell}
        \ge c(1-q^\ell).
    \end{equation*}
    Moreover, by concavity of $1-e^{-\lambda t}$,
    \begin{equation*}
        1-q^\ell
        \ge \frac{\ell}{N}(1-\mathsf{s}),
        \qquad 1\le\ell\le N.
    \end{equation*}
    Consequently, for $0\le k\le N-j$,
    \begin{equation*}
        (q^{\mathsf{a}+1};q)_k
        \ge c^k(q;q)_k,
        \qquad
        (q;q)_k
        \ge \Big(\frac{1-\mathsf{s}}{N}\Big)^k k!.
    \end{equation*}
    For all sufficiently large $N$, we also have
    \begin{equation*}
        \big|(q^{-N+j},q^{N+\mathsf{a}+\mathsf{b}-j+1};q)_k\big|
        \le \big|(q^{-N+j};q)_k\big|
        \le (\mathsf{s}^{-1}-1)^k.
    \end{equation*}
    Thus, for each fixed $R>0$ and all $|u|\le R$,
    \begin{equation*}
        \Big|
        \frac{(q^{-N+j},q^{N+\mathsf{a}+\mathsf{b}-j+1};q)_k}
             {(q,q^{\mathsf{a}+1};q)_k}
        \frac{q^k u^k}{(\kappa N)^{2k}}
        \Big|
        \le
        \frac{1}{(k!)^2}
        \Big(\frac{R}{c\kappa^2\mathsf{s}(1-\mathsf{s})}\Big)^k.
    \end{equation*}
    The right-hand side is summable over $k\ge0$. Extending the
    summand by zero for $k>N-j$, dominated convergence therefore gives \eqref{eqn:hard edge polynomial} uniformly for $|u|\le R$. Since $R$ is arbitrary, this proves locally uniform convergence on $\mathbb{C}$. The convergence of derivatives of any fixed order follows from Cauchy's integral
    formula.
\end{proof}

We next use the polynomial asymptotics to derive the hard-edge limit of the correlation kernel. To obtain uniform convergence up to the origin for the full range $d_1>-1$, we first remove the factor $(uv)^{d_1/2}$ from the rescaled kernels.

\begin{prop}\label{prop:hard edge uniform}
    Let $q=e^{-\lambda/N}$ for a fixed $\lambda>0$, and let
    $\mathsf{a}=d_1>-1$ and $\mathsf{b}=c_2N+d_2$.
    For a fixed $M>0$, define
    \begin{equation}
        \mathfrak{L}_{M}^{\textup{hard}}
        =
        \Big\{u\in(0,M):
        \frac{u}{(\kappa N)^2}\in L_q\Big\},
    \end{equation}
    where $L_q=\{1,q,q^2,\ldots\}$. Then, as $N\to\infty$,
    \begin{equation}
        \sup_{u,v\in\mathfrak{L}_{M}^{\textup{hard}}}
        (uv)^{-d_1/2}
        \Big|
        \frac{1}{(\kappa N)^2}
        K_N^{(\mathsf{a},\mathsf{b})}
        \Big(\frac{u}{(\kappa N)^2},
             \frac{v}{(\kappa N)^2};q\Big)
        -
        K_{d_1}^{\textup{Bessel}}(u,v)
        \Big|
        \to 0.
    \end{equation} 
    In particular, if $d_1\ge0$, the same convergence holds without
    the factor $(uv)^{-d_1/2}$.
\end{prop}

\begin{proof}
    Comparing coefficients in the basic hypergeometric representation,
    we obtain
    \begin{equation}
        p_N^{(\mathsf{a},\mathsf{b})}(x;q)
        -
        p_{N-1}^{(\mathsf{a},\mathsf{b})}(x;q)
        =
        -\frac{q^{1-N}(1-q^{2N+\mathsf{a}+\mathsf{b}})}
              {1-q^{\mathsf{a}+1}}
        x p_{N-1}^{(\mathsf{a}+1,\mathsf{b})}(x;q).
    \end{equation}
    For convenience, write
    \begin{equation*}
        f_N(u)
        =
        p_N^{(\mathsf{a},\mathsf{b})}
        \Big(\frac{u}{(\kappa N)^2};q\Big),
        \qquad
        g_N(u)
        =
        u p_{N-1}^{(\mathsf{a}+1,\mathsf{b})}
        \Big(\frac{u}{(\kappa N)^2};q\Big).
    \end{equation*}
    Then we have 
    \begin{align*}
 & \quad  \frac{(\kappa N)^2}{u-v}
            \Big[
            p_N^{(\mathsf{a},\mathsf{b})}
            \Big(\frac{u}{(\kappa N)^2};q\Big)
            p_{N-1}^{(\mathsf{a},\mathsf{b})}
            \Big(\frac{v}{(\kappa N)^2};q\Big) -  p_{N-1}^{(\mathsf{a},\mathsf{b})}
            \Big(\frac{u}{(\kappa N)^2};q\Big)
            p_N^{(\mathsf{a},\mathsf{b})}
            \Big(\frac{v}{(\kappa N)^2};q\Big)
            \Big] 
            \\
            &=  \frac{q^{1-N}(1-q^{2N+\mathsf{a}+\mathsf{b}})}
                 {1-q^{\mathsf{a}+1}}
            \frac{f_N(u)g_N(v)-g_N(u)f_N(v)}{u-v}.
    \end{align*} 

    By \eqref{eqn:hard edge polynomial}, applied with parameters
    $\mathsf{a}$ and $\mathsf{a}+1$, respectively, we have 
    \begin{equation*}
        \begin{split}
            f_N(u)
            &=
            \Gamma(d_1+1)2^{d_1}
            u^{-d_1/2}J_{d_1}(\sqrt{u})+o(1),
            \\
            g_N(u)
            &=
            \Gamma(d_1+2)2^{d_1+1}
            u^{-d_1/2}\sqrt{u}\,J_{d_1+1}(\sqrt{u})+o(1),
        \end{split}
    \end{equation*}
    locally uniformly in $u$, together with derivatives of any fixed
    order. To control the quotient also when $u$ and $v$ approach
    each other, observe that
    \begin{equation*}
        \begin{split}
            \frac{f_N(u)g_N(v)-g_N(u)f_N(v)}{u-v}
            &=
            g_N(u)\int_0^1 f_N'\big(v+t(u-v)\big)\,dt
           -
            f_N(u)\int_0^1 g_N'\big(v+t(u-v)\big)\,dt.
        \end{split}
    \end{equation*}
    The right-hand side is defined also at $u=v$.
    The uniform convergence of the functions and their first
    derivatives therefore yields uniform convergence of this
    quotient on $[0,M]^2$.

    Using
    \begin{equation*}
        \frac{q^{1-N}(1-q^{2N+\mathsf{a}+\mathsf{b}})}
             {1-q^{\mathsf{a}+1}}
        =
        \frac{N}{\lambda}
        \frac{1-\mathsf{s}^{c_2+2}}
             {\mathsf{s}(d_1+1)}
        \big(1+O(N^{-1})\big)
    \end{equation*}
    and the Bessel identity \cite[Eq.~(10.6.2)]{NIST}
    \begin{equation*}
        J_\nu'(z)
        =
        -J_{\nu+1}(z)+\frac{\nu}{z}J_\nu(z),
    \end{equation*}
    we obtain
    \begin{equation}\label{eqn:hard wronskian asymptotic}
        \begin{split}
            &\quad \frac{q^{1-N}(1-q^{2N+\mathsf{a}+\mathsf{b}})}
                  {1-q^{\mathsf{a}+1}}
            \frac{f_N(u)g_N(v)-g_N(u)f_N(v)}{u-v}
            \\
            & =
            -\frac{N}{\lambda}
            \frac{1-\mathsf{s}^{c_2+2}}{\mathsf{s}}
            \Gamma(d_1+1)^2 2^{2d_1+2}
            \Big(
            (uv)^{-d_1/2}K_{d_1}^{\textup{Bessel}}(u,v)+o(1)
            \Big).
        \end{split}
    \end{equation}
    The error is uniform on $[0,M]^2$, where the weighted Bessel
    kernel is understood by its continuous extension when either
    variable is zero.

    The explicit formulae for the squared norms and leading
    coefficients give
    \begin{equation}\label{eqn:hard constant asymptotic}
        \begin{split}
            h_{N-1}^{(\mathsf{a},\mathsf{b})}
            &=
            \frac{\mathsf{s}^{d_1+1}\Gamma(d_1+1)^2}
            {(1-\mathsf{s})^{d_1}
             (1-\mathsf{s}^{c_2+1})^{d_1}
             (1-\mathsf{s}^{c_2+2})}
            \Big(\frac{\lambda}{N}\Big)^{2d_1+1}
            \big(1+O(N^{-1})\big),
            \\
            \frac{\gamma_{N-1}^{(\mathsf{a},\mathsf{b})}}
                 {\gamma_N^{(\mathsf{a},\mathsf{b})}}
            &=
            -\frac{\mathsf{s}(1-\mathsf{s})
                   (1-\mathsf{s}^{c_2+1})}
                  {(1-\mathsf{s}^{c_2+2})^2}
            +O(N^{-1}).
        \end{split}
    \end{equation}
    On the other hand, it follows from \eqref{eqn:def of qjacobi weight} that
    \begin{equation}\label{eqn:hard weight asymptotic}
        \Big(
        w^{(\mathsf{a},\mathsf{b})}
        \Big(\frac{u}{(\kappa N)^2};q\Big)
        w^{(\mathsf{a},\mathsf{b})}
        \Big(\frac{v}{(\kappa N)^2};q\Big)
        \Big)^{1/2}
        =
        \Big(\frac{uv}{(\kappa N)^4}\Big)^{d_1/2}
        \big(1+O(N^{-1})\big)
    \end{equation}
    uniformly for $0<u,v\le M$.

    Substituting \eqref{eqn:hard wronskian asymptotic},
    \eqref{eqn:hard constant asymptotic} and
    \eqref{eqn:hard weight asymptotic} into the
    Christoffel--Darboux formula \eqref{eqn:CD kernel}, and using \eqref{def of kappa}, the proposition follows. 
\end{proof}

We now pass from the weighted uniform convergence of the correlation
kernels to convergence of the associated Fredholm determinants.
As in the proof of Proposition~\ref{Prop_duality Gnm q-JUE}, it is
convenient to regard the $q$-JUE as a determinantal point process on
$\mathbb{Z}_{\geq0}$ with respect to counting measure. Its correlation
kernel is
\begin{equation}\label{eqn:def of mathsfK}
    \mathsf{K}_N^{(\mathsf{a},\mathsf{b})}(k,l;q)
    =
    (1-q)q^{(k+l)/2}
    K_N^{(\mathsf{a},\mathsf{b})}(q^k,q^l;q).
\end{equation}
For a fixed $s>0$, define
\begin{equation}
    I_N(s)
    =
    \{k\in\mathbb{Z}_{\geq0}:(\kappa N)^2q^k\le s\}.
\end{equation}
Then we have the following.

\begin{prop}\label{prop:hard edge Fredholm}
    Let $q=e^{-\lambda/N}$ for a fixed $\lambda>0$, and let
    $\mathsf{a}=d_1>-1$ and $\mathsf{b}=c_2N+d_2$.
    Then, for each fixed $s>0$,
    \begin{equation} \label{eq for hard edge Fred conv in Prop}
        \lim_{N\to\infty}
        \det
        \big(\Id-\mathsf{K}_N^{(\mathsf{a},\mathsf{b})}\big)_{\ell^2(I_N(s))}
        =
        \det
        \big(\Id-K_{d_1}^{\textup{Bessel}}\big)_{L^2(0,s)}.
    \end{equation}
\end{prop}

\begin{proof}
    Write
    \begin{equation*}
        u_k=(\kappa N)^2q^k,
        \qquad
        H(u,v)=(uv)^{-d_1/2}K_{d_1}^{\textup{Bessel}}(u,v).
    \end{equation*}
    The power series of the Bessel function shows that $H$ extends
    continuously to $[0,s]^2$. 

    By Proposition~\ref{prop:hard edge uniform},
    \begin{equation}\label{eqn:hard edge mathsfK uniform}
        \mathsf{K}_N^{(\mathsf{a},\mathsf{b})}(k,l;q)
        =
        (1-q)(u_ku_l)^{(d_1+1)/2}
        \big(H(u_k,u_l)+o(1)\big),
    \end{equation}
    uniformly for $k,l\in I_N(s)$.
    The kernel $\mathsf{K}_N^{(\mathsf{a},\mathsf{b})}$ defines an
    orthogonal projection of rank $N$ on
    $\ell^2(\mathbb{Z}_{\geq0})$. Its compression to $\ell^2(I_N(s))$
    is therefore positive and trace-class. Since $H$ is bounded on
    $[0,s]^2$, there is a constant $C>0$ such that
    \begin{equation*}
        0\le
        \mathsf{K}_N^{(\mathsf{a},\mathsf{b})}(k,k;q)
        \le C(1-q)u_k^{d_1+1},
        \qquad k\in I_N(s),
    \end{equation*}
    for all sufficiently large $N$. Consequently,
    \begin{equation}\label{eqn:hard diagonal sum}
        \begin{split}
            \sum_{k\in I_N(s)}
            \mathsf{K}_N^{(\mathsf{a},\mathsf{b})}(k,k;q)
            &\le C(1-q)\sum_{k\in I_N(s)}u_k^{d_1+1}
           =
            C\frac{1-q}{1-q^{d_1+1}}
            u_{\min I_N(s)}^{d_1+1}
            \le C_s,
        \end{split}
    \end{equation}
    where $C_s$ is independent of $N$. Here, we used
    $u_{\min I_N(s)}\le s$. 

    We next prove convergence of each term in the Fredholm expansion.
    For a fixed $m\ge1$, factoring the powers of $u_{k_j}$ from the
    rows and columns in \eqref{eqn:hard edge mathsfK uniform} gives
    \begin{equation*}
        \begin{split}
        \det\big[
            \mathsf{K}_N^{(\mathsf{a},\mathsf{b})}(k_j,k_l;q)
            \big]_{j,l=1}^m =
            \prod_{j=1}^m\big((1-q)u_{k_j}^{d_1+1}\big)
            \Big(
            \det\big[H(u_{k_j},u_{k_l})\big]_{j,l=1}^m
            +o(1)
            \Big),
        \end{split}
    \end{equation*}
    uniformly for $k_1,\ldots,k_m\in I_N(s)$.
    Since
    \begin{equation*}
        \sum_{k_1,\ldots,k_m\in I_N(s)}
        \prod_{j=1}^m\big((1-q)u_{k_j}^{d_1+1}\big)
        =
        \Big(
        \frac{1-q}{1-q^{d_1+1}}
        u_{\min I_N(s)}^{d_1+1}
        \Big)^m
    \end{equation*}
    is bounded independently of $N$, the error remains $o(1)$ after
    summation over all tuples.

    To identify the limit of the leading term, associate to each
    $k\in I_N(s)$ the interval $(u_{k+1},u_k]$. These intervals satisfy
    \begin{equation*}
        \bigcup_{k\in I_N(s)}(u_{k+1},u_k]
        =
        (0,u_{\min I_N(s)}],
        \qquad
        qs<u_{\min I_N(s)}\le s
    \end{equation*}
    for all sufficiently large $N$, and their maximal length is at
    most $(1-q)s$, which tends to zero. Furthermore,
    \begin{equation*}
        (1-q)u_k^{d_1+1}
        =
        \frac{(d_1+1)(1-q)}{1-q^{d_1+1}}
        \int_{u_{k+1}}^{u_k}u^{d_1}\,du,
    \end{equation*}
    where the prefactor tends to one. Thus these sums are Riemann
    sums with respect to the finite measure $u^{d_1}\,du$ on $(0,s)$.
    The continuity of $H$ on $[0,s]^2$ therefore yields
    \begin{equation}\label{eqn:hard determinant convergence}
        \begin{split}
            \lim_{N\to\infty}
            \sum_{k_1,\ldots,k_m\in I_N(s)}
            \det\big[
            \mathsf{K}_N^{(\mathsf{a},\mathsf{b})}(k_j,k_l;q)
            \big]_{j,l=1}^m
            &=
            \int_{(0,s)^m}
            \det\big[H(u_j,u_l)\big]_{j,l=1}^m
            \prod_{j=1}^m u_j^{d_1}\,du_j
            \\
            &=
            \int_{(0,s)^m}
            \det\big[K_{d_1}^{\textup{Bessel}}(u_j,u_l)\big]_{j,l=1}^m
            \,du_1\cdots du_m.
        \end{split}
    \end{equation}

    It remains to interchange the limit with the Fredholm series.
    We have
    \begin{equation}\label{eqn:hard Fredholm expansion}
        \begin{split}
       \det 
            \big(\Id-\mathsf{K}_N^{(\mathsf{a},\mathsf{b})}\big)_{\ell^2(I_N(s))}
        =
            \sum_{m=0}^{\infty}\frac{(-1)^m}{m!}
            \sum_{k_1,\ldots,k_m\in I_N(s)}
            \det\big[
            \mathsf{K}_N^{(\mathsf{a},\mathsf{b})}(k_j,k_l;q)
            \big]_{j,l=1}^m,
        \end{split}
    \end{equation}
    with the $m=0$ term understood to be one.
    Positivity and Hadamard's inequality imply
    \begin{equation*}
        0\le
        \det\big[
        \mathsf{K}_N^{(\mathsf{a},\mathsf{b})}(k_j,k_l;q)
        \big]_{j,l=1}^m
        \le
        \prod_{j=1}^m
        \mathsf{K}_N^{(\mathsf{a},\mathsf{b})}(k_j,k_j;q).
    \end{equation*}
    Hence, by \eqref{eqn:hard diagonal sum},
    \begin{equation*}
        \begin{split}
       \frac{1}{m!}
            \sum_{k_1,\ldots,k_m\in I_N(s)}
            \det\big[
            \mathsf{K}_N^{(\mathsf{a},\mathsf{b})}(k_j,k_l;q)
            \big]_{j,l=1}^m
       \le
            \frac{1}{m!}
            \Big(
            \sum_{k\in I_N(s)}
            \mathsf{K}_N^{(\mathsf{a},\mathsf{b})}(k,k;q)
            \Big)^m
            \le \frac{C_s^m}{m!}.
        \end{split}
    \end{equation*}
    This bound is summable over $m$ and independent of $N$.
    Applying dominated convergence to
    \eqref{eqn:hard Fredholm expansion} and using
    \eqref{eqn:hard determinant convergence}, we obtain the desired convergence \eqref{eq for hard edge Fred conv in Prop}.
\end{proof}

\begin{proof}[Proof of Theorem~\ref{Thm_qJUE scaling limits}~(i)]
    Fix $s>0$. By the definition of $I_N(s)$, the event
    $x_{\min}>s/(\kappa N)^2$ is precisely the event that there
    are no particles at the sites $q^k$ with $k\in I_N(s)$.
    Thus, the determinantal gap identity gives
    \begin{equation*}
        \mathbb{P}\Big(x_{\min}>\frac{s}{(\kappa N)^2}\Big)
        =
        \det
        \big(\Id-\mathsf{K}_N^{(\mathsf{a},\mathsf{b})}\big)_{\ell^2(I_N(s))}.
    \end{equation*}
    Taking $N\to\infty$ and applying
    Proposition~\ref{prop:hard edge Fredholm}, we obtain \eqref{hard edge scaling limit main}. 
\end{proof}

\subsection{Soft edge scaling limit}

We now turn to the soft-edge regime $c_1>0$ and begin the proof of Theorem~\ref{Thm_qJUE scaling limits}~(ii) by establishing local uniform convergence of the rescaled correlation kernel to the Airy kernel. Let $a$ be the left endpoint defined in \eqref{eqn:def of endpoints ab}, and let $\theta$ be the corresponding soft-edge scaling constant \eqref{def of theta}. For $u,v\in\mathbb{R}$, set
\begin{equation}
    \widehat{x}
    =a-\frac{u}{(\theta N)^{2/3}},
    \qquad
    \widehat{y}
    =a-\frac{v}{(\theta N)^{2/3}}.
\end{equation}

\begin{prop}\label{prop:soft kernel uniform} Let $q=e^{-\lambda/N}$ for a fixed $\lambda>0$, and fix $c_1>0$ and $c_2\geq0$, with $d_2>-1$ when $c_2=0$. For a fixed $M>0$, define
    \begin{equation}
        \mathfrak{L}_{M}^{\textup{soft}}
        =
        \Big\{u\in\mathbb{R}:|u|<M,\,
        a-\frac{u}{(\theta N)^{2/3}}\in L_q\Big\},
    \end{equation}
     where $L_q=\{1,q,q^2,\ldots\}$. Then, as $N\to\infty$,
    \begin{equation}
        \sup_{u,v\in\mathfrak{L}_{M}^{\textup{soft}}}
        \Big|
        \frac{1}{(\theta N)^{2/3}}
        K_N^{(\mathsf{a},\mathsf{b})}
        (\widehat{x},\widehat{y};q)
        -
        K_{\Ai}(u,v)
        \Big|
        \to 0.
    \end{equation}
\end{prop}

To prove Proposition~\ref{prop:soft kernel uniform}, we apply the method of steepest descent to a contour integral representation of the little $q$-Jacobi polynomials. For $n\in\mathbb{Z}_{\geq0}$, define
\begin{equation}\label{eqn:def of mathfrakF}
    \mathfrak{F}_{n}(x,z)
    :=z^{n}\frac{(-q/z,-q^{n+\mathsf{a}+\mathsf{b}+2}x/z;q)_{\infty}}{(-qx/z,-q^{n+\mathsf{a}+1}/z;q)_{\infty}}.
\end{equation}
For fixed $x\in(0,1)$, the possible poles of $\mathfrak{F}_n(x,z)$
as a function of $z$, away from the origin, are
\begin{equation*}
    -q^{n+\mathsf{a}+1+j},
    \qquad
    -xq^{1+j},
    \qquad
    j\in\mathbb{Z}_{\geq0}.
\end{equation*}
The following representation provides the starting point for
our asymptotic analysis.

\begin{lem}
    For $n\in\mathbb{Z}_{\geq0}$ and $x\in(0,1)$, we have the contour integral representation
    \begin{equation}\label{eqn:littleq integral}
        p_n^{(\mathsf{a},\mathsf{b})}(x;q)
        =
        \frac{(q;q)_n}{(q^{\mathsf{a}+1};q)_n}
        q^{-n(n+1)/2}
        \oint_{|z|=\mathsf{r}}
        \mathfrak{F}_n(x,z)\frac{dz}{2\pi i z},
    \end{equation}
    where the circle is positively oriented and its radius satisfies $ \max\{q^{n+\mathsf{a}+1},qx\}<\mathsf{r}<1.$ 
\end{lem}

\begin{proof}
    We use the coefficient integral
    \cite[Eq.~(3.3)]{IZ16}
    \begin{equation}
        \frac{(A;q)_k}{(B;q)_k}
        =
        \frac{(q,B/A;q)_\infty}{(B,q/A;q)_\infty}
        \oint_{|z|=\mathsf{r}}
        \frac{(Az,q/(Az);q)_\infty}
             {(z,B/(Az);q)_\infty}
        \frac{dz}{2\pi i z^{k+1}},
        \qquad
        |B/A|<\mathsf{r}<1.
    \end{equation}
    We apply this identity with $A=q^{-n}$ and
    $B=q^{\mathsf{a}+1}$, choosing $  \max\{q^{n+\mathsf{a}+1},qx\}<\mathsf{r}<1.$ 
    Since $(q^{-n};q)_k=0$ for $k>n$, the terminating series
    \eqref{eqn:littleq jacobi hypergeometric} may be written as
    a sum over all $k\ge0$. After substituting the coefficient
    integral, the resulting series converges absolutely and
    uniformly on $|z|=\mathsf{r}$ because $|qx/z|<1$.
    We may therefore interchange the sum and the integral.
    Simplifying the $q$-Pochhammer factors gives
    \begin{equation*}
        \begin{split}
            p_n^{(\mathsf{a},\mathsf{b})}(x;q)
            &=
            \frac{(q;q)_n}{(q^{\mathsf{a}+1};q)_n}
            \oint_{|z|=\mathsf{r}}
            \frac{(q^{-n}z;q)_n(q^{n+1}/z;q)_\infty}
                 {(q^{n+\mathsf{a}+1}/z;q)_\infty} 
            \sum_{k=0}^{\infty}
            \frac{(q^{n+\mathsf{a}+\mathsf{b}+1};q)_k}
                 {(q;q)_k}
            \Big(\frac{qx}{z}\Big)^k
            \frac{dz}{2\pi i z}
            \\
            &=
            \frac{(q;q)_n}{(q^{\mathsf{a}+1};q)_n}
            \oint_{|z|=\mathsf{r}}
            \frac{(q^{-n}z;q)_n
                  (q^{n+1}/z,q^{n+\mathsf{a}+\mathsf{b}+2}x/z;q)_\infty}
                 {(q^{n+\mathsf{a}+1}/z,qx/z;q)_\infty}
            \frac{dz}{2\pi i z},
        \end{split}
    \end{equation*}
    where the second equality follows from the $q$-binomial theorem.
    Finally, using
    $$
   (q^{-n}z;q)_n
          =
            (-z)^n q^{-n(n+1)/2}(q/z;q)_n, \qquad     (q/z;q)_n(q^{n+1}/z;q)_\infty
             =
            (q/z;q)_\infty,
    $$ 
    and making the change of variables $z\mapsto-z$, we obtain
    \eqref{eqn:littleq integral}.
\end{proof}

In view of \eqref{eqn:littleq integral}, we introduce the normalised polynomial
\begin{equation}\label{eqn:mathfrak p integral}
    \mathfrak{p}_{n}(x)
    :=\oint_{|z|=\mathsf{r}}\mathfrak{F}_{n}(x,z)\frac{dz}{2\pi iz}
    =\frac{(q^{\mathsf{a}+1};q)_{n}}{(q;q)_{n}}q^{\frac{n(n+1)}{2}}p_{n}^{(\mathsf{a},\mathsf{b})}(x;q), 
\end{equation}
where $\mathsf{r}$ is any admissible radius in \eqref{eqn:littleq integral}. The Christoffel--Darboux formula \eqref{eqn:CD kernel} then becomes
\begin{equation} \label{eqn:CD kernel v2}
    K_{N}^{(\mathsf{a},\mathsf{b})}(x,y;q)  =\mathsf{C}_{N}^{(\mathsf{a},\mathsf{b})}\sqrt{w^{(\mathsf{a},\mathsf{b})}(x;q)w^{(\mathsf{a},\mathsf{b})}(y;q)}\frac{\mathfrak{p}_{N-1}(x)\mathfrak{p}_{N}(y)-\mathfrak{p}_{N}(x)\mathfrak{p}_{N-1}(y)}{x-y},
\end{equation}
with the positive normalisation constant 
\begin{equation*}
    \mathsf{C}_{N}^{(\mathsf{a},\mathsf{b})}
    =\frac{q^{-N^{2}-\mathsf{a}(N-1)}}{(1-q)(1-q^{2N+\mathsf{a}+\mathsf{b}})}\frac{(q^{\mathsf{a}+N},q^{\mathsf{b}+N};q)_{\infty}}{(q^{N+1},q^{\mathsf{a}+\mathsf{b}+N+1};q)_{\infty}}.
\end{equation*}

To analyse the contour integral, we introduce the phase and amplitude functions
\begin{align}
    \Phi(x,z)
    &:=\log z+\frac{1}{\lambda}\Big(\Li(-x/z)+\Li(-\mathsf{s}^{c_{1}+1}/z)-\Li(-1/z)-\Li(-\mathsf{s}^{c_{1}+c_{2}+1}x/z)\Big), \label{def of Phi}
    \smallskip\\
    \Psi(x,z)
    &:=\Big(\frac{x+z}{1+z}\Big)^{\frac{1}{2}}\Big(1+\frac{\mathsf{s}^{c_{1}+1}}{z}\Big)^{d_{1}+\frac{1}{2}}\Big(1+\frac{\mathsf{s}^{c_{1}+c_{2}+1}x}{z}\Big)^{-(d_{1}+d_{2}+\frac{3}{2})}. \label{def of Psi}
\end{align}
The branches are chosen so that $\Phi$ is real and $\Psi$ is positive for $x,z>0$, and are continued analytically near the saddle point introduced below. These functions describe the leading exponential and multiplicative factors in the asymptotics of $\mathfrak{F}_N$.
We use subscripts to denote partial derivatives; for example,
$\Phi_z=\partial_z\Phi$, $\Phi_{zz}=\partial_z^2\Phi$, and
$\Phi_{xz}=\partial_x\partial_z\Phi$.

It is also convenient to introduce the rational function
\begin{equation*}
    \phi(x,z)
    :=\frac{(x+z)(z+\mathsf{s}^{c_{1}+1})}{(1+z)(z+\mathsf{s}^{c_{1}+c_{2}+1}x)}. 
\end{equation*}
Direct differentiation gives 
\begin{equation}\label{eqn:Phizphi}
    \mathsf{s}\exp[\lambda z\Phi_{z}(x,z)]
    =\phi(x,z).
\end{equation}
In particular, this identity allows us to locate the relevant saddle
point explicitly. Set
\begin{equation} \label{def of zstar}
     z_*
    :=
    \frac{\mathsf{s}(1-\mathsf{s}^{c_1})
          -a(1-\mathsf{s}^{c_1+c_2+2})}
         {2(1-\mathsf{s})}.
\end{equation}
The following lemma identifies the double critical point at the left
edge and relates the local derivatives of the phase to the scaling
constant $\theta$.

\begin{lem}\label{lem:zast relation}
    We have $z_*>0$ and
    \begin{equation}
        \Phi_z(a,z_*)=\Phi_{zz}(a,z_*)=0.
    \end{equation} 
     Moreover,
    \begin{equation}\label{eqn:zast2}
        z_*^2
        =
        \frac{\mathsf{s}^{c_1+1}a(1-\mathsf{s}^{c_2+1})}
             {1-\mathsf{s}},
    \end{equation}
    and the relevant derivatives satisfy
    \begin{equation}
        \Phi_{xz}(a,z_*)
        =
        \frac{1-\mathsf{s}^{c_1+c_2+1}}
             {\lambda(z_*+a)(z_*+\mathsf{s}^{c_1+c_2+1}a)}
        >0,
    \end{equation}
    \begin{equation}\label{eqn:Phizzz}
        \Phi_{zzz}(a,z_*)
        =
        \frac{2(1-\mathsf{s})}
             {\lambda\mathsf{s}z_*(z_*+1)
              (z_*+\mathsf{s}^{c_1+c_2+1}a)}
        >0.
    \end{equation}
    Finally, the soft-edge scaling constant $\theta$ in \eqref{def of theta} is related to these
    derivatives by
    \begin{equation}\label{eqn:theta2}
        \theta^2
        =
        \frac{2\Phi_{xz}(a,z_*)^3}{\Phi_{zzz}(a,z_*)}.
    \end{equation}
\end{lem}

\begin{proof}
    For $x,z>0$, \eqref{eqn:Phizphi} shows that the saddle equation
    $\Phi_z(x,z)=0$ is equivalent to $\phi(x,z)=\mathsf{s}$.
    Clearing denominators gives a quadratic equation in $z$ with discriminant
    \[
        (1-\mathsf{s}^{c_1+c_2+2})^2(x-a)(x-b).
    \]
    By \eqref{eqn:def of endpoints ab}, at $x=a$ this equation
    has the positive double root $z=z_\ast$, and hence
    \begin{equation}\label{eqn:phi-s}
        \phi(a,z)-\mathsf{s}
        =\frac{(1-\mathsf{s})(z-z_\ast)^2}
        {(z+1)(z+\mathsf{s}^{c_1+c_2+1}a)}.
    \end{equation}
    Comparing constant terms gives \eqref{eqn:zast2}.
    The derivative identities follow by differentiating
    \eqref{eqn:Phizphi} and using \eqref{eqn:phi-s}.
    Finally, \eqref{eqn:theta2} follows by direct substitution
    using \eqref{eqn:def of endpoints ab} and \eqref{def of theta}.
\end{proof}

We next derive the local asymptotic expansion of the integrand.

\begin{lem}\label{lem:mathfrakF asymp}
    Let $q=e^{-\lambda/N}$ for a fixed $\lambda>0$.
    Uniformly on a sufficiently small fixed compact neighbourhood
    of $(a,z_\ast)$ in $\mathbb{C}^2$, we have
    \begin{equation}\label{eqn:mathfrak F asymp}
        \mathfrak{F}_{N}(x,z)
        =\Psi(x,z)e^{N\Phi(x,z)}(1+O(N^{-1})).
    \end{equation}
    The relative error and all its partial derivatives of any
    fixed order are uniformly $O(N^{-1})$.
\end{lem}

\begin{proof}
    For fixed $d\in\mathbb{R}$, the Euler--Maclaurin formula gives
    \begin{equation}\label{eqn:qPochhammer asymptotic}
        \log(zq^d;q)_\infty
        =-\frac{N}{\lambda}\Li(z)
        +\Big(\frac12-d\Big)\log(1-z)+O(N^{-1}),
    \end{equation}
    uniformly on compact subsets of $\mathbb{C}\setminus[1,\infty)$;
    see also \cite{Mc99}.
    Applying this expansion to the four $q$-Pochhammer symbols
    in \eqref{eqn:def of mathfrakF} yields
    \eqref{eqn:mathfrak F asymp}.
    The derivative bounds follow from Cauchy's estimates applied
    to the analytic relative error on a slightly larger neighbourhood.
\end{proof}

We next construct a contour for the steepest descent analysis.
We use the principal branches of the logarithm and dilogarithm; see \cite[Section~25.12(i)]{NIST}.
On the negative real axis, $\Re\Phi(a,z)$ denotes the common boundary value from the upper and lower half-planes, which agrees by \cite[Eq.~(25.12.4)]{NIST}. Recall that $z_*$ is given by \eqref{def of zstar}.

\begin{lem}\label{lem:soft contour}
Assume $c_1>0$ and $c_2\geq0$.
There exists a simple, positively oriented, piecewise smooth
closed contour $\Gamma_{\textup{s}}$ with the following properties:
    \begin{enumerate}
        \item $\Gamma_{\textup{s}}$ strictly encloses
        $[-\mathsf{s}^{c_1+1},0]$ and $[-a,0]$.
    \smallskip 
    \item For some $\varepsilon>0$,
        \begin{equation*}
            \Gamma_{\textup{s}}\cap\{z:|z-z_\ast|<\varepsilon\}
            =\{z_\ast+te^{\pm\pi i/3}:0\leq t<\varepsilon\}.
        \end{equation*}
        \item $\Re\Phi(a,z)<\Phi(a,z_\ast)$ for
        $z\in\Gamma_{\textup{s}}\setminus\{z_\ast\}$.
    \smallskip
    \item $|z|\geq z_\ast$ for $z\in\Gamma_{\textup{s}}$.
    \end{enumerate}
\end{lem}

\begin{proof} We first assume $c_2>0$ and deform the circle $|z|=z_{\ast}$. 
    By \eqref{eqn:zast2},
    \begin{equation*}
        z_\ast^2-\mathsf{s}^{c_1+c_2+1}a
        =\frac{\mathsf{s}^{c_1+1}a(1-\mathsf{s}^{c_2})}
        {1-\mathsf{s}}>0.
    \end{equation*}
    Substituting $z=z_\ast e^{it}$ into \eqref{eqn:phi-s}
    therefore gives $\Im\phi(a,z_\ast e^{it})>0$ for $0<t<\pi$.
    Hence, by \eqref{eqn:Phizphi},
    \begin{equation*}
        \frac{d}{dt}\Re\Phi(a,z_\ast e^{it})
        =-\frac{1}{\lambda}\arg\phi(a,z_\ast e^{it})<0.
    \end{equation*}
    By continuity and conjugation symmetry, the circle satisfies
    the third property.

    Write $R=\max\{\mathsf{s}^{c_1+1},a\}$.
    If $z_\ast>R$, the first property also holds.
    Otherwise, \eqref{eqn:zast2} implies
    \[
        \min\{\mathsf{s}^{c_1+1},a\}<z_\ast\leq R<1,
        \qquad
        \mathsf{s}^{c_1+c_2+1}a<z_\ast.
    \]
    For $z_\ast<r<R$, \eqref{eqn:phi-s} then gives
    $0<\phi(a,-r)<\mathsf{s}$, so
    \begin{equation*}
        \frac{d}{dr}\Re\Phi(a,-r)
        =\frac{1}{\lambda r}
        \log\frac{\phi(a,-r)}{\mathsf{s}}<0.
    \end{equation*}
    By continuity, a thin outward deformation near the negative
    real axis can thus extend slightly beyond $-R$, ensuring
    the first property while preserving the third.

    Finally, near $z_\ast$, replace the circular arc by the rays
    $z_\ast+te^{\pm\pi i/3}$, $0\leq t\leq\varepsilon$,
    joined to the original circle by arcs of
    $|z-z_\ast|=\varepsilon$.
    The upper connecting arc is parametrised by
    $z=z_\ast+\varepsilon e^{it}$, where
    \[
        \frac{\pi}{3}\leq t\leq t_\varepsilon,
        \qquad
        t_\varepsilon
        =\arccos\Big(-\frac{\varepsilon}{2z_\ast}\Big)
        =\frac{\pi}{2}+O(\varepsilon).
    \]
    Since $\Phi_z(a,z_\ast)=\Phi_{zz}(a,z_\ast)=0$ and
    $\Phi_{zzz}(a,z_\ast)>0$, Taylor expansion gives
    \begin{equation*}
    \begin{split}
        \Re\Phi(a,z_\ast+\varepsilon e^{it})-\Phi(a,z_\ast)
        &=\frac{\varepsilon^3}{6}\Phi_{zzz}(a,z_\ast)
        \cos(3t)+O(\varepsilon^4),\\
        \frac{d}{dt}\Re\Phi(a,z_\ast+\varepsilon e^{it})
        &=-\frac{\varepsilon^3}{2}\Phi_{zzz}(a,z_\ast)
        \sin(3t)+O(\varepsilon^4).
    \end{split}
    \end{equation*}
    For a fixed sufficiently small $\delta>0$, the first expression
    is negative on $[\pi/3,\pi/2-\delta]$, whereas the second
    is positive on $[\pi/2-\delta,t_\varepsilon]$.
    The latter part of the arc is therefore bounded above by
    its endpoint value on $|z|=z_\ast$, which is strictly below
    $\Phi(a,z_\ast)$.
    The same Taylor expansion gives strict descent along the rays.
    Reflecting the construction across the real axis establishes
    the second and third properties.
    All modifications lie in $|z|\geq z_\ast$, proving the fourth.

    It remains to consider $c_2=0$. In this case,
\eqref{eqn:zast2} gives $z_{\ast}^{2}=\mathsf{s}^{c_1+1}a$, and
\[
    0\leq\phi(a,z_{\ast}e^{it})
    =\frac{a+\mathsf{s}^{c_1+1}+2z_{\ast}\cos t}
    {1+z_{\ast}^{2}+2z_{\ast}\cos t}
    <\mathsf{s},
    \qquad 0<t<2\pi.
\]
Thus $\Re\Phi(a,z)$ is constant on $|z|=z_{\ast}$, while
\eqref{eqn:Phizphi} shows that a small outward radial
deformation away from $z_{\ast}$ decreases its real part.
Combining this deformation with the local Airy rays and
connecting arcs constructed above, and with the deformation
near the negative real axis when necessary, gives a contour
satisfying all four properties.
\end{proof}

To justify exponential localisation near $z_\ast$, we first note that,
uniformly for $z\in\Gamma_{\textup{s}}$,
\begin{equation*}
    \frac{1}{N}\log|(-q/z;q)_\infty|
    \leq \frac{1}{\lambda}\int_0^\infty
    \log|1+e^{-t}/z|\,dt+o(1)
    =-\frac{1}{\lambda}\Re\Li(-1/z)+o(1).
\end{equation*}
Indeed, the upper bound follows by taking uniform Riemann sums after replacing the logarithm by $$\max\{\log|1+e^{-t}/z|,-L\},$$ and then letting $L\to\infty$.
The integral identity follows from \cite[Eq.~(25.12.2)]{NIST}.
For the other three $q$-Pochhammer factors in
\eqref{eqn:def of mathfrakF},
\eqref{eqn:qPochhammer asymptotic} applies uniformly on the contour
for $x$ near $a$. Consequently,
$$\log|\mathfrak{F}_N(x,z)|\leq N\Re\Phi(x,z)+o(N)$$ uniformly there.
By compactness and the third property of
Lemma~\ref{lem:soft contour}, for any fixed small neighbourhood
$U$ of $z_\ast$ and $M>0$, there exists $C_{\textup{out}}>0$ such that
\begin{equation}\label{eqn:soft out bound}
    \sup_{\substack{u\in\mathbb{C}\\ |u|\leq M}}
    \sup_{z\in\Gamma_{\textup{s}}\setminus U}
    \big|e^{-N\Phi(\widehat{x},z_\ast)}
    \mathfrak{F}_N(\widehat{x},z)\big|
    \leq e^{-C_{\textup{out}}N}
\end{equation}
for all sufficiently large $N$.

By the first property of Lemma~\ref{lem:soft contour},
$\Gamma_{\textup{s}}$ encloses all poles of both
$\mathfrak{F}_N(x,\cdot)$ and $\mathfrak{F}_{N-1}(x,\cdot)$
for $x$ near $a$ and sufficiently large $N$.
It is therefore admissible in the integral representations of
$\mathfrak{p}_N$ and $\mathfrak{p}_{N-1}$.
To remove their common leading saddle contribution, observe that
\begin{equation} \label{def of mathfrak R}
    \mathfrak{F}_{N-1}(x,z)
    =\mathfrak{F}_N(x,z)\mathfrak{R}_N(x,z),
    \qquad
    \mathfrak{R}_N(x,z)
    :=\frac{z+q^{N+\mathsf{a}+\mathsf{b}+1}x}
    {z(z+q^{N+\mathsf{a}})}.
\end{equation}
Near $(a,z_\ast)$, we have
\begin{equation*}
    \mathfrak{R}_N(x,z)
    =\mathfrak{R}(x,z)+O(N^{-1}),
    \qquad
    \mathfrak{R}(x,z)
    :=\frac{z+\mathsf{s}^{c_1+c_2+1}x}
    {z(z+\mathsf{s}^{c_1+1})},
\end{equation*}
locally uniformly together with all partial derivatives of any
fixed order.
Define
\begin{equation}\label{eqn:mathfrak q integral}
\begin{split}
    \mathfrak{q}_N(x)
    &:=\mathfrak{p}_{N-1}(x)
    -\mathfrak{R}_N(x,z_\ast)\mathfrak{p}_N(x)\\
    &=\oint_{\Gamma_{\textup{s}}}
    \mathfrak{F}_N(x,z)
    \big(\mathfrak{R}_N(x,z)-\mathfrak{R}_N(x,z_\ast)\big)
    \frac{dz}{2\pi iz}.
\end{split}
\end{equation}
The Christoffel--Darboux formula \eqref{eqn:CD kernel v2} then becomes
\begin{equation}\label{eqn:soft kernel pq}
\begin{split}
    &\frac{K_N^{(\mathsf{a},\mathsf{b})}(x,y;q)}
    {\mathsf{C}_N^{(\mathsf{a},\mathsf{b})}
    \sqrt{w^{(\mathsf{a},\mathsf{b})}(x;q)
    w^{(\mathsf{a},\mathsf{b})}(y;q)}} =
    \frac{\mathfrak{q}_N(x)\mathfrak{p}_N(y)
    -\mathfrak{p}_N(x)\mathfrak{q}_N(y)}{x-y}
    +\frac{q^{N+\mathsf{a}+\mathsf{b}+1}}
    {z_\ast(z_\ast+q^{N+\mathsf{a}})}
    \mathfrak{p}_N(x)\mathfrak{p}_N(y),
\end{split}
\end{equation}
with the diagonal understood by continuity.

\begin{lem}\label{lem:soft pq asymptotic}
Under the assumptions of Proposition~\ref{prop:soft kernel uniform},
as $N\to\infty$, 
    \begin{align}
        N^{1/3}e^{-N\Phi(\widehat{x},z_\ast)}
        \mathfrak{p}_N(\widehat{x})
        &=\frac{\Psi(a,z_\ast)}{z_\ast}
        \Big(\frac{2}{\Phi_{zzz}(a,z_\ast)}\Big)^{1/3}
        \Ai(u)+o(1),\\
        N^{2/3}e^{-N\Phi(\widehat{x},z_\ast)}
        \mathfrak{q}_N(\widehat{x})
        &=-\frac{\Psi(a,z_\ast)\mathfrak{R}_z(a,z_\ast)}{z_\ast}
        \Big(\frac{2}{\Phi_{zzz}(a,z_\ast)}\Big)^{2/3}
        \Ai'(u)+o(1).
    \end{align}
    Both error terms and their derivatives of any fixed order
    tend to zero locally uniformly for $u\in\mathbb{C}$.
\end{lem}

\begin{proof}
    Fix $M>0$ and a sufficiently small neighbourhood $U$ of $z_\ast$.
    By \eqref{eqn:soft out bound} and the uniform boundedness of
    $\mathfrak{R}_N$ on $\Gamma_{\textup{s}}$, the contributions
    from $\Gamma_{\textup{s}}\setminus U$ are exponentially small
    after normalisation, uniformly for $|u|\leq M$.

    On the local contour, make the change of variables
    \begin{equation*}
        z=z_\ast+\Big(\frac{2}{N\Phi_{zzz}(a,z_\ast)}\Big)^{1/3}w.
    \end{equation*}
    Lemmas~\ref{lem:zast relation} and \ref{lem:mathfrakF asymp},
    in particular \eqref{eqn:theta2}, give
    \begin{equation}\label{eqn:mathfrak F saddle expansion}
        \mathfrak{F}_N(\widehat{x},z)
        =\Psi(a,z_\ast)
        \exp\Big[N\Phi(\widehat{x},z_\ast)+\frac{w^3}{3}-uw\Big]
        \big(1+O(N^{-1/3})\big),
    \end{equation}
    uniformly for $|u|\leq M$ and $w$ in any fixed compact set.
    In the same range, we have 
    $$
    N^{1/3}\frac{dz}{z} \to  \Big(\frac{2}{\Phi_{zzz}(a,z_\ast)}\Big)^{1/3}
        \frac{dw}{z_\ast}, \qquad   N^{1/3}\big(\mathfrak{R}_N(\widehat{x},z)
        -\mathfrak{R}_N(\widehat{x},z_\ast)\big)
    \to 
        \Big(\frac{2}{\Phi_{zzz}(a,z_\ast)}\Big)^{1/3}
        \mathfrak{R}_z(a,z_\ast)w.
    $$ 

    To justify integration over the expanding contour, note that
    Taylor expansion along the local rays gives
    \begin{equation*}
        N\Re\big(\Phi(\widehat{x},z)-\Phi(\widehat{x},z_\ast)\big)
        \leq -c|w|^3+C_M|w|,
    \end{equation*}
    for some $c,C_M>0$ independent of $N$, uniformly for $|u|\leq M$.
    Here the cubic term is negative along the rays, while the
    perturbation $\widehat{x}-a=O(N^{-2/3})$ contributes at most
    $C_M|w|$.
    Moreover,
    $$N^{1/3}|\mathfrak{R}_N(\widehat{x},z)
    -\mathfrak{R}_N(\widehat{x},z_\ast)|\leq C|w|.$$
    Thus, by Lemma~\ref{lem:mathfrakF asymp}, both rescaled
    integrands admit an integrable majorant of the form
    \[
        C(1+|w|)\exp\big(-c|w|^3+C_M|w|\big).
    \]

    Dominated convergence therefore reduces the two limits to
    integrals over the Airy contour $\Gamma_{\Ai}$, oriented from
    $\infty e^{-\pi i/3}$ through $0$ to $\infty e^{\pi i/3}$.
    By \cite[Eq.~(9.5.4)]{NIST} and differentiation under the
    integral sign,
    \begin{equation*}
        \Ai(u)=\int_{\Gamma_{\Ai}}
        e^{w^3/3-uw}\frac{dw}{2\pi i},
        \qquad
        \Ai'(u)=-\int_{\Gamma_{\Ai}}
        w e^{w^3/3-uw}\frac{dw}{2\pi i}.
    \end{equation*}
    Substitution into \eqref{eqn:mathfrak p integral} and
    \eqref{eqn:mathfrak q integral} proves the asserted limits
    locally uniformly for $u\in\mathbb{C}$.
    Analyticity in $u$ and Cauchy's integral formula then give
    convergence of all derivatives of any fixed order.
\end{proof}

We now complete the proof of Proposition~\ref{prop:soft kernel uniform}.

\begin{proof}[Proof of Proposition~\ref{prop:soft kernel uniform}]
    Write $w(x;q)=w^{(\mathsf{a},\mathsf{b})}(x;q)$.
    By \eqref{eqn:qPochhammer asymptotic},
    \begin{equation}\label{eqn:soft Cn asymptotic}
    \begin{split}
        \mathsf{C}_N^{(\mathsf{a},\mathsf{b})}
        &=\frac{N\mathsf{s}^{c_1-d_1}}
        {\lambda(1-\mathsf{s}^{c_1+c_2+2})}
        \exp\Big[-\frac{N}{\lambda}
        \big(\Li(\mathsf{s}^{c_1+1})
        +\Li(\mathsf{s}^{c_2+1})-\Li(\mathsf{s})
        -\Li(\mathsf{s}^{c_1+c_2+1})
        -\lambda^2(1+c_1)\big)\Big]\\
        &\qquad\times
        \frac{(1-\mathsf{s})^{1/2}
        (1-\mathsf{s}^{c_1+c_2+1})^{d_1+d_2+1/2}}
        {(1-\mathsf{s}^{c_1+1})^{d_1-1/2}
        (1-\mathsf{s}^{c_2+1})^{d_2-1/2}}
        (1+O(N^{-1})).
    \end{split}
    \end{equation}
    Similarly, uniformly for $x$ in a fixed neighbourhood of $a$,
    \begin{equation}\label{eqn:soft weight asymptotic}
        w(x;q)
        =\frac{x^{d_1}(1-\mathsf{s}^{c_2}x)^{d_2+1/2}}
        {(1-x)^{1/2}}e^{-NV(x)}(1+O(N^{-1})),
    \end{equation}
    where  $V$ is given by \eqref{eqn:def of potential}.  
    Direct calculation from the definitions and
    Lemma~\ref{lem:zast relation} gives
    \begin{equation}\label{eqn:Phix Vprime}
        2\Phi_x(a,z_\ast)=V'(a)
    \end{equation}
    and the prefactor cancellation
    \begin{equation*}
        -\frac{\mathsf{C}_N^{(\mathsf{a},\mathsf{b})}
        w(a;q)e^{2N\Phi(a,z_\ast)}}{N}
        \frac{2\Psi(a,z_\ast)^2\mathfrak{R}_z(a,z_\ast)}
        {z_\ast^2\Phi_{zzz}(a,z_\ast)}
        =1+O(N^{-1}).
    \end{equation*}
    Moreover, \eqref{eqn:Phix Vprime} and Taylor expansion imply
    \begin{equation*}
        \frac{w(\widehat{x};q)e^{2N\Phi(\widehat{x},z_\ast)}}
        {w(a;q)e^{2N\Phi(a,z_\ast)}}
        =1+O(N^{-1/3}),
    \end{equation*}
    uniformly for bounded real $u$.

    We now substitute Lemma~\ref{lem:soft pq asymptotic}
    into \eqref{eqn:soft kernel pq}.
    Since
    $\widehat{x}-\widehat{y}=-(u-v)/(\theta N)^{2/3}$,
    the preceding cancellation gives the Airy kernel from
    the divided-difference term.
    The convergence of the first derivatives of the normalised
    functions ensures uniformity also at $u=v$.
    The remaining product term is $O(N^{-1/3})$ after rescaling,
    by the same prefactor estimate and
    Lemma~\ref{lem:soft pq asymptotic}.
    Therefore,
    \begin{equation}
        \frac{1}{(\theta N)^{2/3}}
        K_N^{(\mathsf{a},\mathsf{b})}
        (\widehat{x},\widehat{y};q)
        =K_{\Ai}(u,v)+o(1)
    \end{equation}
    uniformly for $u,v\in\mathfrak{L}_M^{\textup{soft}}$,
    as required.
\end{proof}

To pass from local kernel convergence to the Tracy--Widom limit,
we apply Johansson's criterion \cite[Lemma~3.1]{Joh00}
to the counting-measure kernel
$\mathsf{K}_N^{(\mathsf{a},\mathsf{b})}$ defined in
\eqref{eqn:def of mathsfK}.
The following diagonal estimate provides the required tail control
in the void region; its proof is deferred to
Appendix~\ref{app:soft edge diagonal}.

\begin{prop}\label{prop:soft diagonal estimate}
    Under the assumptions of Proposition~\ref{prop:soft kernel uniform},
    there exists a constant $C >0$ such that
    \begin{equation}\label{eqn:soft diagonal estimate}
        0\leq\mathsf{K}_N^{(\mathsf{a},\mathsf{b})}(k,k;q)
        \leq
        \frac{C (\theta N)^{2/3}(1-q)q^k}
        {(1+(\theta N)^{2/3}(a-q^k))^2}
    \end{equation}
    for all sufficiently large $N$ and every
    $k\in\mathbb{Z}_{\geq0}$ with $q^k\leq a$.
\end{prop}

We now combine local kernel convergence with the diagonal estimate to prove Theorem~\ref{Thm_qJUE scaling limits}~(ii).

\begin{proof}[{Proof of Theorem \ref{Thm_qJUE scaling limits} (ii)}]
Recall the counting-measure kernel $\mathsf{K}_N^{(\mathsf{a},\mathsf{b})}$ from \eqref{eqn:def of mathsfK}. 
For $t\in\mathbb{R}$, define the tail index set
\begin{equation}\label{def of JN soft}
   J_N(t)=\Big\{k:(\theta N)^{2/3}(a-q^k)>t\Big\}. 
\end{equation}
Note that by definition, we have 
 \begin{equation*}
        \mathbb{P}\Big(
        x_{\min}\geq a-\frac{s}{(\theta N)^{2/3}}
        \Big)
        =\det 
        (\Id-\mathsf{K}_N^{(\mathsf{a},\mathsf{b})})_{\ell^2(J_N(s))}. 
    \end{equation*} 
Set 
    \begin{equation*}
        \rho_N:=\frac{N^{1/3}}{\lambda a\theta^{2/3}},
        \qquad
        k_N^{\textup{soft}}(u)
        :=\Big\lfloor-\frac{N}{\lambda}\log a+\rho_Nu\Big\rfloor.
    \end{equation*}
Then 
    \begin{equation*}
        (\theta N)^{2/3}(a-q^{k_{N}^{\textup{soft}}(u)})
        =u+O(N^{-\frac{1}{3}}),\qquad
        \rho_{N}(1-q)(\theta N)^{2/3}q^{\frac{k_{N}^{\textup{soft}}(u)+k_{N}^{\textup{soft}}(v)}{2}}
        =1+o(1).
    \end{equation*}
    locally uniformly for $u$, $v\in\mathbb{R}$. Therefore, it follows from \eqref{eqn:def of mathsfK} and Proposition~\ref{prop:soft kernel uniform} that 
    \begin{equation}\label{eqn:soft mathsfK uniform}
        \rho_N\mathsf{K}_N^{(\mathsf{a},\mathsf{b})}
        (k_N^{\textup{soft}}(u),k_N^{\textup{soft}}(v);q)
        =K_{\Ai}(u,v)+o(1)
    \end{equation}
    locally uniformly.
    Moreover, $\mathsf{K}_N^{(\mathsf{a},\mathsf{b})}$ is a
    projection kernel of rank $N$, so every matrix
    $(\mathsf{K}_N^{(\mathsf{a},\mathsf{b})}(k_i,k_j;q))_{i,j=1}^m$
    is positive semidefinite.

We next verify the diagonal conditions in
    \cite[Lemma~3.1]{Joh00}.
    Introduce the scaled lattice coordinates
    \begin{equation*}
        \xi_k:=(\theta N)^{2/3}(a-q^k),
        \qquad
        \xi_{k+1}-\xi_k=(\theta N)^{2/3}(1-q)q^k,
    \end{equation*}
    so that $J_N(t)=\{k:\xi_k>t\}$.
    Choose $N_0$ sufficiently large that
    Proposition~\ref{prop:soft diagonal estimate} applies and
    $(\theta N)^{2/3}(1-q)a<1$ for $N\geq N_0$.
    For $t\geq0$ and $k\in J_N(t)$, we then have
    $\xi_k>t$ and $\xi_{k+1}-\xi_k<1$, whence
    \begin{equation*}
        \mathsf{K}_N^{(\mathsf{a},\mathsf{b})}(k,k;q)
        \leq C 
        \frac{\xi_{k+1}-\xi_k}{(1+\xi_k)^2}
        \leq 2C 
        \int_{\xi_k}^{\xi_{k+1}}\frac{dr}{(1+r)^2}.
    \end{equation*}
    Summing over $k\in J_N(t)$ yields
    \begin{equation}\label{eqn:Johcondi2}
        \sup_{N\geq N_0}
        \sum_{k\in J_N(t)}
        \mathsf{K}_N^{(\mathsf{a},\mathsf{b})}(k,k;q)
        \leq\frac{2C }{1+t},
        \qquad t\geq0.
    \end{equation}
    In particular, the tails vanish uniformly as $t\to\infty$.

    For fixed $M>0$, the set
    $J_N(-M)\setminus J_N(0)$ contains $O(\rho_N)$ indices,
    and \eqref{eqn:soft mathsfK uniform} bounds the corresponding
    diagonal entries by $O(\rho_N^{-1})$.
    Combining this with \eqref{eqn:Johcondi2} at $t=0$ gives
    \begin{equation}\label{eqn:Johcondi1}
        \sup_{N\geq N_0}
        \sum_{k\in J_N(-M)}
        \mathsf{K}_N^{(\mathsf{a},\mathsf{b})}(k,k;q)<\infty.
    \end{equation}
    Any finitely many remaining values of $N$ are covered
    by the identity
    $\sum_{k\geq0}\mathsf{K}_N^{(\mathsf{a},\mathsf{b})}(k,k;q)=N$.

    To relate $J_N(t)$ to the linear cut-offs in Johansson's
    criterion, note that, for fixed $t$,
    \begin{equation*}
        -\frac{N}{\lambda}
        \log\Big(a-\frac{t}{(\theta N)^{2/3}}\Big)
        =-\frac{N}{\lambda}\log a+\rho_Nt+O(N^{-1/3}).
    \end{equation*}
    Thus the preceding trace bounds also give the required
    conditions for linear cut-offs; finitely many initial
    values of $N$ cause no difficulty since their diagonals
    are summable.
    Johansson's criterion, together with monotonicity of gap
    probabilities and continuity of $F_{\textup{TW}}$, therefore gives
    \begin{equation*}
    \begin{split}
        \lim_{N\to\infty}
        \det 
        (\Id-\mathsf{K}_N^{(\mathsf{a},\mathsf{b})})_{\ell^2(J_N(s))}
        &=\det(\Id-K_{\Ai})_{L^2(s,\infty)} =F_{\textup{TW}}(s),
    \end{split}
    \end{equation*}
    which completes the proof. 
\end{proof}

\appendix

\section{Spectral moments of the \texorpdfstring{$q$}{q}-JUE} \label{Appendix_spectral moment}

In this appendix, we present a complementary derivation of the explicit candidate for the equilibrium measure in Theorem~\ref{thm:limiting density}. While this derivation is not needed for the proof given in Section~\ref{Sec_eq msr}, it explains how the explicit form of the limiting density arises naturally. Following the moment method in \cite{BFO26,BJO25,BJM26}, we first derive an explicit formula for the spectral moments of the $q$-JUE and then analyse their large-$N$ asymptotics. The resulting limiting moments determine the Cauchy transform, from which the density $\rho^{(c_1,c_2)}$ can be recovered explicitly. Besides providing a systematic way to identify the candidate equilibrium measure, the explicit finite-$N$ formula for the spectral moments and its asymptotic analysis are also of independent interest.

The central object in the moment method is the \emph{spectral moment} defined by 
\begin{equation}\label{eqn:def of spectral moment}
    m_{N,p}^{(\mathsf{a},\mathsf{b})}
    =\mathbb{E}\Big[\sum_{j=1}^{N}x_{j}^{p}\Big]
\end{equation}
where the expectation is taken with respect to \eqref{eqn:jpdf of qJUE}. Using the orthogonal polynomial representation
of the correlation kernel \eqref{def of kernel OP sum}, the spectral
moment can be written as 
\begin{equation}\label{eqn:OP representation of spectral moment}
    m_{N,p}^{(\mathsf{a},\mathsf{b})}
    =\sum_{j=0}^{N-1}\frac{1}{h_{j}^{(\mathsf{a},\mathsf{b})}}\int_{0}^{1}x^{p}p_{j}^{(\mathsf{a},\mathsf{b})}(x;q)^{2}w^{(\mathsf{a},\mathsf{b})}(x;q)\,d_{q}x.
\end{equation}

Recall that the Legendre polynomials $\{\mathcal{P}_n(x)\}_{n\geq0}$ are characterised by the generating function  
\begin{equation}\label{eqn:Legendre generating}
    \frac{1}{\sqrt{1-2xt+t^{2}}}=\sum_{p=0}^{\infty}t^{p}\mathcal{P}_{p}(x),\qquad
    |t|<\begin{cases}
        \displaystyle\frac{1}{|x|+\sqrt{x^{2}-1}},&\textup{if }|x|>1, \smallskip\\
        1,&\textup{if }|x|\leq 1;
    \end{cases}
\end{equation}
see e.g. \cite[Eq. (9.8.68)]{KLS10}.
The following proposition gives an explicit formula for the spectral moments, together with their large-$N$ asymptotics.  

\begin{prop}[Spectral moments of the $q$-Jacobi unitary ensemble]\label{prop:spectral moment}
    Let $m_{N,p}^{(\mathsf{a},\mathsf{b})}$ be the spectral moment of the $q$-JUE defined in \eqref{eqn:def of spectral moment}.
    \begin{itemize}
    \item [\textup{(A)}] \textup{\textbf{(Closed formula)}} For any nonnegative integer $p$ and parameters $\mathsf{a},\mathsf{b}>-1$, we have
    \begin{align}
    \begin{split}
    \label{eqn:moment sum formula}
        m_{N,p}^{(\mathsf{a},\mathsf{b})}
        &=\sum_{j=0}^{N-1}\sum_{r=0}^{j\wedge p}q^{p j+r(r+\mathsf{a}-p)}\qbinom{p}{r}^{2}\frac{1-q^{2j+\mathsf{a}+\mathsf{b}+1}}{1-q^{2j-2r+p+\mathsf{a}+\mathsf{b}+1}}
        \\
        &\qquad \times \frac{(q^{j-r+\mathsf{b}+1};q)_{r}(q^{j+\mathsf{a}+1},q^{j+\mathsf{a}+\mathsf{b}+1};q)_{p-r}(q^{j-r+1};q)_{r}}{(q^{2j-r+\mathsf{a}+\mathsf{b}+1};q)_{p-r}^{2}(q^{2j-2r+p+\mathsf{a}+\mathsf{b}+2};q)_{r}^{2}}.
    \end{split}
    \end{align}  Here,
\[
\begin{bmatrix}p\\r\end{bmatrix}_q
:=\frac{(q;q)_p}{(q;q)_r(q;q)_{p-r}},
\qquad 0\le r\le p.
\]
    \item [\textup{(B)}] \textup{\textbf{(Large $N$-asymptotic)}} Let $\mathsf{a}$ and $\mathsf{b}$ be scaled as in
\eqref{eqn:scale of alphabet}, and let $q=e^{-\lambda/N}$
for a fixed $\lambda>0$.
Then, for each fixed nonnegative integer $p$,
\begin{equation}\label{asymptotic of spectral moment}
    m_{N,p}^{(\mathsf{a},\mathsf{b})}
    =N\mathcal{M}_{p}+O(1),
    \qquad N\to\infty.
\end{equation}
Here, for $c_1>0$,
\begin{equation}\label{eqn:moment asymptotic}
    \mathcal{M}_{p}
    :=\int_{0}^{1}
    \Big(\frac{\mathsf{s}^{t}(1-\mathsf{s}^{c_1})}
    {1-\mathsf{s}^{2t+c_1+c_2}}\Big)^{p}
    \mathcal{P}_{p}\Big(
    \frac{(1-\mathsf{s}^{t+c_1})(1-\mathsf{s}^{t+c_1+c_2})
    +\mathsf{s}^{c_1}(1-\mathsf{s}^{t})(1-\mathsf{s}^{t+c_2})}
    {(1-\mathsf{s}^{c_1})(1-\mathsf{s}^{2t+c_1+c_2})}
    \Big)\,dt,
\end{equation}
while for $c_1=0$, it is defined by continuous extension as
\begin{equation}\label{eqn:moment asymptotic c1 zero}
    \mathcal{M}_{p}
    :=\binom{2p}{p}\int_{0}^{1}
    \Big(
    \frac{\mathsf{s}^{t}(1-\mathsf{s}^{t})(1-\mathsf{s}^{t+c_2})}
    {(1-\mathsf{s}^{2t+c_2})^{2}}
    \Big)^{p}\,dt.
\end{equation}
    \end{itemize}
\end{prop}

The derivation of \eqref{eqn:moment sum formula} takes a different route from the proofs in \cite{BFO26,BJO25,BJM26}, where explicit moment formulas were obtained by combinatorial enumeration based on the Flajolet--Viennot theory.
Here, we instead exploit parameter shifts and connection coefficients of orthogonal polynomials, following a related idea used in \cite{MS11} for the classical JUE. 

\begin{proof}[Proof of Proposition \ref{prop:spectral moment} (A)]
    We rewrite \eqref{eqn:OP representation of spectral moment} as  
    \begin{equation}\label{eqn:def of mathfrak m}
       m^{(\mathsf{a},\mathsf{b})}_{N,p}
        =\sum_{j=0}^{N-1}\frac{1}{h^{^{(\mathsf{a},\mathsf{b})}}_{j}}\mathfrak{m}_{p,j}^{(\mathsf{a},\mathsf{b})},\qquad  \mathfrak{m}_{p,j}^{(\mathsf{a},\mathsf{b})}
        :=\int_{0}^{1}x^{p}p_{j}^{(\mathsf{a},\mathsf{b})}(x;q)^{2}w^{(\mathsf{a},\mathsf{b})}(x;q)\,d_{q}x. 
    \end{equation} 
By the definition \eqref{eqn:def of qjacobi weight} of the little $q$-Jacobi weight, we have
\begin{equation}\label{eqn:xpw}
    x^{p}w^{(\mathsf{a},\mathsf{b})}(x;q)=w^{(\mathsf{a}+p,\mathsf{b})}(x;q).
\end{equation}
Thus, $\mathfrak{m}_{p,j}^{(\mathsf{a},\mathsf{b})}$ is the squared $\ell^2(w^{(\mathsf{a}+p,\mathsf{b})}\,d_qx)$-norm of $p_j^{(\mathsf{a},\mathsf{b})}$.
To evaluate this norm, we expand the polynomial in the orthogonal basis $\{p_l^{(\mathsf{a}+p,\mathsf{b})}\}_{0\leq l\leq j}$.

For this purpose, recall the connection formula for the little $q$-Jacobi polynomials \cite[Corollary~3.2]{IO11}:
\begin{equation}\label{eqn:connection formula}
    p_{n}^{(\mathsf{a}_{1},\mathsf{b}_{1})}(x;q)
    =\sum_{k=0}^{n}c_{n,k}^{(\mathsf{a}_{1},\mathsf{b}_{1})\rightarrow (\mathsf{a}_{2},\mathsf{b}_{2})}p_{k}^{(\mathsf{a}_{2},\mathsf{b}_{2})}(x;q),
\end{equation} 
where
\begin{equation}
    c_{n,k}^{(\mathsf{a}_{1},\mathsf{b}_{1})\rightarrow (\mathsf{a}_{2},\mathsf{b}_{2})}
    =q^{k(k-n)}\qbinom{n}{k}\frac{(q^{\mathsf{a}_{2}+1},q^{n+\mathsf{a}_{1}+\mathsf{b}_{1}+1};q)_{k}}{(q^{\mathsf{a}_{1}+1},q^{k+\mathsf{a}_{2}+\mathsf{b}_{2}+1};q)_{k}}{}_{3}
    \phi_{2}\Bigg(\begin{matrix}
        q^{k-n},q^{n+k+\mathsf{a}_{1}+\mathsf{b}_{1}+1},q^{k+\mathsf{a}_{2}+1}\\q^{k+\mathsf{a}_{1}+1},q^{2k+\mathsf{a}_{2}+\mathsf{b}_{2}+2}
    \end{matrix};q,q\Bigg).
\end{equation}
When $\mathsf{b}_{1}=\mathsf{b}_{2}=\mathsf{b}$, the terminating ${}_3\phi_2$ series can be evaluated using the $q$-Pfaff--Saalsch\"utz summation formula \cite[Eq.~(1.11.9)]{KLS10},
\begin{equation}
    {}_{3}\phi_{2}\Bigg(\begin{matrix}q^{-n},a,b\\c,abc^{-1}q^{1-n}\end{matrix};q,q\Bigg)=\frac{(a^{-1}c,b^{-1}c;q)_{n}}{(c,a^{-1}b^{-1}c;q)_{n}},\quad n=0,1,2,\ldots.
\end{equation}
Applying this identity, we obtain
\begin{equation}\label{eqn:simple connection coef}
\begin{split}
    c_{n,k}^{(\mathsf{a}_{1},\mathsf{b})\rightarrow(\mathsf{a}_{2},\mathsf{b})}
    &=q^{k(k-n)}\qbinom{n}{k}\frac{(q^{\mathsf{a}_{2}+1},q^{n+\mathsf{a}_{1}+\mathsf{b}+1};q)_{k}}{(q^{\mathsf{a}_{1}+1},q^{k+\mathsf{a}_{2}+\mathsf{b}+1};q)_{k}}
    \frac{(q^{-(n+\mathsf{b})},q^{\mathsf{a}_{1}-\mathsf{a}_{2}};q)_{n-k}}{(q^{k+\mathsf{a}_{1}+1},q^{-(n+k+\mathsf{a}_{2}+\mathsf{b}+1)};q)_{n-k}}.
\end{split}
\end{equation}

We now specialise to $\mathsf{a}_1=\mathsf{a}$ and $\mathsf{a}_2=\mathsf{a}+p$.
Substituting \eqref{eqn:connection formula} into \eqref{eqn:def of mathfrak m} and using \eqref{eqn:xpw}, orthogonality with respect to the shifted weight gives
\begin{equation*}
    \frac{1}{h_{j}^{(\mathsf{a},\mathsf{b})}}\mathfrak{m}_{p,j}^{(\mathsf{a},\mathsf{b})}
    =\sum_{l=0}^{j}\frac{h_{l}^{(\mathsf{a}+p,\mathsf{b})}}{h_{j}^{(\mathsf{a},\mathsf{b})}}\Big(c_{j,l}^{(\mathsf{a},\mathsf{b})\rightarrow(\mathsf{a}+p,\mathsf{b})}\Big)^{2}.
\end{equation*}
Inserting \eqref{eqn:simple connection coef} and changing variables to $l=j-r$, we find
\begin{equation*}
    \frac{1}{h_{j}^{(\mathsf{a},\mathsf{b})}}\mathfrak{m}_{p,j}^{(\mathsf{a},\mathsf{b})}
    =\sum_{r=0}^{j\wedge p}\frac{h_{j-r}^{(\mathsf{a}+p,\mathsf{b})}}{h_{j}^{(\mathsf{a},\mathsf{b})}}q^{-2r(j-r)}\qbinom{j}{r}^{2}\bigg(\frac{(q^{p+\mathsf{a}+1},q^{j+\mathsf{a}+\mathsf{b}+1};q)_{j-r}}{(q^{j-r+p+\mathsf{a}+\mathsf{b}+1};q)_{j-r}}\frac{(q^{-(j+\mathsf{b})},q^{-p};q)_{r}}{(q^{-(2j-r+p+\mathsf{a}+\mathsf{b}+1)};q)_{r}}\frac{1}{(q^{\mathsf{a}+1};q)_{j}}\bigg)^{2}.
\end{equation*}
Here, the restriction $r\leq j$ comes from $l\geq0$, while $(q^{-p};q)_r=0$ for $r>p$.
For $0\leq r\leq j\wedge p$, applying the elementary identity
\begin{equation*}
    (q^{-n};q)_{k}=(-1)^{k}q^{-nk+\frac{k(k-1)}{2}}(q^{n-k+1};q)_{k},
    \qquad n\in\mathbb{R},\quad k\in\mathbb{Z}_{\geq0},
\end{equation*}
to the factors with negative exponents yields 
\begin{equation}
    \Big(c_{j,j-r}^{(\mathsf{a},\mathsf{b})\rightarrow(\mathsf{a}+p,\mathsf{b})}\Big)^{2}
    =q^{r(r+2\mathsf{a}+1)}\qbinom{j}{r}^{2}\bigg(\frac{(q^{p+\mathsf{a}+1},q^{j+\mathsf{a}+\mathsf{b}+1};q)_{j-r}(q^{j-r+\mathsf{b}+1},q^{p-r+1};q)_{r}}{(q^{j-r+p+\mathsf{a}+\mathsf{b}+1};q)_{j-r}(q^{2j-2r+p+\mathsf{a}+\mathsf{b}+2};q)_{r}(q^{\mathsf{a}+1};q)_{j}}\bigg)^{2}.
\end{equation}
On the other hand, the squared norms in \eqref{eqn:def of hn} give
\begin{equation*}
\begin{split}
    \frac{h_{j-r}^{(\mathsf{a}+p,\mathsf{b})}}{h_{j}^{(\mathsf{a},\mathsf{b})}}
    &=q^{pj-r(p+\mathsf{a}+1)}\frac{(q^{\mathsf{a}+1};q)_{p}}{(q^{\mathsf{a}+\mathsf{b}+1};q)_{p}}
    \frac{1-q^{2j+\mathsf{a}+\mathsf{b}+1}}{1-q^{2j-2r+p+\mathsf{a}+\mathsf{b}+1}} \frac{(q^{\mathsf{a}+1},q^{\mathsf{a}+\mathsf{b}+1};q)_j}{(q^{p+\mathsf{a}+1},q^{p+\mathsf{a}+\mathsf{b}+1};q)_{j-r}(q^{j-r+1},q^{j-r+\mathsf{b}+1};q)_{r}}.
\end{split}
\end{equation*}
Substituting these expressions, the powers of $q$ combine to
$q^{pj+r(r+\mathsf{a}-p)}$.
Simplifying the remaining factors using 
\begin{equation*}
  (u;q)_{m+n}=(u;q)_m(uq^m;q)_n, \qquad   \qbinom{j}{r}^{2}\frac{(q^{p-r+1};q)_{r}^{2}}{(q^{j-r+1};q)_{r}}
    =\qbinom{p}{r}^{2}(q^{j-r+1};q)_{r},
\end{equation*}
and then summing over $j=0,\ldots,N-1$, the desired identity \eqref{eqn:moment sum formula} follows. 
\end{proof}

Next, we turn to the asymptotic analysis of the spectral moment. 

\begin{proof}[Proof of Proposition \ref{prop:spectral moment} (B)] The case $p=0$ is immediate, so we assume that $p\geq1$ is fixed.
For $0\leq r\leq p$ and $j\geq r$, write
  \begin{equation}
    f_{r}^{(\mathsf{a},\mathsf{b})}(j)
    =q^{pj}\qbinom{p}{r}^{2}\frac{1-q^{2j+\mathsf{a}+\mathsf{b}+1}}{1-q^{2j-2r+p+\mathsf{a}+\mathsf{b}+1}}A_{p,j,r}^{(\mathsf{a},\mathsf{b})},
\end{equation}
where
\begin{equation}
    A_{p,j,r}^{(\mathsf{a},\mathsf{b})}
    :=\frac{(q^{j-r+\mathsf{b}+1};q)_{r}(q^{j+\mathsf{a}+1},q^{j+\mathsf{a}+\mathsf{b}+1};q)_{p-r}(q^{j-r+1};q)_{r}}{(q^{2j-r+\mathsf{a}+\mathsf{b}+1};q)_{p-r}^{2}(q^{2j-2r+p+\mathsf{a}+\mathsf{b}+2};q)_{r}^{2}}.
\end{equation} 
  Set $f_r^{(\mathsf{a},\mathsf{b})}(j)=0$ for $j<r$.
Then \eqref{eqn:moment sum formula} can be written as 
    \begin{equation}\label{eqn:moment sum}
        m_{N,p}^{(\mathsf{a},\mathsf{b})}
        =\sum_{r=0}^{p}q^{r(r+\mathsf{a}-p)}
        \sum_{j=0}^{N-1}f_{r}^{(\mathsf{a},\mathsf{b})}(j).
    \end{equation}

   Suppose first that $c_1+c_2>0$, and put $t=j/N$.
Since
\begin{equation*}
    1-q^{k}
    =\frac{k\lambda}{N}+O(N^{-2})
\end{equation*}
for each fixed $k>0$, we have
\begin{equation}\label{eqn:fr prefactor limit}
    q^{pj}\qbinom{p}{r}^{2}
    \frac{1-q^{2j+\mathsf{a}+\mathsf{b}+1}}
    {1-q^{2j-2r+p+\mathsf{a}+\mathsf{b}+1}}
    =\mathsf{s}^{pt}\binom{p}{r}^{2}+O(N^{-1}),
\end{equation}
uniformly for $r\leq j\leq N-1$. 
  
    Next, for a fixed $b\in\mathbb{R}$ and a fixed nonnegative integer
$c$, the expansion
\begin{equation*}
    (q^{aN+b};q)_{c}
    =\prod_{k=0}^{c-1}(1-\mathsf{s}^{a}q^{b+k})
    =(1-\mathsf{s}^{a})^{c}+O(N^{-1})
\end{equation*}
holds uniformly for $a$ in bounded subsets of $[0,\infty)$.
Applying this expansion to the finite products defining
$A_{p,j,r}^{(\mathsf{a},\mathsf{b})}$, and using
$1-\mathsf{s}^{2t+c_1+c_2}\geq1-\mathsf{s}^{c_1+c_2}>0$, gives
\begin{equation}\label{eqn:Apjr limit}
    A_{p,j,r}^{(\mathsf{a},\mathsf{b})}
    =(1-\mathsf{s}^{t})^{r}
    \frac{(1-\mathsf{s}^{t+c_{1}})^{p-r}(1-\mathsf{s}^{t+c_{2}})^{r}(1-\mathsf{s}^{t+c_{1}+c_{2}})^{p-r}}{(1-\mathsf{s}^{2t+c_{1}+c_{2}})^{2p}}+O(N^{-1}),
\end{equation}
uniformly for $r\leq j\leq N-1$.

Combining \eqref{eqn:fr prefactor limit} and \eqref{eqn:Apjr limit}, we obtain 
    \begin{equation*}
        f_{r}^{(\mathsf{a},\mathsf{b})}(j)
        =F_{r}^{(\mathsf{a},\mathsf{b})}(j/N)+O(N^{-1}),
    \end{equation*}
where 
\begin{align*}
      F_{r}^{(\mathsf{a},\mathsf{b})}(t)
        =\mathsf{s}^{pt}(1-\mathsf{s}^{t})^{r}\binom{p}{r}^{2}\frac{(1-\mathsf{s}^{t+c_{1}})^{p-r}(1-\mathsf{s}^{t+c_{2}})^{r}(1-\mathsf{s}^{t+c_{1}+c_{2}})^{p-r}}{(1-\mathsf{s}^{2t+c_{1}+c_{2}})^{2p}}. 
\end{align*}
Thus, summing the uniform error and applying the Riemann sum approximation yields
\begin{equation*}
    \sum_{j=0}^{N-1}f_{r}^{(\mathsf{a},\mathsf{b})}(j)
    =N\int_{0}^{1}F_{r}^{(\mathsf{a},\mathsf{b})}(t)\,dt+O(1),
\end{equation*}
where the finitely many omitted terms with $j<r$ contribute only $O(1)$. Since the scaling \eqref{eqn:scale of alphabet} also gives
\begin{equation*}
    q^{r(r+\mathsf{a}-p)}
    =\mathsf{s}^{c_1r}+O(N^{-1}),
\end{equation*}
substitution into \eqref{eqn:moment sum} yields
\begin{equation}\label{eqn:moment limiting sum}
    m_{N,p}^{(\mathsf{a},\mathsf{b})}
    =N\sum_{r=0}^{p}\mathsf{s}^{c_{1}r}
    \int_{0}^{1}F_{r}^{(\mathsf{a},\mathsf{b})}(t)\,dt+O(1).
\end{equation} 
    After simplifying the sum over $r$ using the generating function of the Legendre polynomial
    \begin{equation*}
        \sum_{r=0}^{p}\binom{p}{r}^{2}x^{r}=(1-x)^{p}\mathcal{P}_{p}\Big(\frac{1+x}{1-x}\Big),
    \end{equation*}
    the desired asymptotic behaviour \eqref{asymptotic of spectral moment} follows. 
    
For $c_1=c_2=0$, put $a=d_1$ and $b=d_2$, and denote
the $j$th summand in \eqref{eqn:def of spectral moment} by $T_{p,j}$.
Expanding the finite products in \eqref{eqn:moment sum formula} to first order and using
\[
\sum_{r=0}^p \binom{p}{r}^2=\binom{2p}{p},
\qquad
\sum_{r=0}^p r\binom{p}{r}^2
=\frac p2\binom{2p}{p},
\]
gives
\[
T_{p,j}
=\binom{2p}{p}
\biggl(\frac{q^j}{(1+q^j)^2}\biggr)^p
+O(N^{-1}+j^{-2}),
\]
uniformly for $J\le j<N$, where $J$ is a sufficiently large
constant independent of $N$.
The finitely many remaining terms are bounded by $1$.
Summing the error and applying the Riemann sum approximation proves \eqref{asymptotic of spectral moment}.
\end{proof}

We now explain how the limiting moments in Proposition~\ref{prop:spectral moment} lead to the explicit candidate equilibrium measure. More precisely, we derive the logarithmic expression for the Cauchy transform in Lemma~\ref{lem:Cauchy transform} directly from the limiting moments.
The boundary-value analysis in Section~\ref{Sec_eq msr} then recovers the density, without the need to repeat the inversion calculation. Throughout, we use the notation $B$ and $R$ introduced in \eqref{eqn:EL discriminant} and \eqref{def of B in terms of s}, respectively.
 
Fix a real number $z>1$. Since the particles lie in $[0,1]$,
Proposition~\ref{prop:spectral moment} \textup{(B)} gives
$0\leq\mathcal{M}_p\leq1$, so the series
$\sum_{p=0}^{\infty}\mathcal{M}_p/z^{p+1}$ converges absolutely.
We first assume that $c_1>0$; the resulting expressions extend to
$c_1=0$ by continuity.
Substituting \eqref{eqn:moment asymptotic}, we obtain
\begin{align*}
  \sum_{p=0}^{\infty}\frac{\mathcal{M}_p}{z^{p+1}}
    &=
    \frac{1}{z}\sum_{p=0}^{\infty}
    \int_0^1
    \Big(
        \frac{\mathsf{s}^t(1-\mathsf{s}^{c_1})}
        {z(1-\mathsf{s}^{2t+c_1+c_2})}
    \Big)^p 
    \mathcal{P}_p\bigg(
        \frac{
            (1-\mathsf{s}^{t+c_1})(1-\mathsf{s}^{t+c_1+c_2})
            +\mathsf{s}^{c_1}(1-\mathsf{s}^t)(1-\mathsf{s}^{t+c_2})
        }{
            (1-\mathsf{s}^{c_1})(1-\mathsf{s}^{2t+c_1+c_2})
        }
    \bigg)\,dt\\
    &=
    \frac{1}{\lambda z}\sum_{p=0}^{\infty}
    \int_{\mathsf{s}}^1
    \Big(
        \frac{u(1-\mathsf{s}^{c_1})}
        {z(1-\mathsf{s}^{c_1+c_2}u^2)}
    \Big)^p 
    \mathcal{P}_p\bigg(
        \frac{
            (1-\mathsf{s}^{c_1}u)(1-\mathsf{s}^{c_1+c_2}u)
            +\mathsf{s}^{c_1}(1-u)(1-\mathsf{s}^{c_2}u)
        }{
            (1-\mathsf{s}^{c_1})(1-\mathsf{s}^{c_1+c_2}u^2)
        }
    \bigg)\frac{du}{u}, 
\end{align*} 
where we used the change of variables $u=\mathsf{s}^t$ in the last
equality.
Interchanging the sum and integral, which is justified by the
non-negativity of the summands, and applying the Legendre generating
function \eqref{eqn:Legendre generating} together with the change of variables $ v=\frac{1}{u}+\mathsf{s}^{c_1+c_2}u$, we obtain 
\begin{equation}\label{eqn:moment Cauchy integral}
\begin{split}
   \sum_{p=0}^{\infty}\frac{\mathcal{M}_p}{z^{p+1}}
    &=\frac{1}{\lambda}\int_{1+\mathsf{s}^{c_1+c_2}}^B\frac{dv}{\sqrt{(zv-1-\mathsf{s}^{c_1})^2-4\mathsf{s}^{c_1}(1-z)(1-\mathsf{s}^{c_2}z)}}.
\end{split}
\end{equation}  
Evaluating this elementary integral gives 
\begin{equation}\label{eqn:moment Cauchy identification}
\begin{split}
    \sum_{p=0}^{\infty}\frac{\mathcal{M}_p}{z^{p+1}}
    &=
    \frac{1}{\lambda z}
    \log\frac{Bz-1-\mathsf{s}^{c_1}+R(z)}{2(z-1)} =g(z),
    \qquad z>1,
\end{split}
\end{equation}
where the last equality follows from Lemma~\ref{lem:Cauchy transform}.

This calculation derives the logarithmic Cauchy transform \eqref{eqn:EL Cauchy transform} directly from the limiting moments. The Sokhotski--Plemelj inversion carried out in the proof of Lemma~\ref{lem:Cauchy transform} then recovers the density $\rho^{(c_1,c_2)}$ in \eqref{def of limiting law}, explaining how the explicit candidate equilibrium measure is obtained.

\section{Diagonal estimates at the soft edge}
\label{app:soft edge diagonal}

In this appendix, we prove the diagonal estimate, Proposition~\ref{prop:soft diagonal estimate}. 
We retain the notation and assumptions of
Proposition~\ref{prop:soft kernel uniform}, and let
$\Gamma_{\textup{s}}$ be the contour constructed in
Lemma~\ref{lem:soft contour}.
We begin with a uniform bound for the weighted integrand
in the void region.
 
Recall that $\mathfrak{F}_n$ is given by \eqref{eqn:def of mathfrakF}. 

\begin{lem}\label{lem:wf bound}
    Let $q=e^{-\lambda/N}$ for a fixed $\lambda>0$. There exist constants $c_{\textup{v}}$, $C_{\textup{v}}>0$ such that, for a sufficiently large $N$,
    \begin{equation} 
        \sqrt{w^{(\mathsf{a},\mathsf{b})}(x;q)}|\mathfrak{F}_{N}(x,z)|
        \leq C_{\textup{v}}\sqrt{w^{(\mathsf{a},\mathsf{b})}(a;q)}|\mathfrak{F}_{N}(a,z)|e^{-c_{\textup{v}}N(a-x)(|z|-z_{\ast})},
    \end{equation}
    uniformly for $x\in(0,a]$ and $z\in\Gamma_{\textup{s}}$. Here, $z_*$ is given by \eqref{def of zstar}. 
\end{lem}

\begin{proof}
    Recall the phase function $\Phi$ in
    \eqref{def of Phi} and the external potential $V$ in \eqref{eqn:def of potential}.  
    Define 
    \begin{equation*}
    \begin{split}
        g(x)
        :=\lambda x\Big(2\Phi_x(x,z_{\ast})-V'(x) \Big)  =2\log\frac{z_{\ast}+\mathsf{s}^{c_1+c_2+1}x}{z_{\ast}+x}
        +\log\frac{1-x}{1-\mathsf{s}^{c_2}x}
        +c_1\lambda.
    \end{split}
    \end{equation*}
    This function is strictly decreasing on $(0,a]$, and
    $g(a)=0$ by \eqref{eqn:Phix Vprime}.
    Thus $g(x)\geq0$ and $g(a/2)>0$.

    Recall that the contour $\Gamma_{\textup{s}}$ constructed in
    Lemma~\ref{lem:soft contour} satisfies $|z|\geq z_{\ast}$.
    Since
    \begin{equation*}
    \Phi_{xz}(x,r)
        =\frac{1-\mathsf{s}^{c_1+c_2+1}}
        {\lambda(r+x)(r+\mathsf{s}^{c_1+c_2+1}x)},
    \end{equation*}
    there is a constant $c_{\Phi}>0$ such that
    \begin{equation}\label{eqn:2Phix-V' bound}
    \begin{split}
        2\Phi_x(x,r)-V'(x)
        &=\frac{g(x)}{\lambda x}
        +2\int_{z_{\ast}}^r\Phi_{xz}(x,t)\,dt \geq\frac{g(x)}{\lambda x}
        +c_{\Phi}(r-z_{\ast})
    \end{split}
    \end{equation}
    uniformly for $x\in(0,a]$ and
    $r\in[z_{\ast},\max_{z\in\Gamma_{\textup{s}}}|z|]$.
    Moreover, the saddle identity \eqref{eqn:zast2} implies
    $z_{\ast}^2>\mathsf{s}^{c_1+c_2+1}a^2$.
    Hence, for $z\in\Gamma_{\textup{s}}$ and $r=|z|$,
    \begin{equation}\label{eqn:lambdaxPhix bound}
        \big|\exp[\lambda x\Phi_x(x,z)]\big|
        =\Big|\frac{z+\mathsf{s}^{c_1+c_2+1}x}{z+x}\Big|
        \geq
        \frac{r+\mathsf{s}^{c_1+c_2+1}x}{r+x}.
    \end{equation}
    Combining these inequalities and integrating from $x$ to $a$
    gives
    \begin{equation}\label{eqn:Rephi-rephi bound}
    \begin{split}
        &2\Re\Phi(x,z)-2\Re\Phi(a,z)-(V(x)-V(a)) \leq
        -c_{\Phi}(a-x)(|z|-z_{\ast})
        -\frac{1}{\lambda}\int_x^a\frac{g(t)}{t}\,dt.
    \end{split}
    \end{equation}

   Note that by \eqref{eqn:def of mathfrakF}, cancelling the factors independent of $x$ gives
    \begin{equation*}
        \frac{\mathfrak{F}_N(x,z)}{\mathfrak{F}_N(a,z)}
        =
        \frac{(-q^{N+\mathsf{a}+\mathsf{b}+2}x/z,-qa/z;q)_{\infty}}
        {(-qx/z,-q^{N+\mathsf{a}+\mathsf{b}+2}a/z;q)_{\infty}},
    \end{equation*}
    which also defines the quotient at common zeros. 
    Since $\Gamma_{\textup{s}}$ strictly encloses $[-a,0]$,
    the remaining factors have arguments in a fixed compact
    subset of $\mathbb{C}\setminus[1,\infty)$.
    Then it follows from \eqref{eqn:def of qjacobi weight} and 
    \eqref{eqn:qPochhammer asymptotic} that 
    \begin{equation}\label{eqn:log ratio expansion}
    \begin{split}
        &\log\Bigg(
        \sqrt{\frac{w^{(\mathsf{a},\mathsf{b})}(x;q)}
        {w^{(\mathsf{a},\mathsf{b})}(a;q)}}
        \Big|\frac{\mathfrak{F}_N(x,z)}
        {\mathfrak{F}_N(a,z)}\Big|
        \Bigg) =
        N\Big(
        \Re\Phi(x,z)-\Re\Phi(a,z)
        -\frac{V(x)-V(a)}{2}
        \Big)
        +\frac{d_1}{2}\log\frac{x}{a}+O(1),
    \end{split}
    \end{equation}
    uniformly for $x\in(0,a]$ and $z\in\Gamma_{\textup{s}}$.

    Finally, since $g$ is decreasing and non-negative,
    \begin{equation*}
        \int_x^a\frac{g(t)}{t}\,dt
        \geq g(a/2)\log\frac{a}{2x},
        \qquad 0<x<a/2.
    \end{equation*}
    It follows that
    \begin{equation*}
        -\frac{N}{2\lambda}\int_x^a\frac{g(t)}{t}\,dt
        +\frac{d_1}{2}\log\frac{x}{a}
    \end{equation*}
    is bounded above uniformly for $x\in(0,a]$ and all sufficiently
    large $N$, including when $d_1<0$.
    Applying \eqref{eqn:Rephi-rephi bound} to
    \eqref{eqn:log ratio expansion} and exponentiating proves
    the assertion with $c_{\textup{v}}=c_{\Phi}/2$.
\end{proof}

We next estimate $\mathfrak{p}_N$ and $\mathfrak{q}_N$, defined in
\eqref{eqn:mathfrak p integral} and
\eqref{eqn:mathfrak q integral}, together with their derivatives.
To remove the common contribution to the derivatives in the
Christoffel--Darboux Wronskian, introduce
\begin{equation}\label{def of mathfrak l}
    \mathfrak{l}_N(x,z)
    :=\frac{\partial}{\partial x}
    \log\frac{(-q^{N+\mathsf{a}+\mathsf{b}+2}x/z;q)_{\infty}}
    {(-qx/z;q)_{\infty}}.
\end{equation}
By \eqref{eqn:def of mathfrakF}, this satisfies
$\partial_x\mathfrak{F}_N(x,z)
=\mathfrak{l}_N(x,z)\mathfrak{F}_N(x,z)$.

\begin{lem}\label{lem:soft pq bound}
    For $x\in(0,a]$, set $u=(\theta N)^{2/3}(a-x)$.
    There exists a constant $C>0$,
    independent of $N$ and $x$, such that, for all sufficiently
    large $N$,
    \begin{align}
        \sqrt{w^{(\mathsf{a},\mathsf{b})}(x;q)}
        |\mathfrak{p}_N(x)|
        &\leq
        C
        \sqrt{w^{(\mathsf{a},\mathsf{b})}(a;q)}
        e^{N\Phi(a,z_{\ast})}
        \frac{1}{N^{1/3}(1+u)},
        \label{eqn:soft p bound}\\[3pt]
        \sqrt{w^{(\mathsf{a},\mathsf{b})}(x;q)}
        |\mathfrak{q}_N(x)|
        &\leq
       C
        \sqrt{w^{(\mathsf{a},\mathsf{b})}(a;q)}
        e^{N\Phi(a,z_{\ast})}
        \frac{1}{N^{2/3}(1+u)^2},
        \label{eqn:soft q bound}\\[3pt]
        \sqrt{w^{(\mathsf{a},\mathsf{b})}(x;q)}
        |\mathfrak{p}_N'(x)
        -\mathfrak{l}_N(x,z_{\ast})\mathfrak{p}_N(x)|
        &\leq
        C
        \sqrt{w^{(\mathsf{a},\mathsf{b})}(a;q)}
        e^{N\Phi(a,z_{\ast})}
        \frac{N^{1/3}}{(1+u)^2},
        \label{eqn:soft p' bound}\\[3pt]
        \sqrt{w^{(\mathsf{a},\mathsf{b})}(x;q)}
        |\mathfrak{q}_N'(x)
        -\mathfrak{l}_N(x,z_{\ast})\mathfrak{q}_N(x)|
        &\leq
        C
        \sqrt{w^{(\mathsf{a},\mathsf{b})}(a;q)}
        e^{N\Phi(a,z_{\ast})}
        \frac{1}{(1+u)^3},
        \label{eqn:soft q' bound}
    \end{align}
    uniformly for $x\in(0,a]$.
\end{lem}

\begin{proof}
    We first prove that, for $j=0,1,2$,
    \begin{equation}\label{eqn:soft polynomial estimate}
    \begin{split}
        &\int_{\Gamma_{\textup{s}}}
        |z-z_{\ast}|^j
        \sqrt{w^{(\mathsf{a},\mathsf{b})}(x;q)}
        \Big|\frac{\mathfrak{F}_N(x,z)}{z}\Big|\,|dz| \leq
        C_j\sqrt{w^{(\mathsf{a},\mathsf{b})}(a;q)}
        e^{N\Phi(a,z_{\ast})}
        N^{-(j+1)/3}(1+u)^{-(j+1)},
    \end{split}
    \end{equation}
    uniformly for $x\in(0,a]$ and all sufficiently large $N$,
    where $u=(\theta N)^{2/3}(a-x)$.

    By Lemma~\ref{lem:soft contour}, the contour near $z_{\ast}$
    consists of the two segments
    $z=z_{\ast}+te^{\pm i\pi/3}$, $0\leq t\leq\varepsilon$.
    The saddle identities in Lemma~\ref{lem:zast relation} give,
    for sufficiently small $\varepsilon$,
    \begin{equation*}
    \begin{split}
        \Re\Phi(a,z)-\Phi(a,z_{\ast})
        \leq-\frac{\Phi_{zzz}(a,z_{\ast})}{12}t^3,\qquad
        |z|-z_{\ast} \geq t/2.
    \end{split}
    \end{equation*}
    Combining the local expansion \eqref{eqn:mathfrak F asymp}
    with Lemma~\ref{lem:wf bound}, we therefore obtain
    \begin{equation*}
        \sqrt{w^{(\mathsf{a},\mathsf{b})}(x;q)}
        |\mathfrak{F}_N(x,z)|
        \leq
        C\sqrt{w^{(\mathsf{a},\mathsf{b})}(a;q)}
        e^{N\Phi(a,z_{\ast})}
        e^{-cN(t^3+(a-x)t)}
    \end{equation*}
    for some constants $C,c>0$.
    Since $|z|\geq z_{\ast}>0$, the substitution
    $r=N^{1/3}t$ bounds the contribution of these segments by
    \begin{equation*}
    \begin{split}
        C\sqrt{w^{(\mathsf{a},\mathsf{b})}(a;q)}
        e^{N\Phi(a,z_{\ast})}
        N^{-(j+1)/3}
        \int_0^\infty r^j e^{-c(r^3+ur)}\,dr \leq
        C_j\sqrt{w^{(\mathsf{a},\mathsf{b})}(a;q)}
        e^{N\Phi(a,z_{\ast})}
        N^{-(j+1)/3}(1+u)^{-(j+1)}.
    \end{split}
    \end{equation*}
    Here, the constants may be adjusted, and the last inequality
    follows by considering $0\leq u\leq1$ and $u\geq1$ separately.

    On the remaining part of the contour,
    \eqref{eqn:soft out bound}, applied at $x=a$, and
    Lemma~\ref{lem:wf bound} bound the contribution by
    \begin{equation*}
        C\sqrt{w^{(\mathsf{a},\mathsf{b})}(a;q)}
        e^{N\Phi(a,z_{\ast})}e^{-cN}.
    \end{equation*}
    Since $0\leq u\leq a(\theta N)^{2/3}$, this exponential
    remainder is absorbed by the right-hand side of
    \eqref{eqn:soft polynomial estimate}.
    This proves the claimed integral estimate.

    We now apply it to the contour representations
    \eqref{eqn:mathfrak p integral} and
    \eqref{eqn:mathfrak q integral}.
    Recall $\mathfrak{l}_N$ from \eqref{def of mathfrak l} and
    the multiplier $\mathfrak{R}_N(x,z)$ in \eqref{def of mathfrak R}.  
    The relevant denominators are uniformly separated from zero
    for $x\in[0,a]$ and $z\in\Gamma_{\textup{s}}$.
    Taking differences at $z$ and $z_{\ast}$ in these definitions,
    and using $\sum_{h\geq0}q^h=(1-q)^{-1}=O(N)$
    for the logarithmic derivative, gives
    $$
     |\mathfrak{l}_N(x,z)-\mathfrak{l}_N(x,z_{\ast})|
        \leq CN|z-z_{\ast}|,
    $$
    and 
    $$
  |\mathfrak{R}_N(x,z)-\mathfrak{R}_N(x,z_{\ast})|
         \leq C|z-z_{\ast}|,\qquad
        |\partial_x\mathfrak{R}_N(x,z)
        -\partial_x\mathfrak{R}_N(x,z_{\ast})|
         \leq C|z-z_{\ast}|.
    $$ 
    Consequently, the cases $j=0$ and $j=1$ of
    \eqref{eqn:soft polynomial estimate} yield
    \eqref{eqn:soft p bound} and \eqref{eqn:soft q bound},
    respectively.

    Differentiating the fixed-contour representations and using
    $\partial_x\mathfrak{F}_N=\mathfrak{l}_N\mathfrak{F}_N$,
    we obtain
    \begin{equation*}
    \begin{split}
        &\mathfrak{p}_N'(x)
        -\mathfrak{l}_N(x,z_{\ast})\mathfrak{p}_N(x) =
        \oint_{\Gamma_{\textup{s}}}
        \mathfrak{F}_N(x,z)
        \big(\mathfrak{l}_N(x,z)-\mathfrak{l}_N(x,z_{\ast})\big)
        \frac{dz}{2\pi iz}.
    \end{split}
    \end{equation*}
    Thus the case $j=1$ gives \eqref{eqn:soft p' bound}.
    Similarly,
    \begin{equation*}
    \begin{split}
        \mathfrak{q}_N'(x)
        -\mathfrak{l}_N(x,z_{\ast})\mathfrak{q}_N(x) 
        & =
        \oint_{\Gamma_{\textup{s}}}
        \mathfrak{F}_N(x,z)
        \big(\mathfrak{l}_N(x,z)-\mathfrak{l}_N(x,z_{\ast})\big)
        \big(\mathfrak{R}_N(x,z)-\mathfrak{R}_N(x,z_{\ast})\big)
        \frac{dz}{2\pi iz}\\
        &\qquad+
        \oint_{\Gamma_{\textup{s}}}
        \mathfrak{F}_N(x,z)
        \big(\partial_x\mathfrak{R}_N(x,z)
        -\partial_x\mathfrak{R}_N(x,z_{\ast})\big)
        \frac{dz}{2\pi iz}.
    \end{split}
    \end{equation*}
    Applying \eqref{eqn:soft polynomial estimate}
    with $j=2$ and $j=1$ gives
    \begin{equation*}
    \begin{split}
         \sqrt{w^{(\mathsf{a},\mathsf{b})}(x;q)}
        |\mathfrak{q}_N'(x)
        -\mathfrak{l}_N(x,z_{\ast})\mathfrak{q}_N(x)| \leq
        C\sqrt{w^{(\mathsf{a},\mathsf{b})}(a;q)}
        e^{N\Phi(a,z_{\ast})}
        \big((1+u)^{-3}+N^{-2/3}(1+u)^{-2}\big).
    \end{split}
    \end{equation*}
    Finally, $N^{-2/3}(1+u)\leq1+a\theta^{2/3}$
    allows the second term to be absorbed into the first,
    proving \eqref{eqn:soft q' bound}.
\end{proof}

We now complete the proof of Proposition~\ref{prop:soft diagonal estimate}.

\begin{proof}[Proof of Proposition~\ref{prop:soft diagonal estimate}]
    Fix $k\in\mathbb{Z}_{\geq0}$ with $q^k\leq a$, and set
    $x=q^k$ and $u=(\theta N)^{2/3}(a-x)\geq0$.
    Taking the diagonal limit in \eqref{eqn:soft kernel pq} gives
    \begin{equation*}
    \begin{split}
        \frac{K_N^{(\mathsf{a},\mathsf{b})}(x,x;q)}
        {\mathsf{C}_N^{(\mathsf{a},\mathsf{b})}
        w^{(\mathsf{a},\mathsf{b})}(x;q)}
        &=
        \mathfrak{p}_N(x)\mathfrak{q}_N'(x)
        -\mathfrak{p}_N'(x)\mathfrak{q}_N(x) +
        \partial_x\mathfrak{R}_N(x,z_{\ast})
        \mathfrak{p}_N(x)^2.
    \end{split}
    \end{equation*}
    The common term involving $\mathfrak{l}_N$, defined in
    \eqref{def of mathfrak l}, cancels in the Wronskian:
      \begin{equation*}
        \mathfrak{p}_{N}(x)\mathfrak{q}_{N}'(x)-\mathfrak{p}_{N}'(x)\mathfrak{q}_{N}(x)
        =\mathfrak{p}_{N}(x)(\mathfrak{q}_{N}'(x)-\mathfrak{l}_{N}(x,z_{\ast})\mathfrak{q}_{N}(x))-(\mathfrak{p}_{N}'(x)-\mathfrak{l}_{N}(x,z_{\ast})\mathfrak{p}_{N}(x))\mathfrak{q}_{N}(x).
    \end{equation*}
    Also,
    \begin{equation*}
        \partial_x\mathfrak{R}_N(x,z_{\ast})
        =\frac{q^{N+\mathsf{a}+\mathsf{b}+1}}
        {z_{\ast}(z_{\ast}+q^{N+\mathsf{a}})}
        =O(1),
    \end{equation*}
    while the prefactor calculation in the proof of
    Proposition~\ref{prop:soft kernel uniform} gives
    \begin{equation*}
        \mathsf{C}_N^{(\mathsf{a},\mathsf{b})}
        w^{(\mathsf{a},\mathsf{b})}(a;q)
        e^{2N\Phi(a,z_{\ast})}
        \leq CN.
    \end{equation*}
    Applying Lemma~\ref{lem:soft pq bound}, we therefore obtain
    \begin{equation*}
    \begin{split}
        0
        \leq\frac{K_N^{(\mathsf{a},\mathsf{b})}(x,x;q)}
        {(\theta N)^{2/3}}
        &\leq
        C\big((1+u)^{-4}+N^{-1/3}(1+u)^{-2}\big) \le 2 C(1+u)^{-2},
    \end{split}
    \end{equation*}
    uniformly for $x=q^k\leq a$ and all sufficiently large $N$.
    Here, non-negativity follows from the sum-of-squares
    representation of the correlation kernel. Now the desired estimate \eqref{eqn:soft diagonal estimate} follows from \eqref{eqn:def of mathsfK}.
\end{proof}

\subsection*{Acknowledgement} 
SSB and YGJ were supported by National Research Foundation of Korea grants (RS-2025-00516909, RS-2026-25607520, and RS-2026-25518141). NS gratefully acknowledges support from the Royal Society under grant URF\textbackslash R\textbackslash 231028.
We thank Peter J. Forrester and Kurt Johansson for their interest and helpful comments.   

\subsection*{AI use disclosure} 
The authors declare that generative AI tools (GPT-6 Astra) were used for language polishing and proof checking of the manuscript. All mathematical arguments and proofs were independently developed and written by the authors.


\end{document}